\documentclass[12pt]{article}
\usepackage{graphics,graphicx}
\usepackage{amssymb,amsmath,amsopn,amsfonts}
\usepackage[dvips]{color}
\usepackage{colordvi,multicol}
\usepackage{epsf}
\newfam\msbfam
\font\tenmsb=msbm10 \textfont\msbfam=\tenmsb \font\sevenmsb=msbm7
\scriptfont\msbfam=\sevenmsb \font\fivemsb=msbm5
\scriptscriptfont\msbfam=\fivemsb

\def\th#1{\vspace{1mm}\noindent{\bf #1}\quad}

\def\proof{\vspace{1mm}\noindent{\it Proof}\quad}

\numberwithin{equation}{section}

\def\bc{\begin{center}}
\def\ec{\end{center}}
\def\no{\noindent}
\def\hang{\hangindent\parindent}
\def\textindent#1{\indent\llap{\qquad #1\ \ \enspace}\ignorespaces}
\def\ref{\par\hang\textindent}

\begin{document}

\title{ {\bf Lattice approximation to the dynamical $\Phi_3^4$ model
\thanks{Research supported in part  by  NSFC (No.11301026, No.11401019) and DFG through IRTG 1132 and CRC
701, Key Lab of Random Complex Structures and Data Science, Chinese Academy of Sciences (Grant No.
2008DP173182)
}\\} }
\author{  {\bf Rongchan Zhu}$^{\mbox{a}}$, {\bf Xiangchan Zhu}$^{\mbox{b},}$\thanks{Corresponding author}
\date{}
\thanks{E-mail address:
zhurongchan@126.com(R. C. Zhu), zhuxiangchan@126.com(X. C. Zhu)}\\\\
$^{\mbox{a}}$Department of Mathematics, Beijing Institute of Technology, Beijing 100081,  China\\
$^{\mbox{b}}$School of Science, Beijing Jiaotong University, Beijing 100044, China  }

\maketitle

\noindent {\bf Abstract}

 We study  the lattice approximations to the dynamical $\Phi^4_3$ model by  paracontrolled distributions proposed in [GIP13]. We prove that the solutions to
the lattice systems
converge  to the solution to the dynamical $\Phi_3^4$ model in probability, locally in time.  The dynamical $\Phi_3^4$ model is not well defined in the classical sense.   Renormalisation has to be performed in order to define the non-linear term. Formally, this renormalisation corresponds to adding an infinite mass term to the equation which leads to adding  a drift term in the lattice systems.

\vspace{1mm}
\no{\footnotesize{\bf 2000 Mathematics Subject Classification AMS}:\hspace{2mm} 60H15, 82C28}
 \vspace{2mm}

\no{\footnotesize{\bf Keywords}:\hspace{2mm}  $\Phi_3^4$ model, regularity structure, paracontrolled distribution, space-time white noise, renormalisation}

\section{Introduction}

Recall that the usual continuum Euclidean $\Phi^4_d$-quantum field theory is heuristically described by the following
probability measure:
$$N^{-1}\Pi_{x\in\mathbb{T}^d}d\phi(x)\exp\bigg(-\int_{\mathbb{T}^d}(|\nabla \phi(x)|^2+\phi^2(x)+\phi^4(x))dx\bigg),\eqno(1.1)$$
where $N$
is the normalization constant and $\phi$ is the real-valued field and $\mathbb{T}^d$ is the $d$-dimensional torus. There have been many approaches to the problem of
giving a meaning to the above heuristic measure for $d=2$ and $d=3$ (see Refs. [GRS75], [GJ87] and references
therein).
In [PW81]  Parisi and Wu proposed a program for Euclidean quantum field
theory of getting Gibbs states of classical
statistical mechanics as limiting distributions of stochastic processes,  especially as solutions to non-linear stochastic differential
equations. Then one can use the stochastic differential equations to study the properties of the Gibbs states. This
procedure is called stochastic field quantization (see [JLM85]). The $\Phi_d^4$ model is the simplest non-trivial Euclidean quantum field (see [GJ87] and the reference therein). The issue of the stochastic quantization of the $\Phi^4_d$ model is to solve the following equation:
$$\aligned d \Phi= &( \Delta \Phi-\Phi^3)dt +dW(t) \quad\Phi(0)=\Phi_0.\endaligned\eqno(1.2)$$
where $W$ is a cylindrical Wiener process on $L^2(\mathbb{T}^d)$. The solution $\Phi$ is also called  dynamical $\Phi^4_d$ model.

In
two spatial dimensions, the dynamical $\Phi_2^4$ model was previously treated in [AR91], [DD03] and [MW15]. In three
spatial dimensions this equation (1.1) is ill-posed and the main difficulty in this case is that W
and hence the solutions are so singular that the non-linear term is not well-defined in the classical sense.
It was a long-standing open problem to give a meaning to equation (1.2) in the three dimensional
case. A breakthrough result was achieved recently by Martin Hairer in [Hai14], where he introduced
a theory of regularity structures and gave a meaning to  equation (1.2) successfully. Also by using the paracontrolled distributions proposed by Gubinelli, Imkeller and Perkowski
in [GIP13] existence and uniqueness of local solutions to (1.2) has been obtained in [CC13].
Recently, these two approaches have been successful in giving a meaning to a lot of ill-posed
stochastic PDEs like the Kardar-Parisi-Zhang (KPZ) equation ([KPZ86], [BG97], [Hai13]), the
dynamical $\Phi^4_3$
model ([Hai14], [CC13]), the Navier-Stokes equation driven by space-time white
noise ([ZZ14a], [ZZ14b]), the dynamical sine-Gordon equation ([HS14]) and so on (see [HP14]
for more interesting examples). From a philosophical perspective, the theory of regularity
structures and the paracontrolled distributions are inspired by the theory of controlled rough
paths [Lyo98], [Gub04]. The main difference is that the regularity structure theory considers
the problem locally, while the paracontrolled distribution method is a global approach using
Fourier analysis. In [Kup14] the author also use renormalization group
techniques to make sense of the dynamical  $\Phi^4_3$ model.
.

The lattice approximation is an important tool in constructing and studying the continuum $\Phi^4_3$ field (see [P75, P77, ABZ04]).
 It also preserves Osterwalder-Schrader positivity and also the ferromagnetic nature of the measure (see [GJ87] and the references therein). Let us set
 $\Lambda_\varepsilon:=\{\varepsilon x\in\mathbb{T}^3, x\in\mathbb{Z}^3\}$. Heuristically, the quantities
$\int |\nabla\phi(x)|^2dx$, $\int\phi^2(x)dx$, and $\int\phi^4(x)dx$ can be approximated  by
$\varepsilon\sum_{|x-y|=\varepsilon,x,y\in\Lambda_\varepsilon}(\phi(x)-\phi(y))^2$, $\varepsilon^3\sum_{x\in\Lambda_\varepsilon}\phi(x)^2$
and $\varepsilon^3\sum_{x\in\Lambda_\varepsilon}\phi(x)^4$, respectively, as $\varepsilon$ tends to zero. Thus heuristically
(1.1) can be approximated by the following heuristic probability measure $\mu_\varepsilon$:
$$N^{-1}_{\varepsilon}\Pi_{x\in\Lambda_\varepsilon}d\phi_x\exp\bigg(2\varepsilon\sum_{|x-y|=\varepsilon,x,y\in\Lambda_\varepsilon}
\phi(x)\phi(y)-(\varepsilon^3+12\varepsilon)\sum_{x\in\Lambda_\varepsilon}\phi^2(x)-\varepsilon^3\sum_{x\in\Lambda_\varepsilon}
\phi^4(x)\bigg),\eqno(1.3)$$
where $N_\varepsilon$ is the normalization constant. (1.3) is still just a heuristic expression, but it is indeed not hard to give a rigorous sense to
it (see [GJ87] and the references therein). We call this
the lattice $\Phi^4_3$-field measure. From $\mu_\varepsilon$ by deriving suitable bounds on its moments and choosing
subsequences if necessary one gets limit measures  by weak convergence. These are then the continuum $\Phi^4_3$-field measures.

The following  stochastic PDEs on $\Lambda_\varepsilon$ are the stochastic quantizations associated with the lattice $\Phi^4_3$-field measure:
$$\aligned d \Phi^{\varepsilon}(t,x)=&(\Delta_\varepsilon \Phi^{\varepsilon}(t,x)-(\Phi^{\varepsilon})^3(t,x) +(3C_0^{\varepsilon}-9C_1^\varepsilon)
\Phi^{\varepsilon}(t,x))dt\\&+ \varepsilon^{-3/2}dW_\varepsilon(t,x)
\\\Phi^\varepsilon(0)=&\Phi^\varepsilon_0.\endaligned(1.4)$$
Here  $W_\varepsilon(t)=\{W(t,x)\}_{x\in\Lambda_\varepsilon}$ is a family of independent  Brownian motions and $C_0^{\varepsilon}$ and $C_1^\varepsilon$ are
defined as below. For $x\in\Lambda_\varepsilon$ define
$$\Delta_\varepsilon f(x):=\varepsilon^{-2}\sum_{y\in\Lambda_\varepsilon, y\sim x}(f(y)-f(x)),$$
and the nearest neighbor relation $x\sim y$ is to be understood with periodic boundary conditions on $\Lambda_\varepsilon$. We emphasize that to
make sense of (1.2) we need to renormalise  some ill-defined terms in (1.2). This is done by adding the renormalisation terms
$C_0^{\varepsilon}\Phi^\varepsilon$ and $C_1^\varepsilon\Phi^\varepsilon$ in the approximating equation (1.4).
In this paper we prove that the dynamical lattice approximation converge to the dynamical $\Phi^4_3$ model. This problem is also related to the
convergence of a rescaled discrete spin
system to the solution to the dynamical $\Phi_3^4$ model (see [MW14] for the dynamical $\Phi^4_2$ model).

In one dimensional case, approximations to general stochastic partial differential equations driven by space-time white noise have been very well studied (see
[Gy98, Gy99, DG01, HMW14] and the reference therein). In [GP15] the authors study the Sasamoto-Spohn type discretizations of the
conservative stochastic Burgers equation. In three dimensional case we also study discrete approximations to N-S equations (see [ZZ14a]).

In this paper we use the paracontrolled distribution method to prove that the solutions to the lattice approximating equation  converge to the solution of the dynamical $\Phi^4_3$ model.
The theory of paracontrolled distributions combines the idea of Gubinelli's controlled rough path [Gub04] and Bony's paraproduct [Bon84], which is defined as follows:
Let $\Delta_jf$ be the jth Littlewood-Paley block of a (Scharwtz) distribution $f$. Define for distributions $f$ and $g$
$$\pi_<(f,g)=\pi_>(g,f)=\sum_{j\geq-1}\sum_{i<j-1}\Delta_if\Delta_jg, \quad\pi_0(f,g)=\sum_{|i-j|\leq1}\Delta_if\Delta_jg.$$
Formally $fg=\pi_<(f,g)+\pi_0(f,g)+\pi_>(f,g)$. Observing that if $f$ is regular $\pi_<(f,g)$ behaves like $g$ and is the only
term in the Bony's paraproduct not improving the regularities, the authors in [GIP13]  consider a paracontrolled ansatz of the type
$$u=\pi_<(u',g)+u^\sharp,$$ where $\pi_<(u',g)$ represents the "bad-part" of the solution, $u'$ is some suitable function and $g$ is some functional of the Gaussian field and $u^\sharp$ is regular enough to define the multiplication required. Then to make sense of the product of $uf$ we only need to define $gf$.

Using the  paracontrolled distribution method, to perform the lattice approximation of the dynamical $\Phi^4_3$ model we will meet the projection operators
$P_N$, which do not commute with the paraproduct. Here we use a random operator technique from [GP15] to handle this operators. However, for the $\Phi_3^4$ model
this technique is not enough for our case and we have to estimate an additional error term $D_N$  by stochastic calculations in Section 6.4 (see Remark 4.5).

\textbf{Framework and main result}

For $N\geq1$, let $\Lambda_N=\{-N,-(N-1),...,N\}^3$. Set $\varepsilon=\frac{2}{2N+1}$. Every point $k\in\Lambda_N$ can be identified with $x=\varepsilon k\in\Lambda_\varepsilon=\{x=(x_1,x_2,x_3)\in\varepsilon \mathbb{Z}^3:-1<x_1,x_2,x_3<1\}.$
We view $\Lambda_\varepsilon$ as a discretisation of the continuous three-dimensional torus $\mathbb{T}^3$ identified with $[-1,1]^3$. In the following for similicity we fix
a cylindrical Wiener process in (1.2)  on $L^2(\mathbb{T}^3)$  given by $2^{-3/2}\sum_k\beta_ke^{\iota\pi k\cdot x}$ for $x\in\mathbb{T}^3$ and restrict it on
$L^2(\Lambda_\varepsilon)$ as $W_N=2^{-3/2}\sum_{|k|_\infty\leq N}\beta_ke^{\iota\pi k\cdot x}$ for $x\in\Lambda_\varepsilon$, which is also a cylindrical Wiener process on
$L^2(\Lambda_\varepsilon)$. Here $\{\beta_k\}$ is a family of independent Brownian motions. Then for (1.4) and fixed $N$ it is a finite dimensional SDE and we could easily obtain  existence and uniqueness of solutions to (1.4) by [PR07], which implies that the solution to (1.4) has the same distribution as the solution to the following equation:
$$\aligned d \Phi^{\varepsilon}(t,x)=&(\Delta_\varepsilon \Phi^{\varepsilon}(t,x)-(\Phi^{\varepsilon})^3(t,x) +(3C_0^{\varepsilon}-9C_1^\varepsilon)
\Phi^{\varepsilon}(t,x))dt\\&+ dW_N(t,x)
\\\Phi^\varepsilon(0)=&\Phi^\varepsilon_0.\endaligned\eqno(1.5)$$

Following [MW14] we discuss a suitable extension of functions  defined on $\Lambda_\varepsilon$ onto all of the torus $\mathbb{T}^3$ (which we identify with the interval $[-1,1]^3$). For any function $Y:\Lambda_\varepsilon\rightarrow\mathbb{R}$, we define the discrete Fourier transform $\hat{Y}$ through
$$\hat{Y}(k)=\left\{\begin{array}{ll}\sum_{x\in\Lambda_\varepsilon}\varepsilon^3Y(x)e^{-\imath\pi k\cdot x},&\ \ \ \ \textrm{ if } k\in\{-N,...,N\}^3\\0&\ \ \ \ \textrm{ if } k\in\mathbb{Z}^3\backslash\{-N,...,N\}^3.\end{array}\right.$$
In this context Fourier inversion states
$$Y(x)=\frac{1}{8}\sum_{k\in\mathbb{Z}^3}\hat{Y}(k)e^{\imath\pi k\cdot x} \textrm{ for all } x\in\Lambda_\varepsilon.\eqno(1.6)$$
It is thus natural to extend $Y$ to all of $\mathbb{T}^3$ by taking (1.6) as a definition of $Y(x)$ for $x\in\mathbb{T}^3\backslash\Lambda_\varepsilon$. More explicitly, for $Y:\Lambda_\varepsilon\rightarrow\mathbb{R}$ we define $(\textrm{Ext} Y):\mathbb{T}^3\rightarrow\mathbb{R}$ as
$$\aligned\textrm{Ext} Y(x)=&\frac{1}{2^3}\sum_{k\in\{-N,...,N\}^3}\sum_{y\in\Lambda_\varepsilon}\varepsilon^3e^{\imath\pi k\cdot(x-y)}Y(y).\endaligned$$
Let $P_t^\varepsilon=\textrm{Ext}e^{t\Delta_\varepsilon}$. By the definition of the operators $\Delta_\varepsilon$ we have
$$\widehat{e^{t\Delta_\varepsilon} v}(k)=\left\{\begin{array}{ll}e^{-|k|^2f(\varepsilon k)}\hat{v}(k),&\ \ \ \ \textrm{ if } k\in\{-N,...,N\}^3
\\0&\ \ \ \ \textrm{ if } k\in\mathbb{Z}^3\backslash\{-N,...,N\}^3.\end{array}\right.$$
Here  $$f(x)=\frac{4}{|x|^2}(\sin^2\frac{x_1\pi}{2}+\sin^2\frac{x_2\pi}{2}+\sin^2\frac{x_3\pi}{2}).$$

Now we extend the solution to all of $\mathbb{T}^3$.  In the following the Fourier transform and the inverse Fourier transform are denoted by $\mathcal{F}$ and $\mathcal{F}^{-1}$. It is easy to see that
 $$\textrm{Ext} \Phi^\varepsilon(t)=P_t^\varepsilon \textrm{Ext}\Phi^\varepsilon_0-\int_0^tP_{t-s}^\varepsilon Q_N[(\textrm{Ext} \Phi^\varepsilon)^3-(3C_0^{\varepsilon}
 -9C_1^\varepsilon) \textrm{Ext}\Phi^{\varepsilon}]ds+\int_0^tP_{t-s}^\varepsilon \textrm{Ext}dW_N,\eqno(1.7)$$
 where $Q_N u(x)=P_N u(x)+\Pi_N u(x)$ with  $$P_N=\mathcal{F}^{-1}1_{|k|_\infty\leq N}\mathcal{F},$$ and $\Pi_N$ is defined for $u$ satisfying supp$\mathcal{F}u\subset \{k:|k|_\infty\leq 3N\}$
 $$\aligned\Pi_Nu(x)=&\sum_{i_1,i_2,i_3\in\{-1,0,1\},\sum_{j=1}^3i_j^2\neq0} e_N^{i_1i_2i_3}\mathcal{F}^{-1}1_{k\in P^{i_1i_2i_3} }
 \mathcal{F}u(x)\\=&\sum_{i_1,i_2,i_3\in\{-1,0,1\},\sum_{j=1}^3i_j^2\neq0}P_N[e^{i_1i_2i_3}_N u]\endaligned$$ with $P^{i_1i_2i_3}=\{k:k^ji_j>N \textrm{ if } i_j=-1,1;|k^j|\leq N,
 \textrm{ if } i_j=0\}$ is a rectangular division of $\mathbb{Z}^3\backslash \{k\in\mathbb{Z}^3, |k|_\infty\leq N\}$,  $e^{i_1i_2i_3}_N=\Pi_{j=1}^3e^{-\imath\pi(2N+1)i_jx^j}$ and $|k|_\infty=\max(|k^1|,|k^2|,|k^3|)$.

Now choose $C_0^\varepsilon$ as in (6.3) and $$C_1^\varepsilon=C_{11}^\varepsilon+\sum_{i_1,i_2,i_3\in\{-1,0,1\},\sum_{j=1}^3i_j^2\neq0}C_{12}^{\varepsilon,i_1i_2i_3},$$
with $C_{11}^\varepsilon,C_{12}^{\varepsilon,i_1i_2i_3}$ as in (6.4) and (6.5) respectively.  In the following  we omit the summation with respect to $i_1,i_2,i_3$ if there's no confusion.
\vskip.10in
The main result of this paper is the following theorem:
\vskip.10in
\th{Theorem 1.1} Let $z\in (1/2,2/3)$  and $\Phi_0\in \mathcal{C}^{-z}$. Let $(\Phi,\tau)$ be the unique ( maximal in time) solution to (1.2)
and let for $\varepsilon\in(0,1)$ the function $\Phi^\varepsilon$ be the unique  solution to (1.5) on $[0,\infty)$. If the initial data satisfy
$\textrm{Ext}\Phi_0^\varepsilon-\Phi_0\rightarrow0 \textrm{ in }\mathcal{C}^{-z}$ then there exists a sequence of random time  $\tau_L$ such that  $\lim_{L\rightarrow\infty}\tau_L=\tau$ and $$\sup_{t\in[0,\tau_L]}\|\textrm{Ext}\Phi^\varepsilon-\Phi\|_{-z}\rightarrow0\quad \textrm{ in probability}, \quad \textrm{ as } \varepsilon\rightarrow0 .$$
\vskip.10in
\th{Remark 1.2} (1) Existence and uniqueness of $(\Phi,\tau)$ has been obtained in [Hai14, CC13]. For the definition of $\mathcal{C}^{-z}$ and norm $\|\cdot\|_{-z}$ see Section 2.

(2) The constant $C_1^\varepsilon$ is divided into two parts: $C_{11}^\varepsilon$ and $C_{12}^\varepsilon$ which correspond to terms with $P_N$  and $\Pi_N$ respectively. In fact $C_0^\varepsilon\backsimeq
\frac{1}{\varepsilon}, C_{11}^\varepsilon\backsimeq-\log\varepsilon$ and $C_{12}^{\varepsilon,i_1i_2i_3}\backsimeq 1$.

\vskip.10in
 The structure of the paper is organized as follows. In Section 2, we recall some basic notions and results for the paracontrolled distribution method. In Section 3 we prove some estimates for the approximating operators. In Section 4 we use the paracontrolled distribution method to prove  uniform bounds for the lattice approximation equations. In Section 5 we give the proof of our main result. In Section 6  convergence of several stochastic terms is proved.
\section{Besov spaces and paraproduct}

In the following we recall the definitions and some properties of Besov spaces and paraproducts. For a general introduction to these theories we refer to [BCD11, GIP13].
First we introduce the following notations. The space of real valued infinitely differentiable functions of compact support is denoted by $\mathcal{D}(\mathbb{R}^d)$ or $\mathcal{D}$. The space of Schwartz functions is denoted by $\mathcal{S}(\mathbb{R}^d)$. Its dual, the space of tempered distributions is denoted by $\mathcal{S}'(\mathbb{R}^d)$.

 Let $\chi,\theta\in \mathcal{D}$ be nonnegative radial functions on $\mathbb{R}^d$, such that

i. the support of $\chi$ is contained in a ball and the support of $\theta$ is contained in an annulus;

ii. $\chi(z)+\sum_{j\geq0}\theta(2^{-j}z)=1$ for all $z\in \mathbb{R}^d$.

iii. $\textrm{supp}(\chi)\cap \textrm{supp}(\theta(2^{-j}\cdot))=\emptyset$ for $j\geq1$ and $\textrm{supp}(\theta(2^{-i}\cdot))\cap \textrm{supp}(\theta(2^{-j}\cdot))=\emptyset$ for $|i-j|>1$.

We call such $(\chi,\theta)$ dyadic partition of unity, and for the existence of dyadic partitions of unity we refer to [BCD11, Proposition 2.10]. The Littlewood-Paley blocks are now defined as
$$\Delta_{-1}u=\mathcal{F}^{-1}(\chi\mathcal{F}u)\quad \Delta_{j}u=\mathcal{F}^{-1}(\theta(2^{-j}\cdot)\mathcal{F}u).$$

For $\alpha\in \mathbb{R}$, the H\"{o}lder-Besov space $\mathcal{C}^\alpha$ is given by $\mathcal{C}^\alpha=B^\alpha_{\infty,\infty}(\mathbb{R}^d,\mathbb{R}^n)$, where for $p,q\in [1,\infty]$ we define
$$B^\alpha_{p,q}(\mathbb{R}^d)=\{u\in\mathcal{S}'(\mathbb{R}^d):\|u\|_{B^\alpha_{p,q}}=(\sum_{j\geq-1}(2^{j\alpha}\|\Delta_ju\|_{L^p})^q)^{1/q}<\infty\},$$
with the usual interpretation as $l^\infty$ norm in case $q=\infty$. We write $\|\cdot\|_{\alpha}$ instead of $\|\cdot\|_{B^\alpha_{\infty,\infty}}$ in the following for simplicity. We also use $C_TE$ to denote $C([0,T];E)$.

We point out that everything above and everything that follows can be applied to distributions on the torus (see [S85, SW71]). More precisely, let $\mathcal{S}'(\mathbb{T}^d)$ be the space of distributions on $\mathbb{T}^d$. Therefore, Besov spaces on the torus with general indices $p,q\in[1,\infty]$ are defined as
$$B^\alpha_{p,q}(\mathbb{T}^d)=\{u\in\mathcal{S}'(\mathbb{T}^d):\|u\|_{B^\alpha_{p,q}}=(\sum_{j\geq-1}(2^{j\alpha}\|\Delta_ju\|_{L^p(\mathbb{T}^d)})^q)^{1/q}<\infty\}.$$
 We  will need the following Besov embedding theorem on the torus (c.f. [GIP13, Lemma 41]):
\vskip.10in
 \th{Lemma 2.1} Let $1\leq p_1\leq p_2\leq\infty$ and $1\leq q_1\leq q_2\leq\infty$, and let $\alpha\in\mathbb{R}$. Then $B^\alpha_{p_1,q_1}(\mathbb{T}^d)$ is continuously embedded in $B^{\alpha-d(1/p_1-1/p_2)}_{p_2,q_2}(\mathbb{T}^d)$.
\vskip.10in

 Now we recall the following paraproduct introduced by Bony (see [Bon81]). In general, the product $fg$ of two distributions $f\in \mathcal{C}^\alpha, g\in \mathcal{C}^\beta$ is well defined if and only if $\alpha+\beta>0$. In terms of Littlewood-Paley blocks, the product $fg$ can be formally decomposed as
 $$fg=\sum_{j\geq-1}\sum_{i\geq-1}\Delta_if\Delta_jg=\pi_<(f,g)+\pi_0(f,g)+\pi_>(f,g),$$
 with $$\pi_<(f,g)=\pi_>(g,f)=\sum_{j\geq-1}\sum_{i<j-1}\Delta_if\Delta_jg, \quad\pi_0(f,g)=\sum_{|i-j|\leq1}\Delta_if\Delta_jg.$$
We also use the notation for $j\geq0$
$$S_jf=\sum_{i\leq j-1}\Delta_if.$$
Moreover define $$\psi_<(k_1,k_2)=\sum_{j\geq-1}\sum_{i<j-1}\theta(2^{-i}k_1)\theta(2^{-j}k_2)$$ and $$\psi_0(k_1,k_2)=\sum_{|i-j|\leq1}\theta(2^{-i}k_1)\theta(2^{-j}k_2).$$
We will use without comment that $\|\cdot\|_\alpha\leq\|\cdot\|_\beta$ for $\alpha\leq\beta$, that $\|\cdot\|_{L^\infty}\lesssim \|\cdot\|_\alpha$ for $\alpha>0$, and that $\|\cdot\|_\alpha\lesssim\|\cdot\|_{L^\infty}$ for $\alpha\leq0$. We will also use that $\|S_ju\|_{L^\infty}\lesssim 2^{-j\alpha}\|u\|_\alpha$ for $\alpha<0, j\geq0$ and $u\in \mathcal{C}^\alpha$, where $\|\cdot\|_\alpha$ denotes the norm in $\mathcal{C}^\alpha, \alpha\in\mathbb{R}$.

\vskip.10in
 The basic result about these bilinear operations is given by the following estimates:
\vskip.10in
 \th{Lemma 2.2}(Paraproduct estimates, [Bon 81, GIP13, Lemma 2]) For any $\beta\in \mathbb{R}$ we have
 $$\|\pi_<(f,g)\|_\beta\lesssim \|f\|_{L^\infty}\|g\|_\beta\quad f\in L^\infty, g\in \mathcal{C}^\beta,$$
 and for $\alpha<0$ furthermore
 $$\|\pi_<(f,g)\|_{\alpha+\beta}\lesssim \|f\|_{\alpha}\|g\|_\beta\quad f\in \mathcal{C}^\alpha, g\in \mathcal{C}^\beta.$$
 For $\alpha+\beta>0$ we have
 $$\|\pi_0(f,g)\|_{\alpha+\beta}\lesssim \|f\|_{\alpha}\|g\|_\beta\quad f\in \mathcal{C}^\alpha, g\in \mathcal{C}^\beta.$$

\vskip.10in
 The following basic commutator lemma is important for our use:
\vskip.10in
 \th{Lemma 2.3}([GIP13, Lemma 5]) Assume that $\alpha\in (0,1)$ and $\beta,\gamma\in \mathbb{R}$ are such that $\alpha+\beta+\gamma>0$ and $\beta+\gamma<0$. Then for smooth $f,g,h,$ the trilinear operator
 $$C(f,g,h)=\pi_0(\pi_<(f,g),h)-f\pi_0(g,h)$$ allows for the bound
 $$\|C(f,g,h)\|_{\alpha+\beta+\gamma}\lesssim\|f\|_\alpha\|g\|_\beta\|h\|_\gamma.$$
 Thus, $C$ can be uniquely extended to a bounded trilinear operator from $\mathcal{C}^\alpha\times \mathcal{C}^\beta \times \mathcal{C}^\gamma$ to $ \mathcal{C}^{\alpha+\beta+\gamma}$.

\vskip.10in

Now we recall the following estimate for heat semigroup $P_t:=e^{t\Delta}$.
\vskip.10in
\th{Lemma 2.4}([GIP13, Lemma 47]) Let $u\in \mathcal{C}^\alpha$ for some $\alpha\in \mathbb{R}$. Then for every $\delta\geq0$
$$\|P_tu\|_{\alpha+\delta}\lesssim t^{-\delta/2}\|u\|_\alpha.$$

\vskip.10in

\th{Lemma 2.5} ([CC13, Lemma A.1]) Let $u\in \mathcal{C}^\alpha$ for some $\alpha<1$ and $v\in \mathcal{C}^\beta$ for some $\beta\in \mathbb{R}$. Then for $\delta\geq \alpha+\beta$
$$\|P_t\pi_{<}(u,v)-\pi_<(u,P_t v)\|_{\delta}\lesssim t^{\frac{\alpha+\beta-\delta}{2}}\|u\|_\alpha\|v\|_\beta.$$
\vskip.10in

\th{Lemma 2.6} ([CC13, Lemma 2.5])  Let $u\in \mathcal{C}^{\alpha+\delta}$ for some $\alpha\in \mathbb{R},\delta>0$. Then for every $ t\geq0$
$$\|(P_t-I)u\|_{\alpha}\lesssim  t^{\delta/2}\|u\|_{\alpha+\delta}.$$
\vskip.10in
We also have the following result.
\vskip.10in
\th{Lemma 2.7}(Bernstein type lemma) Let $u\in \mathcal{C}^{\alpha}$ for some $\alpha\in\mathbb{R}$.

1) If $supp\mathcal{F}u\subset\{k:|k|\leq CN\}$ for some $C>0$ then for $\beta>\alpha$
$$\|u\|_{\beta}\lesssim N^{\beta-\alpha}\|u\|_\alpha.$$

2)If $supp\mathcal{F}u\subset\{k:|k|>CN\}$ for some $C>0$ then for $\beta<\alpha$
$$\|u\|_{\beta}\lesssim N^{\beta-\alpha}\|u\|_\alpha.$$
Here all the constants we omit are independent of $N$.

\proof We have $$\|u\|_\beta=\sup_j2^{j\beta}\|\Delta_j u\|_{L^\infty}=\sup_j2^{j(\beta-\alpha)}2^{j\alpha}\|\Delta_j u\|_{L^\infty}.$$
For the first case we have $\Delta_ju\neq0$ iff  $2^j\lesssim N$, which implies the first result.
If $supp\mathcal{F}u\subset\{k:|k|>CN\}$ we have $\Delta_ju\neq0$ iff $2^j\gtrsim N$ which implies the second result. $\hfill\Box$

\section{Estimates for the approximated operators}

Now we prove the following estimates for the approximated operators on $\mathbb{T}^3$. First we consider estimate for $P_N$ and $\Pi_N$:

\vskip.10in
\th{Lemma 3.1} Let $u\in \mathcal{C}^{\alpha}$ for some $\alpha\in\mathbb{R}$. Then for any $\kappa>0$ small enough we have the following estimate:

(1) (Estimate for $P_N$)
$$\|P_Nu\|_{\alpha-\kappa}\lesssim \|u\|_{\alpha},
\quad\|(I-P_N)u\|_{\alpha-\kappa}\lesssim N^{-\kappa/2}\|u\|_{\alpha}.$$

(2) (Estimate for $\Pi_N$) If $\alpha>0$  then for $u$ satisfying $supp\mathcal{F}u\subset \{k:|k|_\infty\leq 3N\}$
$$\|\Pi_Nu\|_{\alpha-\kappa}\lesssim N^{-\kappa/2}\|u\|_{\alpha}.$$
If $\alpha<0$ and $supp\mathcal{F}u\subset\{k:|k|_\infty\leq N\}$  then
$$\|e_N^{i_1i_2i_3}u\|_{\alpha-\kappa}\lesssim N^{-\kappa/2}\|u\|_{\alpha}.$$
Here all the constants we omit are independent of $N$.

\proof We have  for $p$ large enough $$\|P_Nu\|_{\alpha-\kappa}\lesssim \|P_Nu\|_{B^{\alpha}_{p,\infty}}\lesssim \|u\|_{B^{\alpha}_{p,\infty}}\lesssim \|u\|_\alpha,$$
where in the first inequality we used Lemma 2.1 and in the second inequality we used that $1_{|k|_\infty\leq N}$ is an $L^p$ multiplier.
Similarly
$$\|(I-P_N)u\|_{\alpha-\kappa}\lesssim N^{-\kappa/2}\|(I-P_N)u\|_{\alpha-\kappa/2}\lesssim N^{-\kappa/2}\|u\|_\alpha,$$
where in the first inequality we used Lemma 2.7 and in the second inequality we used the result for $P_N$.
Moreover for $\alpha>5\kappa/4$
 $$\aligned&\|\Pi_Nu\|_{\alpha-\kappa}\lesssim N^{\alpha-5\kappa/4}\|\Pi_Nu\|_{\kappa/4}
\lesssim N^{\alpha-\kappa}\|\mathcal{F}^{-1}1_{k\in P^{i_1i_2i_3} }\mathcal{F}u\|_{\kappa/4}
 \\\lesssim&N^{-\kappa/2}\|\mathcal{F}^{-1}1_{k\in P^{i_1i_2i_3} }\mathcal{F}u\|_{\alpha-\kappa/4}\lesssim N^{-\kappa/2}\|u\|_{\alpha}.\endaligned$$
 Here in the first and third inequalities we used Lemma 2.7, in the second inequality we used that $\|e_N^{i_1i_2i_3}\|_{\kappa/2}\lesssim N^{\kappa/2}$ and in the last inequality we used similar argument  for $P_N$ since $1_{k\in P^{i_1i_2i_3} }$ is an $L^p$ multiplier.
Similarly for $\alpha<0$ $$\|e_N^{i_1i_2i_3}u\|_{\alpha-\kappa}\lesssim N^{\alpha-3\kappa/2}\|e_N^{i_1i_2i_3}u\|_{\kappa/2}\lesssim N^{\alpha-\kappa}\|u\|_{\kappa/2}\lesssim N^{-\kappa/2}\|u\|_{\alpha}.$$
Here in the first inequality we used $supp\mathcal{F}(e_N^{i_1i_2i_3}u)\subset\{k:|k|>N\}$ and Lemma 2.7.
Thus the result follows.
$\hfill\Box$
\vskip.10in
Now want to prove estimates for $P_t^\varepsilon=P_Ne^{t\Delta_\varepsilon}$. In fact, $$P_t^\varepsilon=\mathcal{F}^{-1}1_{|k|_\infty\leq N}e^{-t|k|^2f(\varepsilon k)}\varphi(\varepsilon k)\mathcal{F}=P_N\tilde{P}_t^\varepsilon,$$
with $$\tilde{P}_t^\varepsilon:=\mathcal{F}^{-1}e^{-t|k|^2f(\varepsilon k)}\varphi(\varepsilon k)\mathcal{F},$$
where $\varphi$ is a smooth function and equals $1$ on $\{|x|_\infty\leq 1\}$ with support in $\{|x|\leq 1.8\}$. Then by a similar argument as [GIP13, Lemma 47] we have the following result:
\vskip.10in
\th{Lemma 3.2} Let $u\in \mathcal{C}^\alpha$ for some $\alpha\in \mathbb{R}$. Then for every $\delta\geq0,\kappa>0, t>0,$
$$\|P^\varepsilon_tu\|_{\alpha+\delta-\kappa}\lesssim t^{-\delta/2}\|u\|_\alpha,$$
$$\|(P^\varepsilon_t-P_t)u\|_{\alpha+\delta-\kappa}\lesssim \varepsilon^{\kappa/2} t^{-\delta/2}\|u\|_{\alpha}.$$
Here the constants we omit are independent of $N$.

\proof For the first result by Lemma 3.1 it suffices to prove $$\|\tilde{P}^\varepsilon_tu\|_{\alpha+\delta}\lesssim t^{-\delta/2}\|u\|_\alpha.\eqno(3.1)$$
In the following we consider (3.1) and have for $j\geq0$ $$\aligned &\|\Delta_j\tilde{P}^\varepsilon_tu\|_{L^\infty}=\|\mathcal{F}^{-1}\theta_j\phi^\varepsilon\mathcal{F}u\|_{L^\infty}
=\|\mathcal{F}^{-1}\theta_j\tilde{\theta}(2^{-j}\cdot)\phi^\varepsilon\mathcal{F}u\|_{L^\infty}\\\leq& \|\mathcal{F}^{-1}(\phi^\varepsilon\tilde{\theta}(2^{-j}\cdot))\|_{L^1}\|\Delta_ju\|_{L^\infty}.\endaligned$$
Here $$\phi^\varepsilon(\xi)=e^{-t|\xi|^2f(\varepsilon \xi)}\varphi(\varepsilon \xi),$$ and $\tilde{\theta}$ be a smooth function supported in an annulus such that $\tilde{\theta}\theta=\theta.$
Then we get that for $\delta\geq0$ $$\aligned&\|\mathcal{F}^{-1}(\phi^\varepsilon\tilde{\theta}(2^{-j}\cdot))\|_{L^1}
=\|\mathcal{F}^{-1}(\phi^\varepsilon(2^{j}\cdot)\tilde{\theta})\|_{L^1}
\lesssim \|(1-\Delta)^2((\phi^\varepsilon(2^{j}\cdot)\tilde{\theta})\|_{L^1}\\\lesssim& \sum_{0\leq|k|\leq 4}2^{j|k|}\|(D_k\phi^\varepsilon)(2^j\cdot)|_{\cdot\in\textrm{supp} \tilde{\theta}}\|_{L^\infty}\lesssim \sum_{0\leq|k|\leq 4}2^{j|k|}\frac{1}{2^{j|k|}(2^j\sqrt{t})^{\delta}}\lesssim (2^j\sqrt{t})^{-\delta}.\endaligned$$
Here in the third inequality we used that  $f(\varepsilon \xi)\geq c>0$ and $|\varepsilon\xi|\lesssim1$ on the support of $\varphi^\varepsilon$ which implies that for any multiindices $k$ satisfying $|k|\leq4$ and every $\delta\geq0$ $|D_k\phi^\varepsilon(\xi)|\lesssim \frac{1}{|\xi|^{|k|+\delta}t^{\delta/2}}$.
For $j=-1$ we can use Bernstein's lemma to obtain the estimate. Thus (3.1) follows.

For the second result we have $$P^\varepsilon_t-P_t=P_N(\tilde{P}_t^\varepsilon-P_t)+(I-P_N)P_t.$$
By Lemmas 2.4 and 3.1 it is sufficient to consider $\tilde{P}_t^\varepsilon-P_t$. Since  $\phi^\varepsilon(\xi)-\phi(\xi)=\varphi(\varepsilon\xi)(e^{-t|\xi|^2f(\varepsilon \xi)}-e^{-t|\xi|^2})+(\varphi(\varepsilon\xi)-1)e^{-t|\xi|^2}$  and $|\varphi(\varepsilon\xi)-1|\lesssim |\varepsilon \xi|^\delta,|f(\varepsilon\xi)-1|\lesssim |\varepsilon \xi|^\delta$,
we obtain that for any multiindices $k$ satisfying $|k|\leq4$ and every $\delta\geq0,0<\eta<1$ $|D_k(\phi^\varepsilon-\phi)(\xi)|\leq \frac{(\varepsilon|\xi|)^\eta}{|\xi|^{|k|+\delta}t^{\delta/2}}$. Thus the second result follows by a similar argument as the calculation for (3.1).
$\hfill\Box$

\vskip.10in
Now we prove a commutator estimate for $P_t^\varepsilon$. However $P_N$ does not commute with paraproduct, we could only obtain the commutator estimate for $\tilde{P}_t^\varepsilon$.
\vskip.10in
\th{Lemma 3.3} Let $u\in \mathcal{C}^\alpha$ for some $\alpha<1$ and $v\in \mathcal{C}^\beta$ for some $\beta\in \mathbb{R}$. Then for $\delta\geq \alpha+\beta$ and any $\kappa>0$
$$\|P_t^\varepsilon\pi_{<}(u,v)-P_N\pi_<(u,\tilde{P}_t^\varepsilon v)\|_{\delta-\kappa}\lesssim t^{\frac{\alpha+\beta-\delta}{2}}\|u\|_\alpha\|v\|_\beta,\eqno(3.2)$$
$$\|(P_t^\varepsilon-P_t)\pi_{<}(u,v)-P_N\pi_<(u,\tilde{P}_t^\varepsilon v)-\pi_<(u,P_tv)\|_{\delta-\kappa}\lesssim \varepsilon^{\kappa/2} t^{\frac{\alpha+\beta-\delta}{2}}\|u\|_\alpha\|v\|_\beta.\eqno(3.3)$$
Here the constants we omit are independent of $N$.

\proof We have $$P_t^\varepsilon\pi_{<}(u,v)-P_N\pi_<(u,\tilde{P}_t^\varepsilon v)=P_N(\tilde{P}_t^\varepsilon\pi_{<}(u,v)-\pi_<(u,\tilde{P}_t^\varepsilon v)).$$
By Lemma 3.1 it suffices to prove that $$\|\tilde{P}_t^\varepsilon\pi_{<}(u,v)-\pi_<(u,\tilde{P}_t^\varepsilon v)\|_{\delta}\lesssim t^{\frac{\alpha+\beta-\delta}{2}}\|u\|_\alpha\|v\|_\beta.\eqno(3.4)$$
In fact, we have
$$\tilde{P}_t^\varepsilon\pi_{<}(u,v)-\pi_<(u,\tilde{P}_t^\varepsilon v)=\sum_{j=-1}^{\infty}(\tilde{P}_t^\varepsilon(S_{j-1}u\Delta_jv)-S_{j-1}u\tilde{P}_t^\varepsilon \Delta_jv).$$
We also have that the Fourier transform of $\tilde{P}_t^\varepsilon(S_{j-1}u\Delta_jv)-S_{j-1}u\tilde{P}_t^\varepsilon \Delta_jv$ has its support in a suitable annulus $2^j\mathcal{A}$. Let $\psi\in \mathcal{D}(\mathbb{R}^3)$ with support in an annulus $\tilde{\mathcal{A}}$ be such that $\psi=1$ on $\mathcal{A}$.

Thus by the same argument as the proof of [CC13, Lemma A.1] we obtain that
$$\aligned &\|[(\psi(2^{-j}\cdot)\phi^\varepsilon)(D),S_{j-1}u]\Delta_jv\|_{L^\infty}
\\\lesssim&\sum_{\eta\in \mathbb{N}^d,|\eta|=1}\|x^\eta\mathcal{F}^{-1}(\psi(2^{-j}\cdot)\phi^\varepsilon)\|_{L^1}\|\partial^\eta S_{j-1}u\|_{L^\infty}\|\Delta_jv\|_{L^\infty}.\endaligned$$
Here $(\psi(2^{-j}\cdot)\phi^\varepsilon)(D)u=\mathcal{F}^{-1}(\psi(2^{-j}\cdot)\phi^\varepsilon\mathcal{F}u)$, $[(\psi(2^{-j}\cdot)\phi^\varepsilon)(D),S_{j-1}u]$ denotes the commutator.

Now we have that
$$\aligned &\|x^\eta\mathcal{F}^{-1}(\psi(2^{-j}\cdot)\phi^\varepsilon)\|_{L^1}\\\leq& 2^{-j}\|\mathcal{F}^{-1}(\partial^\eta\psi)(2^{-j}\cdot)\phi^\varepsilon)\|_{L^1}+\|\mathcal{F}^{-1}(\psi(2^{-j}\cdot)\partial^\eta\phi^\varepsilon)\|_{L^1}
\\=&2^{-j}\|\mathcal{F}^{-1}(\partial^\eta\psi(\cdot)\phi^\varepsilon(2^j\cdot))\|_{L^1}+\|\mathcal{F}^{-1}(\psi(\cdot)\partial^\eta\phi^\varepsilon(2^j\cdot))\|_{L^1}
\\\lesssim&2^{-j}\|(1+|\cdot|^2)^{2}\mathcal{F}^{-1}(\partial^\eta\psi(\cdot)\phi^\varepsilon(2^j\cdot))\|_{L^\infty}+ \|(1+|\cdot|^2)^{2}\mathcal{F}^{-1}(\psi(\cdot)\partial^\eta\phi^\varepsilon( 2^j\cdot))\|_{L^\infty} \\=&2^{-j}\|\mathcal{F}^{-1}((1-\Delta)^{2}(\partial^\eta\psi(\cdot)\phi^\varepsilon( 2^j\cdot)))\|_{L^\infty}+ \|\mathcal{F}^{-1}((1-\Delta)^{2}(\psi(\cdot)\partial^\eta\phi^\varepsilon( 2^j\cdot)))\|_{L^\infty}
\\\lesssim&2^{-j}\|(1-\Delta)^{2}(\partial^\eta\psi(\cdot)\phi^\varepsilon( 2^j\cdot))\|_{L^1}+ \|(1-\Delta)^{2}(\psi(\cdot)\partial^\eta\phi^\varepsilon( 2^j\cdot))\|_{L^1}
 \\\lesssim&2^{-j}\sum_{0\leq|m|\leq 4}( 2^j)^{|m|}\frac{t^{-\mu}2^{-2j\mu}}{( 2^j)^{|m|}}+ \sum_{|m|\leq 5}( 2^j)^{|m|}\frac{t^{-\mu}2^{-2j\mu}}{( 2^j)^{|m|+1}}\\\lesssim&2^{-j}t^{-\mu}2^{-2j\mu},\endaligned$$
 where  in the fourth inequality we used $|D^m\phi^\varepsilon(\xi)|\lesssim |\xi|^{-|m|}t^{-\mu}|\xi|^{-2\mu}, \mu\geq0$ for any multiindices $m$ satisfying $|m|\leq5$.
Hence we get that
$$\|[\psi(2^{-j}\cdot)\phi^\varepsilon(D),S_{j-1}u]\Delta_jv\|_{L^\infty}\lesssim t^{\frac{\alpha+\beta-\delta}{2}}2^{j(\alpha+\beta-\delta)}2^{-j(\alpha+\beta)}\|u\|_\alpha\|v\|_\beta,$$which yields (3.4) by the same argument as in the proof of [CC13, Lemma A.1].

 Moreover we have$$\aligned &(P_t^\varepsilon-P_t)\pi_{<}(u,v)-P_N\pi_<(u,\tilde{P}_t^\varepsilon v)-\pi_<(u,P_tv)\\=&P_N[(\tilde{P}_t^\varepsilon-P_t)\pi_{<}(u,v)-\pi_<(u,(\tilde{P}_t^\varepsilon-P_t) v)]-(I-P_N)(P_t\pi_{<}(u,v)-\pi_<(u,P_tv)).\endaligned$$
 The estimate for the second term can be obtained by Lemmas 2.5, 2.7 and 3.1. By a similar argument as Lemma 3.2 we obtain that for any multiindices $k$ satisfying $|k|\leq5$ and every $\delta\geq0,0<\eta<1$ $|D_k(\phi^\varepsilon-\phi)(\xi)|\leq \frac{(\varepsilon|\xi|)^\eta}{|\xi|^{|k|+\delta}t^{\delta/2}}$
 Thus (3.3) follows by a similar argument as the proof of (3.4).
$\hfill\Box$
\vskip.10in

Now we prove the following continuity result for $P_t^\varepsilon$.

\vskip.10in
\th{Lemma 3.4} Let $u\in \mathcal{C}^{\alpha+\delta}$ for some $\alpha\in \mathbb{R},0<\delta<1$. Then for every $\varepsilon\in (0,1),\kappa>0, t>s>0$
$$\|(P^\varepsilon_t-P^\varepsilon_s)u\|_{\alpha-\kappa}\lesssim  (t-s)^{\delta/2}\|u\|_{\alpha+\delta}.$$
Here the constants are independent of $N$.

\proof We have $(P^\varepsilon_t-P^\varepsilon_s)u=P_N(\tilde{P}^\varepsilon_t-\tilde{P}^\varepsilon_s)u$. By Lemma 3.1 it suffices to prove that
$$\|(\tilde{P}^\varepsilon_t-\tilde{P}^\varepsilon_s)u\|_{\alpha}\lesssim  (t-s)^{\delta/2}\|u\|_{\alpha+\delta}.$$
By $|1-e^{-(t-s)f(\varepsilon \xi)|\xi|^2}|\leq (t-s)^{\delta/2}|\xi|^\delta$  we obtain that
for any multiindices $k$ satisfying $|k|\leq4$ and any $\delta\geq0$ $|D_k(\phi^\varepsilon_t-\phi_s^\varepsilon)(\xi)|\lesssim \frac{(t-s)^{\delta/2}|\xi|^\delta}{|\xi|^{|k|}}$. Thus by a similar argument as Lemma 3.2  the result follows.$\hfill\Box$

\section{Paracontrolled analysis for the approximating equations}

Now let $u^\varepsilon=\textrm{Ext}\Phi^\varepsilon$ for simplicity and we have the following equation:
$$u^\varepsilon(t)=P_t^\varepsilon \textrm{Ext}\Phi_0^\varepsilon-\int_0^tP_{t-s}^\varepsilon Q_N[(u^\varepsilon)^3-(3C_0^{\varepsilon}-9C_1^\varepsilon) u^{\varepsilon}]ds+\int_0^tP_{t-s}^\varepsilon P_NdW.\eqno(4.1)$$

Therefore it suffices to prove the convergence result for solutions to (4.1). In this section we give an uniform estimate for solutions to (4.1) by using paracontrolled analysis.

In this section we fix $\delta,\beta,\kappa,\gamma>0$ satisfying
$$2z-1\geq\delta>2\kappa,\quad\beta>\frac{\delta}{2},\quad \beta+\frac{\delta}{2}+\kappa<\gamma,\quad 5\kappa+\frac{\delta}{2}+\beta+3\gamma<2-3z.$$

Now we split (4.1) into the following three equations:
$$ u_1^{\varepsilon}=\int_{-\infty}^t P_{t-s}^\varepsilon P_N dW,$$
$$u_2^{\varepsilon}=-\int_0^tP_{t-s}^\varepsilon Q_N[(u_1^{\varepsilon})^{\diamond,3}]ds$$
and
$$\aligned u_3^{\varepsilon}(t)=&P^\varepsilon_t(\textrm{Ext}\Phi^\varepsilon_0-u_1^{\varepsilon}(0))-\int_0^t
P_{t-s}^\varepsilon\bigg{[} Q_N[6u_1^{\varepsilon}\diamond u_2^{\varepsilon}u_3^{\varepsilon}+3u_1^{\varepsilon}(u_3^{\varepsilon})^2+3u_1^{\varepsilon}\diamond (u_2^{\varepsilon})^{2}+(u_2^{\varepsilon}+u_3^{\varepsilon})^3]\\&+P_N[3(u_1^{\varepsilon})^{\diamond,2} \diamond(u_2^{\varepsilon}+u_3^{\varepsilon})+3e_N^{i_1i_2i_3}(u_1^{\varepsilon})^{\diamond,2} \diamond(u_2^{\varepsilon}+u_3^{\varepsilon})-9\varphi^\varepsilon u^\varepsilon]\bigg{]}ds.\endaligned\eqno(4.2)$$
Here
$$(u^{\varepsilon}_1)^{\diamond,2} :=(u^{\varepsilon}_1)^{2}-C^{\varepsilon}_0,$$
$$(u^{\varepsilon}_1)^{\diamond,3} :=(u^{\varepsilon}_1)^3-3C^{\varepsilon}_0u^{\varepsilon}_1,$$
 $$u_1^{\varepsilon}\diamond u_2^{\varepsilon}:=u_2^\varepsilon u_1^{\varepsilon},$$
 $$u_1^{\varepsilon}\diamond (u_2^{\varepsilon})^2:=(u_2^\varepsilon)^2 u_1^{\varepsilon},$$
 $$\aligned (u_1^{\varepsilon})^{\diamond,2}\diamond u_2^{\varepsilon}:=&\pi_<(u_2^{\varepsilon},(u_1^{\varepsilon})^{\diamond,2})+\pi_>(u_2^{\varepsilon},(u_1^{\varepsilon})^{\diamond,2})+\pi_{0,\diamond}(u_2^{\varepsilon},(u_1^{\varepsilon})^{\diamond,2})   \\=&u_2^{\varepsilon}(u_1^{\varepsilon})^{\diamond,2}+3(C_{11}^\varepsilon+\varphi_1^\varepsilon) u_1^\varepsilon,\endaligned$$
 $$\aligned (u_1^{\varepsilon})^{\diamond,2}\diamond u_3^{\varepsilon}:=&\pi_<(u_3^{\varepsilon},(u_1^{\varepsilon})^{\diamond,2})+\pi_>(u_3^{\varepsilon},(u_1^{\varepsilon})^{\diamond,2})+\pi_{0,\diamond}(u_3^{\varepsilon},(u_1^{\varepsilon})^{\diamond,2})
 \\=&u_3^{\varepsilon}(u_1^{\varepsilon})^{\diamond,2}+3(C_{11}^\varepsilon+\varphi_1^\varepsilon) (u_2^\varepsilon+u_3^\varepsilon),\endaligned$$
 $$\aligned e_N^{i_1i_2i_3}(u_1^{\varepsilon})^{\diamond,2}\diamond u_2^{\varepsilon}:=&\pi_<(u_2^{\varepsilon},e_N^{i_1i_2i_3}(u_1^{\varepsilon})^{\diamond,2})+\pi_>(u_2^{\varepsilon},e_N^{i_1i_2i_3}(u_1^{\varepsilon})^{\diamond,2})+\pi_{0,\diamond}(u_2^{\varepsilon},e_N^{i_1i_2i_3}(u_1^{\varepsilon})^{\diamond,2})   \\=&e_N^{i_1i_2i_3}u_2^{\varepsilon}(u_1^{\varepsilon})^{\diamond,2}+3(C_{12}^{\varepsilon,i_1i_2i_3}+\varphi_2^{\varepsilon,i_1i_2i_3}) u_1^\varepsilon,\endaligned$$
 $$\aligned e_N^{i_1i_2i_3}(u_1^{\varepsilon})^{\diamond,2}\diamond u_3^{\varepsilon}:=&\pi_<(u_3^{\varepsilon},(u_1^{\varepsilon})^{\diamond,2}e_N^{i_1i_2i_3})+\pi_>(u_3^{\varepsilon},(u_1^{\varepsilon})^{\diamond,2}e_N^{i_1i_2i_3})+\pi_{0,\diamond}(u_3^{\varepsilon},(u_1^{\varepsilon})^{\diamond,2}e_N^{i_1i_2i_3})
 \\=&u_3^{\varepsilon}(u_1^{\varepsilon})^{\diamond,2}e_N^{i_1i_2i_3}+3(C_{12}^{\varepsilon,i_1i_2i_3}+\varphi_2^{\varepsilon,i_1i_2i_3}) (u_2^\varepsilon+u_3^\varepsilon),\endaligned$$
 $$\varphi^\varepsilon:=\varphi_1^\varepsilon+\varphi_2^\varepsilon,$$
 where $C^\varepsilon_0,C_{1i}^\varepsilon, \varphi_i^\varepsilon$ are defined in Section 6.  Moreover there exist $\varphi\in C((0,T];\mathbb{R})$ such that for $\rho>0$ small enough $\sup_{t\in[0,T]}t^\rho|\varphi^\varepsilon-\varphi|\rightarrow0$ as $\varepsilon\rightarrow0$.

 Define $$K^{\varepsilon}(t):=\int_0^tP_{t-s}^\varepsilon (u_1^{\varepsilon})^{\diamond,2}ds,\quad \tilde{K}^{\varepsilon}(t):=\int_0^t\tilde{P}_{t-s}^\varepsilon (u_1^{\varepsilon})^{\diamond,2}ds,$$
 and $$K_1^{\varepsilon}(t):=\int_0^tP_{t-s}^\varepsilon [e^{i_1i_2i_3}_N(u_1^{\varepsilon})^{\diamond,2}]ds,\quad \tilde{K}_1^{\varepsilon}(t):=\int_0^t\tilde{P}_{t-s}^\varepsilon [e^{i_1i_2i_3}_N(u_1^{\varepsilon})^{\diamond,2}]ds.$$
 Also define $$\pi_{0,\diamond}(K^{\varepsilon},(u_1^{\varepsilon})^{\diamond,2}):=\pi_{0}(K^{\varepsilon},(u_1^{\varepsilon})^{\diamond,2})-C_{11}^\varepsilon-\varphi_1^\varepsilon,$$
 and $$\pi_{0,\diamond}(K_1^{\varepsilon},e^{i_1i_2i_3}_N(u_1^{\varepsilon})^{\diamond,2}):=\pi_{0}(K_1^{\varepsilon},e^{i_1i_2i_3}_N(u_1^{\varepsilon})^{\diamond,2})-C_{12}^{\varepsilon,i_1i_2i_3}-\varphi_2^{\varepsilon,i_1i_2i_3}.$$
 Now we introduce the following notations:$$\aligned C^\varepsilon_W(T):=&\sup_{t\in[0,T]}(\|u_1^{\varepsilon}\|_{-1/2-\delta/2}+\|(u_1^{\varepsilon})^{\diamond,2}\|_{-1-\delta/2}+\|u_2^{\varepsilon}\|_{1/2-\delta}+\|\pi_{0}( u_2^{\varepsilon},u_1^{\varepsilon})\|_{-\delta}\\&+\|\pi_{0,\diamond}( u_2^{\varepsilon},(u_1^{\varepsilon})^{\diamond,2})\|_{-1/2-\delta}+\|\pi_{0,\diamond}
(K^{\varepsilon},(u_1^{\varepsilon})^{\diamond,2})\|_{-\delta})+\|u_2^{\varepsilon}\|_{C^{1/4-\delta-\kappa/2}_T\mathcal{C}^{\kappa/2}},\endaligned$$
and
$$\aligned E^\varepsilon_W(T):=&\sup_{t\in[0,T]}(\|(u_1^{\varepsilon})^{\diamond,2}e_N^{i_1i_2i_3}\|_{-1-\delta/2}+\|\pi_{0}( u_2^\varepsilon,e^{i_1i_2i_3}_N u^\varepsilon_1)\|_{-\delta}+\|\pi_{0,\diamond}(u_2^{\varepsilon},e^{i_1i_2i_3}_N(u_1^{\varepsilon})^{\diamond,2})\|_{-1/2-\delta}
\\&+\|\pi_0(K^\varepsilon,e_N^{i_1i_2i_3}(u_1^{\varepsilon})^{\diamond,2})\|_{-\delta}+\|\pi_0(K_1^\varepsilon,(u_1^{\varepsilon})^{\diamond,2})\|_{-\delta}+\|\pi_{0,\diamond}(K_1^\varepsilon,e_N^{i_1i_2i_3}(u_1^{\varepsilon})^{\diamond,2})\|_{-\delta})
.\endaligned$$
Here $E^\varepsilon_W$ appears as an error term for the lattice approximations which goes to $0$ in probability (see Section 6.2).

Then Lemma 3.2 and (3.1) implies that for $t\in[0,T]$
$$\|K^\varepsilon(t)\|_{1-\delta}+\|\tilde{K}^\varepsilon(t)\|_{1-\delta}\lesssim t^{\delta/4}C_W^\varepsilon,\quad \|K_1^\varepsilon(t)\|_{1-\delta}+\|\tilde{K}_1^\varepsilon(t)\|_{1-\delta}\lesssim  t^{\delta/4}E_W^\varepsilon.\eqno(4.3)$$

Now we could write the paracontrolled ansatz as follows: $$u_3^{\varepsilon}=-3P_N[\pi_<(u_2^{\varepsilon}+u_3^{\varepsilon},\tilde{K}^{\varepsilon}+\tilde{K}_1^{\varepsilon})]+u^{\varepsilon,\sharp}$$
 with $u^{\varepsilon,\sharp}(t)\in \mathcal{C}^{1+\beta}$. This yields that
 $$\aligned\|u_3^{\varepsilon}(t)\|_{1/2+\delta}\lesssim& \|u_2^{\varepsilon}(t)+u_3^{\varepsilon}(t)\|_{\gamma}(C_W^\varepsilon+E_W^\varepsilon)+\|u^{\varepsilon,\sharp}(t)\|_{1/2+\delta},\endaligned(4.4)$$
 and
$$\aligned\|u_3^{\varepsilon}(t)\|_{1-\delta}\lesssim& \|u_2^{\varepsilon}(t)+u_3^{\varepsilon}(t)\|_{\gamma}(C_W^\varepsilon+E_W^\varepsilon)+\|u^{\varepsilon,\sharp}(t)\|_{1-\delta}.\endaligned\eqno(4.5)$$
Then $u_3^\varepsilon$ solves (4.2) if and only if $u^{\varepsilon,\sharp}$ solves the following equation:
$$\aligned u^{\varepsilon,\sharp}=&P^\varepsilon_t(\textrm{Ext}u_0-u^{\varepsilon}_1(0)) -\int_0^t P_{t-s}^\varepsilon \bigg{[}Q_N[6u_1^{\varepsilon}\diamond                                                                                    u_2^{\varepsilon}u_3^{\varepsilon}+3u_1^{\varepsilon}(u_3^{\varepsilon})^2+3u_1^{\varepsilon}\diamond  (u_2^{\varepsilon})^{2}+(u_2^{\varepsilon}+u_3^{\varepsilon})^3]\\&+3P_N[(\pi_{>}+\pi_{0,\diamond})(u_2^{\varepsilon}+u_3^{\varepsilon},(u_1^{\varepsilon})^{\diamond,2}+e_N^{i_1i_2i_3}(u_1^{\varepsilon})^{\diamond,2})]-9\varphi^\varepsilon u^\varepsilon\bigg{]}ds
\\&-3\int_0^tP^\varepsilon_{t-s} P_N[\pi_<(u_2^{\varepsilon}+u_3^{\varepsilon},(u_1^{\varepsilon})^{\diamond,2}+e_N^{i_1i_2i_3}(u_1^{\varepsilon})^{\diamond,2})]ds+3P_N[\pi_<(u_2^{\varepsilon}+u_3^{\varepsilon}
,\tilde{K}^{\varepsilon}+\tilde{K}_1^{\varepsilon})]
\\:=&P^\varepsilon_t(\textrm{Ext}\Phi_0^\varepsilon-u^{\varepsilon}_1(0))+\int_0^tP^\varepsilon_{t-s}[Q_N\phi^{\varepsilon,\sharp}_1+P_N\phi^{\varepsilon,\sharp}_2+9\varphi^\varepsilon u^\varepsilon]ds+F^{\varepsilon},\endaligned\eqno(4.6)$$
where $F^{\varepsilon}$ represents the last two terms.
\vskip.10in
First we prove the estimate for $\phi^{\varepsilon,\sharp}_1$.
\vskip.10in

\th{Proposition 4.2} For $\phi_1^{\varepsilon,\sharp}$ defined above, the following estimate holds:
$$\aligned&\|Q_N\phi_1^{\varepsilon,\sharp}\|_{-1/2-\delta-2\kappa}
\lesssim C(C^\varepsilon_W,E^\varepsilon_W)(1+\|u^{\varepsilon,\sharp}\|_{1/2+\delta}\|u_3^{\varepsilon}\|_{\gamma}+\|u_3^{\varepsilon}\|^2_{\gamma})
+\|u_3^\varepsilon\|_\gamma^3.\endaligned$$
Here the constant we omit is independent of $N$.

\proof
Since
$$\aligned\Pi_N[u_3^\varepsilon u_2^\varepsilon u_1^\varepsilon]
=P_N[u_3^\varepsilon e^{i_1i_2i_3}_Nu_2^\varepsilon u_1^\varepsilon]
,\endaligned$$
we have for $\delta>2\kappa$
$$\aligned& \|\Pi_N[u_3^\varepsilon u_2^\varepsilon u_1^\varepsilon]\|_{-1/2-\delta/2-2\kappa}\lesssim \|u_3^\varepsilon u_2^\varepsilon u_1^\varepsilon e_N^{i_1i_2i_3}\|_{-1/2-\delta/2-\kappa}
\\\lesssim&(\|e^{i_1i_2i_3}_Nu_1^\varepsilon\|_{-1/2-\delta/2-\kappa}\|u^\varepsilon_2\|_{1/2-\delta}+\|\pi_{0}( u_2^\varepsilon,e^{i_1i_2i_3}_Nu^\varepsilon_1)\|_{-\delta})\|u_3^\varepsilon\|_{1/2+\delta}
 \\\lesssim&(N^{-\kappa/2}\|u_2^\varepsilon\|_{1/2-\delta}\|u^\varepsilon_1\|_{-1/2-\delta/2}+\|\pi_{0}( u_2^\varepsilon,e^{i_1i_2i_3}_Nu^\varepsilon_1)\|_{-\delta}) \|u_3^\varepsilon\|_{1/2+\delta},\endaligned$$
 where in the first and last inequalities we used Lemma 3.1.

Using paraproduct one has
$$\aligned\Pi_N[u_1^{\varepsilon}(u_3^{\varepsilon})^2]=&P_N[u_1^{\varepsilon}e^{i_1i_2i_3}_N(u_3^{\varepsilon})^2]
\\=&P_N[\pi_<(( u_3^{\varepsilon} )^2,e^{i_1i_2i_3}_N u_1^{\varepsilon})+\pi_0(( u_3^{\varepsilon} )^2,e^{i_1i_2i_3}_Nu_1^{\varepsilon})+\pi_>(( u_3^{\varepsilon} )^2,e^{i_1i_2i_3}_Nu_1^{\varepsilon}) ]
\\=&P_N[\pi_<((u_3^{\varepsilon} )^2,e^{i_1i_2i_3}_Nu_1^{\varepsilon})+\pi_0(\pi_0(u_3^{\varepsilon} ,u_3^{\varepsilon}), e^{i_1i_2i_3}_Nu_1^{\varepsilon})\\&+\pi_>(( u_3^{\varepsilon} )^2,e^{i_1i_2i_3}_N u_1^{\varepsilon})
+2C( u_3^{\varepsilon} ,u_3^{\varepsilon},e^{i_1i_2i_3}_Nu_1^{\varepsilon})+2u_3^{\varepsilon}\pi_0(u_3^{\varepsilon},e^{i_1i_2i_3}_Nu_1^{\varepsilon})].\endaligned$$
Here $C( u_3^{\varepsilon} ,u_3^{\varepsilon},e^{i_1i_2i_3}_Nu_1^{\varepsilon})$ is defined in Lemma 2.3.
Then by using Lemmas 2.3 and 3.1 we obtain
$$\aligned \|\Pi_N[(u_3^{\varepsilon})^2u_1^{\varepsilon}]
\|_{-1/2-\delta/2-2\kappa}\lesssim N^{-\kappa/2}\|u_3^{\varepsilon}\|_{1/2+\delta}\|u_3^{\varepsilon}\|_{\gamma}\|u_1^{\varepsilon}\|_{-1/2-\delta/2}.\endaligned\eqno(4.7)$$
Moreover by a similar argument as  (4.7) we have
$$\aligned &\|\Pi_N[(u_2^{\varepsilon})^2\diamond u_1^{\varepsilon}]
\|_{-1/2-\delta/2-2\kappa}
\\\lesssim& N^{-\kappa/2}\|u_2^{\varepsilon}\|_{1/2-\delta}^2\|u_1^{\varepsilon}\|_{-1/2-\delta/2}+\|u_2^{\varepsilon}\|_{1/2-\delta}\|\pi_{0,\diamond}(u_2^{\varepsilon},u_1^{\varepsilon}e^{i_1i_2i_3}_N)\|_{-\delta}.\endaligned$$
Furthermore Lemma 3.1 implies that
$$\|Q_N[(u_2^{\varepsilon}+u_3^{\varepsilon})^3]\|_{\gamma-\kappa}\lesssim \|u_2^{\varepsilon}+u_3^{\varepsilon}\|_{\gamma}^3.$$
The estimate for the terms containing $P_N$ can be obtained similarly. Hence the result follows from (4.4) and the above estimates. $\hfill\Box$
\vskip.10in
Now we consider $\phi^{\varepsilon,\sharp}_2$. To prove an estimate for $\pi_0(u_3^{\varepsilon},(u_1^{\varepsilon})^{\diamond,2}+e_N^{i_1i_2i_3}(u_1^{\varepsilon})^{\diamond,2})$ we have to use paracontrolled ansatz. However, the Fourier cutoff operator $P_N$ does not commute with the paraproduct. Here we follow the technique from [GP15, Lemma 8.16] and prove the following result.
\vskip.10in
\th{Lemma 4.3} Let $\alpha+\beta+\gamma>0,\beta+\gamma<0$, assume that $\alpha\in(0,1)$, and let $\varphi\in\mathcal{C}^\alpha,\psi\in\mathcal{C}^\beta,\chi\in\mathcal{C}^\gamma$. Define the operator for any $f\in\mathcal{C}^\alpha$
$$A_N^1(\psi,\chi)(f):=\pi_0((I-P_N)\pi_<(f,P_N\psi),\chi),$$
and $$A^2_N(\psi,\chi)(f):=\pi_0(P_N\pi_<(f,(P_{3N}-P_N)\psi),\chi).$$
Then for all $\eta<0$
$$\aligned&\|\pi_0(P_N\pi_<(\varphi,P_{3N}\psi),\chi)-\varphi\pi_0(P_N\psi,\chi)\|_{\eta}\\\lesssim & \|\varphi\|_\alpha\|P_N\psi\|_\beta\|\chi\|_\gamma+(\|A^1_N(\psi,\chi)\|_{L(\mathcal{C}^\alpha,\mathcal{C}^\eta)}+\|A^2_N(\psi,\chi)\|_{L(\mathcal{C}^\alpha,\mathcal{C}^\eta)})\|\varphi\|_\alpha.\endaligned$$
Here the constant we omit is independent of $N$.

\proof We have that
$$\pi_0(P_N\pi_<(\varphi,P_{3N}\psi),\chi)=A^2_N(\psi,\chi)(\varphi)+\pi_0(\pi_<(\varphi,P_N\psi),\chi)-A^1_N(\psi,\chi)(\varphi).$$
Thus the result follows from Lemma 2.3.$\hfill\Box$

\vskip.10in
\th{Proposition 4.4}
 For $\phi_2^{\varepsilon,\sharp}$ defined in (4.6), the following estimate holds:
$$\aligned\|P_N\phi_2^{\varepsilon,\sharp}\|_{-1/2-2\delta-\kappa}
\lesssim& C(C^\varepsilon_W,E^\varepsilon_W,A_N,D_N)(1+\|u_3^\varepsilon\|_\gamma+\|u^{\varepsilon,\sharp}\|_{1+\beta}).\endaligned$$
with $$\aligned A_N:=&\|A^1_N(K^\varepsilon+K_1^\varepsilon,(u_1^{\varepsilon})^{\diamond,2}+e_N^{i_1i_2i_3}(u_1^{\varepsilon})^{\diamond,2})\|_{C_TL(\mathcal{C}^{1-\delta},\mathcal{C}^{-1/2-2\delta})}\\&+\|A^2_N(\tilde{K}^\varepsilon+\tilde{K}_1^\varepsilon,(u_1^{\varepsilon})^{\diamond,2}+e_N^{i_1i_2i_3}(u_1^{\varepsilon})^{\diamond,2})\|_{C_TL(\mathcal{C}^{1-\delta},\mathcal{C}^{-1/2-2\delta})}\endaligned$$
and $$\aligned D_N:=&\sup_{t\in[0,T]}(\|\pi_0((I-P_N)\pi_<(u_2^\varepsilon,K^\varepsilon+K_1^\varepsilon),(u_1^{\varepsilon})^{\diamond,2}+e_N^{i_1i_2i_3}(u_1^{\varepsilon})^{\diamond,2})\|_{-\delta}\\&+\|\pi_0(P_N\pi_<(u_2^\varepsilon,(P_{3N}-P_N)(\tilde{K}^\varepsilon+\tilde{K}_1^\varepsilon)),(u_1^{\varepsilon})^{\diamond,2}+e_N^{i_1i_2i_3}(u_1^{\varepsilon})^{\diamond,2})\|_{-\delta}).\endaligned$$
\proof First we consider $\pi_0(u_3^{\varepsilon},(u_1^{\varepsilon})^{\diamond,2}+e_N^{i_1i_2i_3}(u_1^{\varepsilon})^{\diamond,2})$.
By paracontrolled ansatz we obtain
$$\aligned&\pi_0(u_3^{\varepsilon},(u_1^{\varepsilon})^{\diamond,2}+e_N^{i_1i_2i_3}(u_1^{\varepsilon})^{\diamond,2})\\=&-3\pi_0(P_N(\pi_<(u_2^\varepsilon+u_3^{\varepsilon},P_{3N}(\tilde{K}^\varepsilon+\tilde{K}_1^\varepsilon)),
(u_1^{\varepsilon})^{\diamond,2}+e_N^{i_1i_2i_3}(u_1^{\varepsilon})^{\diamond,2})+\pi_0(u^{\varepsilon,\sharp},(u_1^{\varepsilon})^{\diamond,2}+e_N^{i_1i_2i_3}(u_1^{\varepsilon})^{\diamond,2}).\endaligned$$
Here in the equality we used $P_{3N}\tilde{K}^\varepsilon=\tilde{K}^\varepsilon.$
Then by using Lemma 4.3 and $P_{N}\tilde{K}^\varepsilon=K^\varepsilon$ we obtain that for $\beta>\delta/2$
$$\aligned&\|\pi_0(u_3^{\varepsilon},(u_1^{\varepsilon})^{\diamond,2}+e_N^{i_1i_2i_3}(u_1^{\varepsilon})^{\diamond,2})\|_{-1/2-2\delta}
\\\lesssim&\|u_2^\varepsilon+u_3^{\varepsilon}\|_{1/2-\delta}(\|K^\varepsilon+K_1^\varepsilon\|_{1-\delta}\|(u_1^{\varepsilon})^{\diamond,2}+e_N^{i_1i_2i_3}(u_1^{\varepsilon})^{\diamond,2}\|_{-1-\delta/2}+\|\pi_{0,\diamond}(K^\varepsilon+K_1^\varepsilon,(u_1^{\varepsilon})^{\diamond,2}+e_N^{i_1i_2i_3}(u_1^{\varepsilon})^{\diamond,2})\|_{-\delta})
\\&+A_N\|u_3^{\varepsilon}\|_{1-\delta}+D_N
+\|u^{\varepsilon,\sharp}\|_{1+\beta}\|(u_1^{\varepsilon})^{\diamond,2}+e_N^{i_1i_2i_3}(u_1^{\varepsilon})^{\diamond,2}\|_{-1-\delta/2} .\endaligned$$
The estimate for $\pi_>(u_2^{\varepsilon}+u_3^{\varepsilon},(u_1^{\varepsilon})^{\diamond,2}+e_N^{i_1i_2i_3}(u_1^{\varepsilon})^{\diamond,2})$ can be obtained by Lemma 2.2. Thus the result follows by (4.3), (4.4) and (4.5).$\hfill\Box$
\vskip.10in
\th{Remark 4.5} In our case to use the random operator technique, it requires that  $u_3^\varepsilon\in \mathcal{C}^{1/2+\beta+\kappa}$. However the best regularity we can obtain for $u_2^\varepsilon$ is in $\mathcal{C}^{1/2-\delta}$. Thus for the error terms including $u_2^\varepsilon$ we have to calculate it directly which corresponds to $D_N$.
\vskip.10in
\no\textbf{Estimate for $F^\varepsilon$} We now turn to $F^\varepsilon$: Here we divide $F^\varepsilon$ into two parts.
$$\aligned&\|F^\varepsilon\|_{1+\beta}\\\lesssim& \|\int_0^tP_{t-s}^\varepsilon\pi_<(u_2^{\varepsilon}(s)+u_3^{\varepsilon}(s)-(u_2^{\varepsilon}(t)+u_3^{\varepsilon}(t)),(u_1^{\varepsilon})^{\diamond,2}(s)+e^{i_1i_2i_3}_N(u_1^{\varepsilon})^{\diamond,2}(s))ds\|_{1+\beta}
\\&+\|\int_0^tP^\varepsilon_{t-s}\pi_<(u_2^{\varepsilon}(t)+u_3^{\varepsilon}(t),(u_1^{\varepsilon})^{\diamond,2}(s)+e^{i_1i_2i_3}_N(u_1^{\varepsilon})^{\diamond,2}(s))ds-P_N\pi_<(u_2^{\varepsilon}(t)+u_3^{\varepsilon}(t),\tilde{K}^\varepsilon+\tilde{K}_1^\varepsilon )\|_{1+\beta}
\\=&I_1+I_2.\endaligned$$
 The estimate for $I_2$ can be obtained by Lemma 3.3: $$I_2\lesssim t^{\frac{\gamma-\beta-\frac{\delta}{2}-\kappa}{2}}\|u_2^\varepsilon(t)+u_3^\varepsilon(t)\|_{\gamma}(C^\varepsilon_W+E^\varepsilon_W),\eqno(4.8)$$
where by the condition on $\beta$ we have  $\frac{\gamma-\beta-\frac{\delta}{2}-\kappa}{2}>0.$

For $I_1$ we will use the regularity of $u^\varepsilon_2+u^\varepsilon_3$ with respect to time to control it.  Lemmas 2.2 and 3.2 yield that for $5\delta/4+\beta/2+\kappa<1/4$
$$\aligned I_1\lesssim&\int_0^t(t-s)^{-1-\frac{\delta/2+\beta+\kappa}{2}}\|(u_1^\varepsilon)^{\diamond,2}(s)+e^{i_1i_2i_3}_N(u_1^{\varepsilon})^{\diamond,2}(s)\|_{-1-\delta/2}\|u_2^\varepsilon(t)+u^\varepsilon_3(t)-u^\varepsilon_2(s)-u^\varepsilon_3(s)\|_{\kappa/2}ds
\\\lesssim&(C_W^\varepsilon+E_W^\varepsilon)(C_W^\varepsilon+\int_0^t(t-s)^{-1-\frac{\delta/2+\beta+\kappa}{2}}\|u^\varepsilon_3(t)-u^\varepsilon_3(s)\|_{\kappa/2}ds),\endaligned$$
and we note that by Lemmas 3.2 and 3.4 that for $t>s>0$
$$\aligned &\|u^\varepsilon_3(t)-u^\varepsilon_3(s)\|_{\kappa/2}\\\lesssim
&\|(P^\varepsilon_{\frac{t}{2}}-P^\varepsilon_{\frac{s}{2}})(P^\varepsilon_{\frac{t}{2}}+P^\varepsilon_{\frac{s}{2}})(\textrm{Ext}\Phi^\varepsilon_0-u_1^\varepsilon(0))\|_{\frac{\kappa}{2}}+\|\int_0^s(P^\varepsilon_{t-r}-P^\varepsilon_{s-r})G^\varepsilon(r)dr\|_{\frac{\kappa}{2}}+\|\int_s^tP^\varepsilon_{t-r} G^\varepsilon(r)dr\|_{\frac{\kappa}{2}}
\\\lesssim &(t-s)^{b_0}s^{-\frac{z+2\kappa+2b_0}{2}}\|\textrm{Ext}\Phi^\varepsilon_0-u_1^\varepsilon(0)\|_{-z}+(t-s)^{b}\int_0^s(s-r)^{-\frac{1+\delta+\kappa+2b}{2}}\|G^\varepsilon(r)\|_{-1-\delta}dr
\\&+(t-s)^{b_1}(\int_s^t(t-r)^{-\frac{1+\delta+\kappa}{2(1-b_1)}}\| G^\varepsilon(r)\|_{-1-\delta}^{\frac{1}{1-b_1}}dr)^{1-b_1},\endaligned$$
where in the last inequality for the third term we used H\"{o}lder's inequality. Here $\frac{\delta}{2}+\beta+2\kappa<2b_0<2-z-2\kappa$, $\frac{\delta}{2}+\beta+2\kappa<2b<1-\kappa-\delta$, $\frac{1}{2}(\frac{\delta}{2}+\beta+2\kappa)<b_1<[1-\frac{3(\gamma+z+\kappa)}{2}]\wedge \frac{1}{2}(1-\delta-\kappa)$ and $$\aligned G^\varepsilon=&Q_N[3u_1^{\varepsilon}\diamond(u_2^{\varepsilon})^2+6u_1^{\varepsilon}\diamond u_2^{\varepsilon}u_3^{\varepsilon}+3u_1^{\varepsilon}(u_3^{\varepsilon})^2+3(u_1^{\varepsilon})^{\diamond,2}\diamond u_2^{\varepsilon}+3(u_1^{\varepsilon})^{\diamond,2}\diamond u_3^{\varepsilon}+(u_2^{\varepsilon}+u_3^{\varepsilon})^3].\endaligned$$ Moreover, by Propositions 4.2 and 4.3 one has the following estimate
$$\|G^\varepsilon\|_{-1-\delta}\lesssim C(C^\varepsilon_W,E^\varepsilon_W,A_N,D_N)(1+\|u^{\varepsilon,\sharp}\|_{1/2+\delta}\|u_3^{\varepsilon}\|_{\gamma}+\|u_3^{\varepsilon}\|^3_{\gamma}+\|u^{\varepsilon,\sharp}\|_{1+\beta}).\eqno(4.9)$$
Thus we obtain that
$$\aligned I_1
\lesssim &(C_W^\varepsilon+E_W^\varepsilon)\bigg(C_W^\varepsilon+t^{-\frac{\delta/2+\beta+z}{2}-2\kappa}\|\textrm{Ext}\Phi_0^\varepsilon-u_1^\varepsilon(0)\|_{-z}
\\&+\int_0^t\int_r^t(t-s)^{-1-\frac{\delta/2+\beta+2\kappa}{2}+b}(s-r)^{-\frac{1+\delta+\kappa+2b}{2}}ds
\|G^\varepsilon(r)\|_{-1-\delta}dr
\\&+(\int_0^t(t-s)^{-1-\frac{\delta/2+\beta+2\kappa}{2}+b_1}ds)^{b_1}(\int_0^t\int_0^r(t-s)^{-1-\frac{\delta/2+\beta+2\kappa}{2}+b_1}(t-r)^{-\frac{1+\delta+\kappa}{2(1-b_1)}}\\&\| G^\varepsilon(r)\|_{-1-\delta}^{\frac{1}{1-b_1}}dsdr)^{1-b_1}\bigg),
\endaligned$$
where for the last term we used H\"{o}lder's inequality.
Then by changing variable $s=r+(t-r)\sigma$ for the third term and using (4.9) we have
$$\aligned I_1\lesssim &(C_W^\varepsilon+E_W^\varepsilon)t^{-\frac{\delta/2+\beta+z}{2}-2\kappa}\|\textrm{Ext}\Phi_0^\varepsilon-u_1^\varepsilon(0)\|_{-z}+C(C^\varepsilon_W,E^\varepsilon_W,A_N,D_N)
\\&+C(C^\varepsilon_W,E^\varepsilon_W,A_N,D_N)\int_0^t(t-r)^{-\frac{1}{2}-\frac{3\delta/2+\beta+3\kappa}{2}}(\|u^{\varepsilon,\sharp}\|_{1/2+\delta}\|u_3^{\varepsilon}\|_{\gamma}+\|u^{\varepsilon,\sharp}\|_{1+\beta}+\|u_3^{\varepsilon}\|^3_{\gamma})dr
\\&+C(C^\varepsilon_W,E^\varepsilon_W,A_N,D_N)\big{[}\int_0^t(t-r)^{-\frac{1+\delta+\kappa}{2(1-b_1)}}(\|u^{\varepsilon,\sharp}\|_{1+\beta}+\|u^{\varepsilon,\sharp}\|_{1/2+\delta}\|u_3^{\varepsilon}\|_{\gamma}+\|u_3^\varepsilon\|_\gamma^3)^{\frac{1}{1-b_1}}dr\big{]}^{1-b_1}.
\endaligned\eqno(4.10)$$
Combining (4.8) and (4.10) we could control $\|F^\varepsilon\|_{1+\beta}$ by the right hand side of (4.8) and (4.10).
Now we also want to estimate  $\|F^\varepsilon\|_{1/2+\delta}$ and $\|F^\varepsilon\|_\gamma.$ The estimates for these two terms are much easier. We do not need to use Lemma 3.3. We can obtain the following estimates by Lemmas 2.2 and 3.2 directly:
$$\aligned&\|F^\varepsilon\|_{1/2+\delta}\\\lesssim& (C^\varepsilon_W+E^\varepsilon_W)\int_0^t(t-s)^{-\frac{3+3\delta+\kappa}{4}}\|u_3^{\varepsilon}\|_{\gamma}ds
+C(C^\varepsilon_W,E^\varepsilon_W)
+t^{\frac{1-3\delta-\kappa}{4}}\|u_2^\varepsilon(t)+u_3^\varepsilon(t)\|_{\gamma}(C^\varepsilon_W+E^\varepsilon_W),\endaligned\eqno(4.11)$$
and $$\aligned&\|F^\varepsilon\|_{\gamma}\\\lesssim&
C^\varepsilon_W\int_0^t(t-r)^{-\frac{1+\delta/2+\kappa+\gamma}{2}}\|u_3^{\varepsilon}\|_{\gamma}dr+C(C^\varepsilon_W,E^\varepsilon_W)
+t^{\frac{2-\delta-2\gamma-\kappa}{4}}\|u_2^\varepsilon(t)+u_3^\varepsilon(t)\|_{\gamma}(C^\varepsilon_W+E^\varepsilon_W).\endaligned\eqno(4.12)$$

\no \textbf{Uniform estimate of the solution}
\vskip.10in
Now we introduce the following random time: Define for any $L\geq1$
$$\tau^\varepsilon_L:=\inf\{ t\geq0:\|u^\varepsilon(t)\|_{-z}\geq L\}\wedge L\quad\rho_L^\varepsilon:=\inf\{ t\geq0:C^\varepsilon_W(t)+E_W^\varepsilon+A_N+D_N\geq L\}.$$
\vskip.10in
\th{Proposition 4.6}  For any $L,L_1\geq1$, we have
$$\aligned &\sup_{t\in[0,\tau^\varepsilon_L\wedge \rho_{L_1}^\varepsilon]}(t^{\frac{3(\gamma+z+\kappa)}{2}}\|u^{ \varepsilon ,\sharp}\|_{1+\beta}+t^{\frac{1/2+\delta+z+\kappa}{2}}\|u^{ \varepsilon ,\sharp}\|_{1/2+\delta}+t^{\frac{\gamma+z+\kappa}{2}}\|u^{ \varepsilon,\sharp}(t)\|_{\gamma})\lesssim C(L,L_1)
.\endaligned$$
Moreover before $\tau^\varepsilon_L\wedge \rho_{L_1}^\varepsilon$  one has $u_3^{ \varepsilon}$ depends in a  Lipschitz continuous way on the data $\textrm{Ext}\Phi^\varepsilon_0$ and terms in $(C_W^\varepsilon,E_W^\varepsilon,A_N).$
Here we consider $u_3^{ \varepsilon}$ with respect to $$\sup_{t\in[0,\tau^\varepsilon_L\wedge \rho_{L_1}^\varepsilon]}\|u_3^{ \varepsilon}(t)\|_{-z}.$$

\proof By paracontrolled ansatz Lemma 2.2 and (4.3) one then has for $t\in[0,\tau^\varepsilon_L\wedge \rho_{L_1}^\varepsilon]$
$$\aligned\|u_3^{\varepsilon}(t)\|_{\gamma}\lesssim& t^{\delta/4}L_1\|u_2^{\varepsilon}(t)+u_3^{\varepsilon}(t)\|_{\gamma}+\|u^{\varepsilon,\sharp}(t)\|_{\gamma},\endaligned$$
which shows that for $t$ small enough (depending on $L_1$)
$$\aligned\|u_3^{\varepsilon}\|_{\gamma}\lesssim& L_1^2+\|u^{\varepsilon,\sharp}\|_{\gamma}.\endaligned$$
Then  it follows from Propositions 4.2 4.4 and (4.8) (4.10) that for $\frac{3(\gamma+z+\kappa)}{2}<1$
and $t$ small enough (depending on $L_1$)
$$\aligned &t^{\frac{3(\gamma+z+\kappa)}{2}}\|u^{ \varepsilon,\sharp}(t)\|_{1+\beta}\\\lesssim & C\|\textrm{Ext}\Phi_0^\varepsilon-u_1^{\varepsilon}(0)\|_{-z}+t^{\frac{3(\gamma+z+\kappa)}{2}}C\int_0^t(t-r)^{-\frac{3}{4}-\delta-\frac{\beta}{2}-\kappa}(r^{-\frac{3(\gamma+z+\kappa)}{2}}U^\varepsilon +r^{-\frac{(\gamma+z+\kappa)}{2}-\rho}\|u^{\varepsilon,\sharp}\|_\gamma)dr
+C
\\&+Ct^{\frac{3(\gamma+z+\kappa)}{2}}\int_0^t(t-r)^{-\frac{1}{2}-\frac{3\delta/2+\beta+3\kappa}{2}}r^{-\frac{3(\gamma+z+\kappa)}{2}}U^\varepsilon(r)dr
+t^{\frac{3(\gamma+z+\kappa)}{2(1-b_1)}}\int_0^t(t-r)^{-\frac{1+\delta+\kappa}{2(1-b_1)}}r^{-\frac{3(\gamma+z+\kappa)}{2(1-b_1)}}U^\varepsilon(r)^{\frac{1}{1-b_1}}dr
.\endaligned\eqno(4.13)$$
Here and in the following $C=C(L_1)$ and $$U^\varepsilon(r)=r^{3(\gamma+z+\kappa)/2}(\|u^{ \varepsilon,\sharp}(r)\|_{1+\beta}+\|u^{ \varepsilon,\sharp}(r)\|_{1/2+\delta}\|u^{ \varepsilon,\sharp}(r)\|_\gamma
+\|u^{ \varepsilon,\sharp}(r)\|_\gamma^3).$$
A similar argument as (4.13) and using (4.11) (4.12) one also has that for $t$ small enough (depending on $L_1$) and $0<6\kappa<\frac{3}{2}-2z-2\delta-3\gamma$
$$\aligned &t^{\frac{1/2+\delta+z+\kappa}{2}}\|u^{ \varepsilon,\sharp}(t)\|_{1/2+\delta}\\\lesssim & \|\textrm{Ext}\Phi_0^\varepsilon-u_1^{\varepsilon}(0)\|_{-z}+
t^{\frac{1/2+\delta+z+\kappa}{2}}C\int_0^t(t-s)^{-\frac{1+3\delta+2\kappa}{2}}(s^{-\frac{3(\gamma+z+\kappa)}{2}}U^\varepsilon+s^{-\frac{(\gamma+z+\kappa)}{2}-\rho}\|u^{\varepsilon,\sharp}\|_{\gamma})ds
\\&+C+Ct^{\frac{1/2+\delta+z+\kappa}{2}}\int_0^t(t-r)^{-\frac{3+3\delta+\kappa}{4}}r^{-\frac{\gamma+z+\kappa}{2}}r^{\frac{\gamma+z+\kappa}{2}}\|u^{ \varepsilon,\sharp}\|_\gamma dr
,\endaligned\eqno(4.14)$$
and $$\aligned &t^{\frac{\gamma+z+\kappa}{2}}\|u^{ \varepsilon,\sharp}(t)\|_{\gamma}\\\lesssim & \|\textrm{Ext}\Phi_0^\varepsilon-u_1^{\varepsilon}(0)\|_{-z}+
t^{\frac{\gamma+z+\kappa}{2}}C\int_0^t(t-s)^{-\frac{1}{4}-\delta-\frac{\gamma}{2}-\kappa}(s^{-\frac{3(\gamma+z+\kappa)}{2}}U^\varepsilon+s^{-\frac{(\gamma+z+\kappa)}{2}-\rho}\|u^{\varepsilon,\sharp}\|_{\gamma})ds
\\&+C+Ct^{\frac{\gamma+z+\kappa}{2}}\int_0^t(t-r)^{-\frac{1+\delta/2+\gamma+\kappa}{2}}r^{-\frac{(\gamma+z+\kappa)}{2}}r^{\frac{\gamma+z+\kappa}{2}}\|u^{ \varepsilon,\sharp}\|_\gamma dr
.\endaligned\eqno(4.15)$$
Since $\frac{1/2+\delta+\kappa+z}{2}\leq \gamma+z+\kappa$,  combining  (4.13-4.15) we get that by Bihari's inequality there exists some $T_0$ (depending on $L_1$) such that
$$\aligned \sup_{t\in[0,T_0]}(t^{\frac{3(\gamma+z+\kappa)}{2}}\|u^{ \varepsilon ,\sharp}\|_{1+\beta}+t^{\frac{1/2+\delta+z+\kappa}{2}}\|u^{ \varepsilon ,\sharp}\|_{1/2+\delta}+t^{\frac{\gamma+z+\kappa}{2}}\|u^{ \varepsilon,\sharp}(t)\|_{\gamma})\lesssim C(L,L_1)
,\endaligned$$
which combining with Propositions 4.2 and 4.4 implies that
$$\sup_{t\in[0,T_0]}t^{3(\gamma+z+\kappa)/2}\|Q_N\phi_1^{ \varepsilon,\sharp}+P_N\phi_2^{ \varepsilon,\sharp}\|_{-1/2-2\delta-\kappa}\lesssim C(L,L_1).\eqno(4.16)$$
Moreover by paracontrolled ansatz  we also obtain
$$\aligned\|u_3^{ \varepsilon}\|_{-z}\lesssim& t^{\delta/4}\|u_2^{ \varepsilon}+u_3^{ \varepsilon}\|_{-z}L_1+\|u^{ \varepsilon,\sharp}\|_{-z},\endaligned$$
which combining with (4.16) implies that for $t$ small enough and  $t\in [0,T_0]$
$$\aligned\|u_3^{ \varepsilon}(t)\|_{-z}\lesssim&  C+\|u^{ \varepsilon,\sharp}(t)\|_{-z}
\\\lesssim & C+\|\textrm{Ext}\Phi_0^\varepsilon-u_1^\varepsilon(0)\|_{-z}+\|F^\varepsilon(t)\|_{-z}+\int_0^t(t-s)^{-\frac{1/2+2\delta+3\kappa-z}{2}}s^{-\rho}\|u^\varepsilon\|_{-1/2-\delta}ds\\&+\int_0^t(t-s)^{-\frac{1/2+2\delta+3\kappa-z}{2}}s^{-\frac{3(\gamma+z+\kappa)}{2}}s^{\frac{3(\gamma+z+\kappa)}{2}}\|Q_N\phi_1^{ \varepsilon,\sharp}+P_N\phi_2^{ \varepsilon,\sharp}\|_{-1/2-2\delta-\kappa}ds \\\lesssim&C(L,L_1)+t^{\frac{1}{2}-\frac{\delta}{4}-\kappa}\|u_3^{ \varepsilon}(t)\|_{-z}C.\endaligned$$
Here in the last inequality we used  $$\aligned&\|F^{ \varepsilon}(t)\|_{-z}\\\lesssim& C\int_0^t(t-s)^{-\frac{1+\delta/2+\kappa-z}{2}}s^{-\frac{\gamma+\kappa+z}{2}}ds\sup_{s\in[0,t]}s^{\frac{\gamma+\kappa+z}{2}}\|u^{ \varepsilon}_2+u^{ \varepsilon}_3\|_{\gamma}
+t^{\frac{1}{2}-\frac{\delta}{4}-\kappa}\|u_2^{ \varepsilon}(t)+u_3^{ \varepsilon}(t)\|_{-z}C.\endaligned$$
Hence before $T_0$  one has $u_3^{ \varepsilon}$ depends in a  Lipschitz continuous way on the data $\textrm{Ext}\Phi^\varepsilon_0$ and terms in $(C_W^\varepsilon,E_W^\varepsilon,A_N).$
Furthermore we can extend the time to $\tau^\varepsilon_L\wedge \rho_{L_1}^\varepsilon$ as we did in [ZZ14].
$\hfill\Box$

\section{Proof of main result}
In [CC13] it is obtained that the solution to (1.1) can be obtained as limit of solutions $\bar{\Phi}^\varepsilon$ to the following equation:
$$d \bar{\Phi}^{\varepsilon}=\Delta \bar{\Phi}^{\varepsilon}dt+P_N dW-(\bar{\Phi}^{\varepsilon})^3dt+(3 \bar{C}_0^\varepsilon-9\bar{C}_1^\varepsilon)\bar{\Phi}^{\varepsilon}dt,$$
$$\bar{\Phi}^\varepsilon(0)=\Phi_0.$$
Here $\bar{C}_0^\varepsilon$ and $\bar{C}_1^\varepsilon$ are defined in Section 6.1. For $L\geq0$ define $\tau_L:=\inf\{ t\geq0:\|\Phi(t)\|_{-z}\geq L\}\wedge L$ and then $\tau_L$ increases to the explosion time $\tau$. Moreover define $\bar{\tau}^\varepsilon_L:=\inf\{ t\geq0:\|\bar{\Phi}^\varepsilon(t)\|_{-z}\geq L\}\wedge L$ and $\bar{\rho}_L^\varepsilon:=\inf\{ t\geq0:\bar{C}^\varepsilon_W(t)\geq L\}$ with $\bar{C}^\varepsilon_W$ defined similarly as $C^\varepsilon_W$. A similar argument as above implies that
 $$ \sup_{t\in[0,\tau_L\wedge\bar{\rho}_{L_3}^\varepsilon\wedge \bar{\tau}^\varepsilon_{L_4}]}\|\bar{\Phi}^{\varepsilon}(t)-\Phi(t)\|_{-z}\rightarrow^P0.\eqno(5.1)$$
 Here $\Phi$ is the solution to (1.2).
 Define
$$\aligned\delta C_W^{\varepsilon}:=&\sup_{t\in[0,T]}(\|u_1^{\varepsilon}-\bar{u}_1^{\varepsilon}\|_{-1/2-\delta/2}+\|(u_1^{\varepsilon})^{\diamond,2}-(\bar{u}_1^{\varepsilon})^{\diamond,2}\|_{-1-\delta/2}+\|u_2^{\varepsilon}-\bar{u}_2^{\varepsilon}\|_{1/2-\delta}\\&+\|\pi_{0}( u_2^{\varepsilon},u_1^{\varepsilon})-\pi_{0}( \bar{u}_2^{\varepsilon},\bar{u}_1^{\varepsilon})\|_{-\delta}+\|\pi_{0,\diamond}( u_2^{\varepsilon},(u_1^{\varepsilon})^{\diamond,2})-\pi_{0,\diamond}( \bar{u}_2^{\varepsilon},(\bar{u}_1^{\varepsilon})^{\diamond,2})\|_{-1/2-\delta}\\&+\|\pi_{0,\diamond}
(K^{\varepsilon},(u_1^{\varepsilon})^{\diamond,2})-\pi_{0,\diamond}
(\bar{K}^{\varepsilon},(\bar{u}_1^{\varepsilon})^{\diamond,2})\|_{-\delta})+\|u_2^{\varepsilon}-\bar{u}_2^{\varepsilon}\|_{C^{1/4-\delta-\kappa/2}_T\mathcal{C}^{\kappa/2}}.\endaligned$$
Here $\bar{u}^\varepsilon_1,\bar{u}_2^\varepsilon,\bar{u}^\varepsilon_3$ and associated terms are defined similarly as $u^\varepsilon_1,u_2^\varepsilon,u^\varepsilon_3$ and associated terms respectively.  In Section 6 we will prove that $\delta C_W^{\varepsilon}\rightarrow^P0, E_W^{\varepsilon}\rightarrow^P0, A_N\rightarrow^P0$ and $D_N\rightarrow^P0$ as $\varepsilon\rightarrow0$. Then by a similar argument as Section 4 we have
$$ \sup_{t\in[0,\tau_L\wedge \tau_{L_1}^\varepsilon\wedge \rho_{L_2}^\varepsilon\wedge\bar{\rho}_{L_3}^\varepsilon\wedge \bar{\tau}^\varepsilon_{L_4}]}\|u^{\varepsilon}(t)-\bar{\Phi}^\varepsilon(t)\|_{-z}\rightarrow^P0, \quad\varepsilon\rightarrow0.\eqno(5.2)$$
Here $E_W^{\varepsilon},A_N,D_N$ appear as error terms  for lattice approximations.
Then (5.1) and (5.2) implies that
$$ \sup_{t\in[0,\tau_L\wedge \tau_{L_1}^\varepsilon\wedge \rho_{L_2}^\varepsilon\wedge\bar{\rho}_{L_3}^\varepsilon\wedge \bar{\tau}^\varepsilon_{L_4}]}\|u^{\varepsilon}(t)-\Phi(t)\|_{-z}\rightarrow^P0, \quad\varepsilon\rightarrow0.\eqno(5.3)$$
Moreover we have the following estimates:
$$\aligned &P(\sup_{t\in[0,\tau_L]}\|u^\varepsilon-\Phi\|_{-z}>\epsilon)\\\leq& P(\sup_{t\in[0,\tau_L\wedge \tau_{L_1}^\varepsilon\wedge \rho_{L_2}^\varepsilon\wedge\bar{\rho}_{L_3}^\varepsilon\wedge \bar{\tau}^\varepsilon_{L_4}]}\|u^\varepsilon-\Phi\|_{-z}>\epsilon)+P(\tau_L\wedge\rho_{L_2}^\varepsilon\wedge\bar{\rho}_{L_3}^\varepsilon\wedge \bar{\tau}^\varepsilon_{L_4}>\tau_{L_1}^\varepsilon)\\&+P(\tau_L\wedge\bar{\rho}_{L_3}^\varepsilon >\bar{\tau}^\varepsilon_{L_4})+P(\tau_L>\rho_{L_2}^\varepsilon)+P(\tau_L>\bar{\rho}_{L_3}^\varepsilon).\endaligned$$
The first term goes to zero  as $\varepsilon\rightarrow0$ by (5.3). Also  for $L_1>L+\epsilon$
$$P(\tau_L\wedge\rho_{L_2}^\varepsilon\wedge\bar{\rho}_{L_3}^\varepsilon\wedge \bar{\tau}^\varepsilon_{L_4}>\tau_{L_1}^\varepsilon)\leq P(\sup_{t\in[0,\tau_L\wedge \tau_{L_1}^\varepsilon\wedge \rho_{L_2}^\varepsilon\wedge\bar{\rho}_{L_3}^\varepsilon\wedge \bar{\tau}^\varepsilon_{L_4}]}\|u^\varepsilon-\Phi\|_{-z}>\epsilon),$$
which goes to zero as $\varepsilon\rightarrow0$ by (5.3). Moreover for $L_4>L+\epsilon$ we have
$$P(\tau_L\wedge\bar{\rho}_{L_3}^\varepsilon >\bar{\tau}^\varepsilon_{L_4})\leq P(\sup_{t\in[0,\tau_L\wedge\bar{\rho}_{L_3}^\varepsilon\wedge\bar{\tau}^\varepsilon_{L_4}]}\|\bar{\Phi}^\varepsilon-\Phi\|_{-z}>\epsilon)$$
which goes to zero by (5.1).
The last two terms go to zero uniformly over $\varepsilon\in(0,1)$ as $L_2,L_3$ go to $\infty.$
Thus the result follows. $\hfill\Box$

\section{Stochastic convergence}

To simplify the arguments below, we assume that $\mathcal{F}{W}(0)=0$ and restrict ourselves to the flow of $\int_{\mathbb{T}^3} u(x)dx=0$. We follow the notations from [GP15, Section 9].
 We represent the white noise in terms of its spatial Fourier transform. More precisely, let $E=\mathbb{Z}^3\backslash\{0\}$ and let $W(s,k)=\langle W(s),e_k\rangle$ for $e_k(x)=2^{-3/2}e^{\imath\pi x\cdot k},x\in\mathbb{T}^3$.
Then $$u_1^{\varepsilon}(t,x)=\int_{\mathbb{R}\times E}e_k(x)P^\varepsilon_{t-s}(k)W(d\eta),\quad \bar{u}_1^{\varepsilon}(t,x)=\int_{\mathbb{R}\times E}e_k(x)\bar{P}_{t-s}^\varepsilon(k)W(d\eta),$$ where $\eta_a=(s_a,k_a)$, $s_{-a}=s_a, k_{-a}=-k_a$ and the measure $d\eta_a=ds_adk_a$ is the product of the Lebesgue measure $ds_a$ on $\mathbb{R}$ and of the counting measure $dk_a$ on $E$ and $p_t^\varepsilon(k)=e^{-|k|^2f(\varepsilon k)t}1_{\{t\geq0\}}$, $P_t^\varepsilon(k)=p_t^\varepsilon(k)1_{\{|k|_\infty\leq N\}}, p_t(k)=e^{-|k|^2t}1_{\{t\geq0\}},  \bar{P}_t^\varepsilon(k)=p_t(k)1_{\{|k|_\infty\leq N\}}$. Moreover,
$$\int P_{t-s}^{\varepsilon}(k)P_{\sigma-s}^{\varepsilon}(k)ds=\frac{e^{-|k|^2f(\varepsilon k)|t-\sigma|}1_{\{|k|_\infty\leq N\}}}{2|k|^2f(\varepsilon k)}:=V_{t-\sigma}^\varepsilon(k),\eqno(6.1)$$
and
$$\int \bar{P}^\varepsilon_{t-s}(k)\bar{P}^\varepsilon_{\sigma-s}(k)ds=\frac{e^{-|k|^2|t-\sigma|}1_{\{|k|_\infty\leq N\}}}{2|k|^2}:=\bar{V}_{t-\sigma}^\varepsilon(k).\eqno(6.2)$$
In this section we will prove that
 $\delta C_W^\varepsilon\rightarrow0, E_W^\varepsilon\rightarrow0, A_N\rightarrow0, D_N\rightarrow0$ in probability as $\varepsilon\rightarrow0$.

 Now we introduce the following notations: $k_{[1...n]}=\sum_{i=1}^nk_i$, $\tilde{k}^{i_1i_2i_3}=(k^j-i_j(2N+1))_{j=1,2,3}$ for $i_j=1,0,-1$ and $\sum_{j=1}^3i_j^2\neq0$. In the following we always omit the superscript of $\tilde{k}$ if there's no confusion. Denote by $$\int_{(\mathbb{R}\times E)^n}f(\eta_{1...n})W(d\eta_{1...n})$$
  a generic element of the $n$-th chaos of $W$. By [GP15, Section 9.2] We know that
  $$E[|\int_{(\mathbb{R}\times E)^n}f(\eta_{1...n})W(d\eta_{1...n})|^2]\leq (n!)\int_{(\mathbb{R}\times E)^n}|f(\eta_{1...n})|^2d\eta_{1...n},$$
  such that for bounding the variance of the chaos it is enough to bound the $L^2$ norm of the unsymmetrized kernels.
  To obtain the results we first recall the following lemma from [ZZ14] for our later use:

\vskip.10in
\th{Lemma 6.1} ([ZZ14, Lemma 3.10]) Let $0<l,m<d,l+m-d>0$. Then we have
$$\sum_{k_1,k_2\in \mathbb{Z}^d\backslash\{0\},k_1+k_2=k}\frac{1}{|k_1|^{l}|k_2|^{m}}\lesssim \frac{1}{|k|^{l+m-d}}.$$
\vskip.10in
By a similar argument as the proof of [ZZ14, Lemma 3.11] we have the following results.
\vskip.10in
  \th{Lemma 6.2} For any $0<\kappa<1$, $i\geq0$, $t\geq0$, $k_1, k_2\in E$ we have
  $$|e^{-|k_{[12]}|^2t}\theta(2^{-i}k_{[12]})-e^{-|k_{2}|^2t}\theta(2^{-i}k_{2})|\lesssim |k_1|^\kappa2^{-i\kappa}.$$

\vskip.10in
  \th{Lemma 6.3} For any $0<\kappa<1$, $i\geq0$, $t\geq0$ we have for $k_1, k_2\in E$ with $|k_{[12]}|_\infty\leq N, |k_2|_\infty\leq N$
  $$|e^{-|k_{12}|^2t{f}(\varepsilon k_{[12]})}\theta(2^{-i}k_{[12]})-e^{-|k_{2}|^2t{f}(\varepsilon k_2)}\theta(2^{-i}k_{2})|\lesssim |k_1|^\kappa2^{-i\kappa}.$$
\vskip.10in

Now we prove the following estimate for the approximating operators:
\vskip.10in

\th{Lemma 6.4} For any $0<\kappa<1$ and $t>0, k\in E$, $\varepsilon>0$
$$(1) |p_{t}^\varepsilon(k)-p_t(k)|\lesssim e^{-|k|^2\bar{c}_ft} |\varepsilon k|^{\kappa},\quad |P_{t}^\varepsilon(k)-p_t(k)|\lesssim e^{-|k|^2\bar{c}_ft} |\varepsilon k|^{\kappa};$$
$$(2) |P_{t}^\varepsilon(k)-\bar{P}^\varepsilon_t(k)|\lesssim e^{-|k|^2\bar{c}_ft} |\varepsilon k|^{\kappa},\quad |V_{t}^\varepsilon(k)-\bar{V}^\varepsilon_t(k)|\lesssim \frac{e^{-|k|^2\bar{c}_ft} |\varepsilon k|^{\kappa}}{|k|^2}.$$
Here $\bar{c}_f=c_f\wedge1$, $c_f=\min\{f(x):|x|\leq 1.8\}$.

\proof The results follow from $|f(\varepsilon k)-1|\lesssim |\varepsilon k|^\kappa$ and
$$\aligned&|e^{-|k|^2tf(\varepsilon k)}-e^{-|k|^2t}|\lesssim e^{-|k|^2\bar{c}_ft}(1\wedge t^\kappa|f(\varepsilon k)-1|^\kappa|k|^{2\kappa})\lesssim e^{-|k|^2\bar{c}_ft} |\varepsilon k|^{\kappa}.\endaligned$$$\hfill\Box$
\vskip.10in

We prove the following two lemmas for dealing with the error terms.

\vskip.10in
\th{Lemma 6.5} For every $q\geq0, 0<r<3$, $$\int_E\theta(2^{-q}\tilde{k})^2\frac{1}{|k|^r}dk\lesssim2^{(3-r) q},\quad \int_E\theta(2^{-q}\tilde{\tilde{k}})^2\frac{1}{|k|^r}dk\lesssim2^{(3-r) q}.$$
\proof We consider the first one, the second can be obtained by a similar argument. We have $$\int\theta(2^{-q}\tilde{k})^2\frac{1}{|k|^r}dk\lesssim\int1_{|k|\leq 2^q}\theta(2^{-q}\tilde{k})^2\frac{1}{|k|^r}dk+\int1_{|k|> 2^q}\theta(2^{-q}\tilde{k})^2\frac{1}{|k|^{r}}dk \lesssim2^{(3-r) q}$$
Here  in the last inequality we used that the cardinality of the $k$ with $\theta(2^{-q}\tilde{k})\neq0$ is of order $2^{3q}$.$\hfill\Box$
\vskip.10in

\th{Lemma 6.6} For every $q\geq0, 0<r<3$, $$\int\theta(2^{-q}\tilde{k})^2\frac{1}{|k|^r}dk\lesssim\varepsilon^\kappa2^{(3-r+\kappa) q}.$$
 Here $\kappa>0$ is small enough.

\proof We have $$\int\theta(2^{-q}\tilde{k})^2\frac{1}{|k|^r}dk\lesssim\int1_{|k|\leq N}\theta(2^{-q}\tilde{k})^2\frac{1}{|k|^r}dk+\varepsilon^\kappa\int1_{|k|\geq N}\theta(2^{-q}\tilde{k})^2\frac{1}{|{k}|^{r-\kappa}}dk \lesssim\varepsilon^\kappa2^{(3-r+\kappa) q},$$
where in the last inequality we used that $|k|\leq N\backsimeq |\tilde{k}|\backsimeq 2^q$ and Lemma 6.5. $\hfill\Box$

\subsection{Convergence for renormalisation terms}
In this subsection we prove $\delta C_W^\varepsilon\rightarrow0$ in probability as $\varepsilon\rightarrow0$.
\vskip.10in

\no\textbf{Convergence for $u_1^\varepsilon-\bar{u}_1^\varepsilon$}

In this part we consider the convergence of $u_1^\varepsilon-\bar{u}_1^\varepsilon$.
$$\aligned &E|\Delta_q[u_1^{\varepsilon}(t)-\bar{u}_1^{\varepsilon}(t)]|^2
\\\lesssim &\int_{\mathbb{R}\times E}\theta(2^{-q}k)^2|e_k(P_{t-s}^\varepsilon(k)-\bar{P}_{t-s}^\varepsilon(k))|^2d\eta
\lesssim\int\theta(2^{-q}k)^2(\varepsilon|k|)^\kappa|k|^{-2}dk\lesssim\varepsilon^\kappa2^{q(\kappa+1)}.\endaligned$$
Here $\kappa>0$ is small enough and in the  second inequality we used Lemma 6.4. Similarly by using
$$|1-e^{-|t_2-t_1|f(\varepsilon k)|k|^2}|\lesssim |t_1-t_2|^\kappa |k|^{2\kappa},$$
we get desired estimates for $E|\Delta_q[(u_1^{\varepsilon}(t_2)-\bar{u}_1^{\varepsilon}(t_2))-(u_1^{\varepsilon}(t_1)-\bar{u}_1^{\varepsilon}(t_1))]|$, which combining with Gaussian hypercontractivity  implies that for $p>1$,  $\epsilon>0$ small enough
$$\aligned &E[\|(u_1^{\varepsilon}(t_2)-\bar{u}_1^{\varepsilon}(t_2))-(u_1^{\varepsilon}(t_1)-\bar{u}_1^{\varepsilon}(t_1))\|_{B^{-1/2-\kappa-\epsilon}_{p,p}}^p]
\\\lesssim&\varepsilon^{p\kappa/2}|t_2-t_1|^{\kappa p/4},\endaligned$$
Then by Lemma 2.1 we obtain that for  $\delta>0, p>1$, $u_1^{\varepsilon}-\bar{u}_1^{\varepsilon}\rightarrow 0$ in $L^p(\Omega;C_T\mathcal{C}^{-1/2-\delta/2})$ as $\varepsilon\rightarrow0$.

\vskip.10in
 \no\textbf{Convergence for $u_1^{\varepsilon}\diamond u_1^{\varepsilon}-\bar{u}_1^{\varepsilon}\diamond \bar{u}_1^{\varepsilon}$}

 In this part we consider the convergence of $u_1^{\varepsilon}\diamond u_1^{\varepsilon}$. Recall that $u_1^{\varepsilon}\diamond u_1^{\varepsilon}=u_1^{\varepsilon} u_1^{\varepsilon}-C^{\varepsilon}_0$ and $\bar{u}_1^{\varepsilon}\diamond \bar{u}_1^{\varepsilon}=\bar{u}_1^{\varepsilon} \bar{u}_1^{\varepsilon}-\bar{C}^{\varepsilon}_0$.

 Take $$C^{\varepsilon}_0=2^{-3}\int_E\frac{1_{\{|k|_\infty\leq N\}}}{2|k|^2f(\varepsilon k)}dk,\quad\bar{C}^{\varepsilon}_0=2^{-3}\int\frac{1_{\{|k|_\infty\leq N\}}}{2|k|^2}dk.\eqno(6.3)$$
 Then we have
 $$\aligned &E|\Delta_q[u_1^{\varepsilon}\diamond u_1^{\varepsilon}(t)-\bar{u}_1^{\varepsilon}\diamond \bar{u}_1^{\varepsilon}(t)]|^2
\\\lesssim &\int_{(\mathbb{R}\times E)^2}\theta(2^{-q}k_{[12]})^2|(P_{t-s_1}^\varepsilon(k_1)P_{t-s_2}^\varepsilon(k_2)-\bar{P}_{t-s_1}^\varepsilon(k_1)\bar{P}_{t-s_2}^\varepsilon(k_2))|^2d\eta_{12}
\\\lesssim & \varepsilon^\kappa\int\theta(2^{-q}k_{[12]})^2\frac{|k_1|^{\kappa}+|k_2|^{\kappa}}{|k_1|^2|k_2|^2}dk_{12}
\lesssim\varepsilon^\kappa2^{(\kappa+2)q}.\endaligned$$
Here $\kappa>0$ is small enough and in the second inequality we used Lemma 6.4 and in the last inequality we used Lemma 6.1. Then by Gaussian hypercontractivity and  Lemma 2.1 we obtain that for  $\delta>0, p>1,$ $u_1^{\varepsilon}\diamond u_1^{\varepsilon}-\bar{u}_1^{\varepsilon}\diamond \bar{u}_1^{\varepsilon}\rightarrow 0$ in $L^p(\Omega;C_T\mathcal{C}^{-1-\delta})$ as $\varepsilon\rightarrow0$.
\vskip.10in

\no\textbf{Convergence for $u_2^{\varepsilon}-\bar{u}_2^{\varepsilon}$}

 In this part we consider the convergence of $u_2^{\varepsilon}$. Recall that
 $$\bar{u}_2^{\varepsilon}(t)-u_2^{\varepsilon}(t)=I_t^3-\bar{I}_t^3+J_t^3.$$
 Here
 $$\aligned & I_t^3=2^{-3}\int_{(\mathbb{R}\times E)^3}e_{k_{[123]}}\int_0^tP_{t-\sigma}^\varepsilon(k_{[123]})P_{\sigma-s_1}^\varepsilon(k_1)P_{\sigma-s_2}^\varepsilon(k_2)P_{\sigma-s_3}^\varepsilon(k_3)d\sigma W(d\eta_{123}),\endaligned$$
and $\bar{I}_t^3$ is defined similarly as $I_t^3$ with $P^\varepsilon_{t-\sigma}(k_{[123]})$ replaced by $p_{t-\sigma}(k_{[123]})$ and  with other $ P^\varepsilon$ replaced by $\bar{P}^\varepsilon$ and $J_t^3$ is defined similarly as $I_t^3$ with $e_{k_{[123]}}, k_{[123]}$ replaced by $e_{\tilde{k}_{[123]}}, \tilde{k}_{[123]}$.
By  Lemma 6.4 and a straightforward calculation we obtain that
$$\aligned &E|\Delta_q(I_t^3-\bar{I}_t^3)|^2
\\\lesssim &\int_{(\mathbb{R}\times E)^2}\theta(2^{-q}k_{[123]})^2\bigg|\int_0^t\bigg(P_{t-\sigma}^\varepsilon(k_{[123]})\Pi_{i=1}^3P_{\sigma-s_i}^\varepsilon(k_i)
-p_{t-\sigma}(k_{[123]})
\Pi_{i=1}^3\bar{P}_{\sigma-s_i}^\varepsilon(k_i)\bigg)d\sigma\bigg|^2d\eta_{123}\\
\lesssim &\int\theta(2^{-q}k_{[123]})\frac{\varepsilon^\kappa\sum_{i=1}^3|k_i|^{\kappa}+|k_{[123]}|^\kappa\varepsilon^\kappa}
{|k_1|^{2}|k_2|^{2}|k_3|^{2}[|k_1|^2+|k_2|^2+|k_3|^2]|k_{[123]}|^{2}}dk_{123}\\\lesssim &\int_E\theta(2^{-q}k)\frac{\varepsilon^\kappa}{|k|^{4-\kappa}}dk
\lesssim \varepsilon^\kappa2^{q(-1+\kappa)},\endaligned$$
where we used Young's inequality and Lemma 6.1 in the second inequality.
 Similar calculations also imply that
 $$\aligned E|\Delta_qJ_t^3|^2
\lesssim &\int\theta(2^{-q}\tilde{k}_{[123]})\frac{1_{\{|k_{123}|>N, 2^q\lesssim N\}}}
{|k_1|^{2}|k_2|^{2}|k_3|^{2}[|k_1|^2+|k_2|^2+|k_3|^2]|\tilde{k}_{[123]}|^{2}}dk_{123}\\\lesssim &\int_E\theta(2^{-q}\tilde{k})\frac{\varepsilon^\kappa1_{\{|k|>N, 2^q\lesssim N\}}}{|k|^{2-2\kappa}|\tilde{k}|^{2}}dk
\lesssim \varepsilon^\kappa2^{q(-1+\kappa)},\endaligned$$
where we used Lemma 6.5  in the last inequality.
 By a similar argument as above  we also obtain that  for $\delta>0, p>1$, $u_2^{\varepsilon}-\bar{u}_2^{\varepsilon} \rightarrow0$ in $L^p(\Omega;C_T\mathcal{C}^{1/2-\delta/2})$. Similarly we obtain $u_2^{\varepsilon}-\bar{u}_2^{\varepsilon} \rightarrow0$ in $L^p(\Omega;C^{1/4-\delta-\kappa/2}([0,T],\mathcal{C}^{\kappa/2}))$.

\vskip.10in
 \no\textbf{Convergence for $\pi_{0,\diamond}(K^\varepsilon,(u_1^\varepsilon)^{\diamond,2})-\pi_{0,\diamond}(\bar{K}^\varepsilon,(\bar{u}_1^\varepsilon)^{\diamond,2})$}

In this part we focus on $\pi_0(K^\varepsilon,(u_1^\varepsilon)^{\diamond,2})$ and prove that $\pi_{0,\diamond}(K^\varepsilon,(u_1^\varepsilon)^{\diamond,2})-\pi_{0,\diamond}(\bar{K}^\varepsilon,(\bar{u}_1^\varepsilon)^{\diamond,2})$ in $C_T\mathcal{C}^{-\delta}$ for every $\delta>0$.
Now we have the following identity: for $t\in[0,T]$,
$$ \pi_0(K^\varepsilon,(u_1^\varepsilon)^{\diamond,2})(t)-\pi_0(\bar{K}^\varepsilon,(\bar{u}_1^\varepsilon)^{\diamond,2})(t)=I_t^1+4I_t^2+2I_t^3
-[\bar{I}_t^1+4\bar{I}_t^2+2\bar{I}_t^3].$$
Here
 $$\aligned I_t^1=&2^{-\frac{9}{2}}\int e_{k_{[1234]}}\psi_0(k_{[12]},k_{[34]})\int_0^td\sigma P^\varepsilon_{t-\sigma}(k_{[12]})P^\varepsilon_{\sigma-s_1}(k_1)P^\varepsilon_{\sigma-s_2}(k_2)P^\varepsilon_{t-s_3}(k_3)P_{t-s_4}^\varepsilon(k_4)W(d\eta_{1234}),\\ I_t^2=&2^{-\frac{9}{2}}\int\int e_{k_{[23]}}\psi_0(k_{[12]},k_{3}-k_1)\int_0^td\sigma P^\varepsilon_{t-\sigma}(k_{[12]})P^\varepsilon_{\sigma-s_2}(k_2)P^\varepsilon_{t-s_3}(k_3)V_{t-\sigma}^\varepsilon(k_1)dk_1W(d\eta_{23}),
 \\I_t^3=&2^{-6}\int_{E^2}
\int_0^td\sigma V^\varepsilon_{t-\sigma}(k_1)V^\varepsilon_{t-\sigma}(k_2)P^\varepsilon_{t-\sigma}(k_{[12]})dk_{12},\endaligned$$
and for $i=1,2,3$,  $\bar{I}_t^i$ is defined similarly with $P_{t-\sigma}^\varepsilon(k_{[12]})$ replaced by $p_{t-\sigma}(k_{[12]})$ and other $P^\varepsilon, V^\varepsilon$ replaced by $ \bar{P}^\varepsilon, \bar{V}^\varepsilon$ respectively.
In fact, choose $$C_{11}^\varepsilon=2^{-5}\int\int_{-\infty}^td\sigma V^\varepsilon_{t-\sigma}(k_1)V^\varepsilon_{t-\sigma}(k_2)P^\varepsilon_{t-\sigma}(k_{[12]})dk_{12}\eqno(6.4)$$
and $\bar{C}_{11}^\varepsilon$ is defined with each $P^\varepsilon, V^\varepsilon$ replaced by $ p, \bar{V}^\varepsilon$ respectively. Choose $\varphi_1^\varepsilon(t)=2I_t^3-C_{11}^\varepsilon$ and  $\bar{\varphi}_1^\varepsilon(t)=2\bar{I}_t^3-\bar{C}_{11}^\varepsilon$  and  $\varphi_1(t)=-2^{-7}\int\frac{e^{-t(|k_1|^2+|k_2|^2+|k_{[12]}|^2)}}{|k_1|^2|k_2|^2(|k_1|^2+|k_2|^2+|k_{[12]}|^2)}dk_{12}$. Then we can easily obtain that
$$\aligned \sup_{t\in[0,T]}t^\rho|\varphi_1^\varepsilon-\varphi_1|\lesssim\varepsilon^\kappa,\quad \sup_{t\in[0,T]}t^\rho|\bar{\varphi}_1^\varepsilon-\varphi_1|\lesssim\varepsilon^\kappa ,\endaligned$$
for every $\rho>0, 0<\kappa<2\rho$.

\no\textbf{Term in the second chaos:}
Now we consider $I_t^2$ and by Lemma 6.4 and (6.1), (6.2) we have the following calculations:
$$\aligned &E|\Delta_q(I_t^2-\bar{I}_t^2)|^2
\\\lesssim& \int\psi_0( k_{[12]},k_3-k_1)\psi_0( k_{[24]},k_3-k_4)\theta(2^{-q}k_{[23]})^2\\&\frac{|\varepsilon k_{[12]}|^{\kappa/2}|\varepsilon k_{[24]}|^{\kappa/2}+|\varepsilon k_{1}|^{\kappa/2}|\varepsilon k_{4}|^{\kappa/2}+|\varepsilon k_{2}|^{\kappa}+|
\varepsilon k_{3}|^{\kappa}}{|k_1|^2|k_2|^2|k_3|^2|k_4|^2(|k_1|^{2}+|k_{[12]}|^{2})
(|k_4|^{2}+|k_{[24]}|^{2})}dk_{1234}
\\\lesssim& \varepsilon^\kappa\int\theta(2^{-q}k_{[23]})^2\frac{2^{-2q+2\kappa}}{|k_2|^{2-\kappa}|k_3|^{2}}dk_{23}
\\\lesssim& \varepsilon^\kappa 2^{q3\kappa},\endaligned$$
with $\kappa>0$ small enough. Here we used $|k_{[i2]}|\gtrsim 2^q$ on the support of $\psi_0( k_{[i2]},k_3-k_i)\theta(2^{-q}k_{[23]})$ for $i=1,4$ in the second inequality and Lemma 6.1 in the last inequality.

\no\textbf{Terms in the fourth chaos}: Now for $I_t^1$ by (6.1), (6.2) and Lemma 6.4 we have the following calculations:
$$\aligned &E[|\Delta_q(I_t^1-\bar{I}_t^1)|^2]\\\lesssim&\varepsilon^\kappa\int\theta(2^{-q}k_{[1234]})^2
\frac{\theta(2^{-q}k)^2\psi_0(k_{[12]},k_{[34]})}{|k_1|^{2}|k_2|^{2}|k_3|^2|k_4|^2|k_{[12]}|^{4}}(|k_{[12]}|^{\kappa}+\sum_{i=1}^4|k_i|^{\kappa})dk_{1234}
\\\lesssim& \int\theta(2^{-q}k_{[1234]})^2\psi_0(k_{[12]},k_{[34]})\bigg(\frac{\varepsilon^\kappa}{|k_{[34]}||k_{[12]}|^{5-\kappa}}+\frac{ \varepsilon^\kappa}{|k_{[34]}|^{1-\kappa}|k_{[12]}|^{5}}\bigg)dk_{[12][34]}
\\\lesssim& \int\theta(2^{-q}k)^22^{-q(2+\kappa)}\frac{ \varepsilon^\kappa}{|k|^{1-2\kappa}}dk
\lesssim\varepsilon^\kappa2^{q\kappa},\endaligned$$
where we used Lemma 6.1 in the second inequality and $|k_{[12]}|\gtrsim 2^q $ on the support of $\theta(2^{-q}k_{[1234]})$ $\psi_0(k_{[12]},k_{[34]})$ in the third inequality.
Now we have that for $\kappa>0$ small enough
$$E[|\Delta_q(I_t^1-\bar{I}_t^1)|^2]\lesssim2^{q\kappa}\varepsilon^\kappa.$$
By a similar calculation as above and Gaussian hypercontractivity and Lemma 2.1  we obtain that for $\delta>0$, $p>1$
$$\pi_{0,\diamond}(K^{\varepsilon},(u_1^{\varepsilon})^{\diamond2})
-\pi_{0,\diamond}(\bar{K}^{\varepsilon},(\bar{u}_1^{\varepsilon})^{\diamond2})\rightarrow 0\textrm{ in } L^p(\Omega;C_T\mathcal{C}^{-\delta}).$$

\no\textbf{Convergence for $\pi_0(u_2^\varepsilon, u_1^\varepsilon)-\pi_0(\bar{u}_2^\varepsilon, \bar{u}_1^\varepsilon)$}

 In this part we focus on $\pi_0(u_2^\varepsilon, u_1^\varepsilon)$ and prove that $\pi_0(u_2^\varepsilon, u_1^\varepsilon)-\pi_0(\bar{u}_2^\varepsilon, \bar{u}_1^\varepsilon)\rightarrow 0$ in $C_T\mathcal{C}^{-\delta}$.
Now we have the following identity: for $t\in[0,T]$,
$$ \pi_0(\bar{u}_2^\varepsilon, \bar{u}_1^\varepsilon)(t)-\pi_0(u_2^\varepsilon, u_1^\varepsilon)(t)=I_t^1+3I_t^2-[\bar{I}_t^1+3\bar{I}_t^2]+J_t^1+3J_t^2.$$
Here
 $$\aligned I_t^1=&2^{-\frac{9}{2}}\int e_{k_{[1234]}}\psi_0( k_{[123]},k_4)\int_0^td\sigma P_{t-\sigma}^\varepsilon(k_{[123]})
P_{\sigma-s_1}^\varepsilon(k_{1})P_{\sigma-s_2}^\varepsilon(k_{2})P_{\sigma-s_3}^\varepsilon(k_{3})P_{t-s_4}^\varepsilon(k_{4})W(d\eta_{1234}),\\ I_t^2=&2^{-\frac{9}{2}}\int \int e_{k_{[23]}}\psi_0( k_{[123]},k_1)\int_0^td\sigma P_{t-\sigma}^\varepsilon(k_{[123]})
P_{\sigma-s_2}^\varepsilon(k_{2})P_{\sigma-s_3}^\varepsilon(k_{3})V_{t-\sigma}^\varepsilon(k_{1})dk_1W(d\eta_{23}),\endaligned$$
and for $i=1,2$, $\bar{I}_t^i$ is defined with $P_{t-\sigma}^\varepsilon(k_{[123]})$ replaced by $p_{t-\sigma}^\varepsilon(k_{[123]})$ and other $ P^\varepsilon, V^\varepsilon$ replaced by $ \bar{P}^\varepsilon, \bar{V}^\varepsilon$ respectively and for $i=1,2,$ $J_t^i$ is  defined similar as $I_t^i$ with each $k_{[123]}, e_{k_{[1234]}},   e_{k_{[23]}}$ replaced by $\tilde{k}_{[123]}, e_{\tilde{k}_{[1234]}}, e_{\tilde{k}_{[23]}}$.

\no\textbf{Terms in the second chaos:}
First we consider $I_t^2$ and we have the following calculations:
$$\aligned &E|\Delta_q(I_t^2-\bar{I}_t^2)|^2
\\\lesssim& \int\psi_0(k_{[123]},k_1)\psi_0( k_{[234]},k_4)\theta(2^{-q}k_{[23]})^2\\&\frac{|\varepsilon k_{[123]}|^{\kappa/2}|\varepsilon k_{[234]}|^{\kappa/2}+|\varepsilon k_{1}|^{\kappa/2}|\varepsilon k_{4}|^{\kappa/2}+|\varepsilon k_{2}|^{\kappa}+|
\varepsilon k_{3}|^{\kappa}}{|k_2|^2|k_3|^2|k_1|^2(|k_1|^{2}+|k_{[123]}|^{2})|k_4|^2
(|k_4|^{2}+|k_{[234]}|^{2})}dk_{1234}
\\\lesssim& \varepsilon^\kappa\int2^{-q(2-2\kappa)}\theta(2^{-q}k_{[23]})^2\frac{1}{|k_2|^{2-\kappa}|k_3|^{2-\kappa}}
dk_{23}\lesssim \varepsilon^\kappa 2^{3q\kappa},\endaligned$$
where $\kappa>0$ are small enough. Here we used (6.1, 6.2), Lemma 6.4  in the  first inequality and $|k_{[123]}|\gtrsim 2^q$, $k_{[234]}\gtrsim 2^q$ in the second inequality
and we used Lemma 6.1 in the last inequality.
By a similar calculation as above, we know that
 $$\aligned E|\Delta_qJ_t^2|^2
\lesssim& \int2^{-q(2-2\kappa)}\theta(2^{-q}\tilde{k}_{[23]})^2\frac{\varepsilon^\kappa}{|k_2|^{2-\kappa}|k_3|^{2}}
dk_{23}
 \lesssim\varepsilon^\kappa2^{3\kappa q}.\endaligned$$
Here $\kappa>0$ is small enough and in the first inequality we used $|k_{[123]}|\backsimeq N$ to deduce that $|k_i|\backsimeq N$ for some $i\in \{1,2,3\}$ and in the last inequality we used Lemmas 6.1 and 6.5.

\no\textbf{Terms in the fourth chaos}: Now for $I_t^1$ we have the following calculations:
$$\aligned &E[|\Delta_q(I_t^1-\bar{I}_t^1)|^2]\\\lesssim&\varepsilon^\kappa\int\frac{\theta(2^{-q}k_{[1234]})^2\psi_0( k_{[123]},k_4)(|k_{[123]}|^\kappa+\sum_{i=1}^4|k_i|^{\kappa})}{|k_1|^{2}|k_2|^{2}|k_3|^2|k_4|^2[|k_1|^2+|k_2|^2+|k_3|^2]|k_{[123]}|^{2}}dk_{1234}
\\\lesssim& \int2^{-q(2-\kappa)}\theta(2^{-q}k)^2\frac{ \varepsilon^\kappa}{|k|}dk
\lesssim\varepsilon^\kappa2^{q\kappa},\endaligned$$
where we used (6.1), (6.2) and Lemma 6.4 in the first inequality, Lemma 6.1 in the second inequality and $|k_{[123]}|\gtrsim 2^q$ in the third inequality. For $J_t^1$, using Lemma 6.5 and by a similar argument we also obtain that
 $$\aligned &E|\Delta_qJ_t^1|^2
 \lesssim\varepsilon^\kappa2^{\kappa q}.\endaligned$$
Now by a similar calculation as above, Gaussian hypercontractivity and Lemma 2.1 we have that for $\delta>0$, $p>1$
$$\pi_{0,\diamond}(u_2^{\varepsilon},u_1^{\varepsilon})-\pi_{0,\diamond}(\bar{u}_2^{\varepsilon},\bar{u}_1^{\varepsilon})\rightarrow 0\textrm{ in } L^p(\Omega;C_T\mathcal{C}^{-\delta}).$$

\no\textbf{Convergence for $\pi_{0,\diamond}(u_2^\varepsilon, (u_1^\varepsilon)^{\diamond2})-\pi_{0,\diamond}(\bar{u}_2^\varepsilon, (\bar{u}_1^\varepsilon)^{\diamond2})$}
 In this part we focus on $\pi_{0,\diamond}(u_2^\varepsilon, (u_1^\varepsilon)^{\diamond2})$ and prove that $\pi_{0,\diamond}(u_2^\varepsilon, (u_1^\varepsilon)^{\diamond2})-\pi_{0,\diamond}(\bar{u}_2^\varepsilon, (\bar{u}_1^\varepsilon)^{\diamond2})\rightarrow 0$ in $C_T\mathcal{C}^{-1/2-\delta/2}$.
Now we have the following identity: for $t\in[0,T]$,
$$ \pi_{0,\diamond}(\bar{u}_2^\varepsilon, (\bar{u}_1^\varepsilon)^{\diamond2})-\pi_{0,\diamond}(u_2^\varepsilon, (u_1^\varepsilon)^{\diamond2})=I_t^1+6I_t^2+6I_t^3-[\bar{I}_t^1+6\bar{I}_t^2+6\bar{I}_t^3]+J_t^1+6J_t^2+6J_t^3.$$
Here
 $$\aligned I_t^1=&2^{-6}\int e_{k_{[12345]}}\psi_0( k_{[123]},k_{[45]})\int_0^td\sigma P^\varepsilon_{t-\sigma}(k_{[123]})\Pi_{i=1}^3P^\varepsilon_{\sigma-s_i}(k_{i})
 \Pi_{i=4}^5P^\varepsilon_{t-s_i}(k_{i})W(d\eta_{12345}),\\ I_t^2=&2^{-6}\int e_{k_{[234]}}\psi_0( k_{[123]},k_{4}-k_1)\int_0^td\sigma P^\varepsilon_{t-\sigma}(k_{[123]})\Pi_{i=2}^3P^\varepsilon_{\sigma-s_i}(k_{i})
 P^\varepsilon_{t-s_4}(k_{4})V_{t-\sigma}^\varepsilon(k_1)W(d\eta_{234}),\\I_t^3=&2^{-6}\int e_{k_3}\psi_0 (k_{[123]},k_{[12]})\int_0^td\sigma P_{\sigma-s_3}^{\varepsilon}(k_3)V^\varepsilon_{t-\sigma}(k_1)V^\varepsilon_{t-\sigma}(k_2)P_{t-\sigma}^\varepsilon(k_{[123]})W(d\eta_3),\endaligned$$
 and for $i=1,2,3$, $\bar{I}_t^i$ is defined similarly with $P_{t-\sigma}^\varepsilon(k_{[123]})$ replaced by $p_{t-\sigma}(k_{[123]})$ and other $ P^\varepsilon, V^\varepsilon$ replaced by $ \bar{P}^\varepsilon, \bar{V}^\varepsilon$ respectively and for $i=1,2,3$,  $J_t^i$ is defined similar as $I_t^i$ with each $k_{[123]}, e_{k_{[12345]}},$ $e_{k_{[234]}}, e_{k_{3}}$ replaced by $\tilde{k}_{[123]}, e_{\tilde{k}_{[12345]}}, e_{\tilde{k}_{[234]}}, e_{\tilde{k}_{3}}$.

 We consider the following term first: $$I_t^3-\bar{I}_t^3-[\tilde{I}_t^3-\tilde{\bar{I}}_t^3]+\tilde{I}_t^3-\tilde{\bar{I}}_t^3-C^\varepsilon(t) u_1^{\varepsilon}(t)-\bar{C}^\varepsilon(t) \bar{u}_1^{\varepsilon}(t),$$
 where $\tilde{I}_t^3, \tilde{\bar{I}}_t^3$ are defined similarly as $I_t^3, {\bar{I}}_t^3$ with $P_{\sigma-s_3}^{\varepsilon}(k_3), \bar{P}_{\sigma-s_3}^\varepsilon(k_3) $ replaced by $P_{t-s_3}^{\varepsilon}(k_3), \bar{P}^\varepsilon_{t-s_3}(k_3)$, respectively and $C(t)^\varepsilon=\frac{1}{2}[C_{11}^\varepsilon+\varphi_1^\varepsilon(t)], \bar{C}(t)^\varepsilon=\frac{1}{2}[\bar{C}_{11}^\varepsilon+\bar{\varphi}_1^\varepsilon(t)]$.

Since $\int|P_{t-s_3}^\varepsilon(k_3)-P_{\sigma-s_3}^\varepsilon(k_3)|^2ds_3\lesssim\frac{(t-\sigma)^{\kappa/2}}{|k_3|^{2-\kappa}}$ and $$\int|P_{t-s_3}^\varepsilon(k_3)-P_{\sigma-s_3}^\varepsilon(k_3)-[\bar{P}_{t-s_3}^\varepsilon(k_3)-\bar{P}_{\sigma-s_3}^\varepsilon(k_3)]|^2ds_3\lesssim\frac{(t-\sigma)^{\kappa/2}\wedge\varepsilon^\kappa}{|k_3|^{2-\kappa}},$$ by a straightforward calculation we obtain that for $\kappa>0$ small enough
$$\aligned& E[|\Delta_q(I_t^3-\bar{I}_t^3-[\tilde{I}_t^3-\tilde{\bar{I}}_t^3])|^2]
\\\lesssim &\int\theta(2^{-q}k_3)^2\bigg{[}\frac{1}{|k_3|^{2-2\kappa}}\bigg(\int_0^t\int
\varepsilon^{\kappa/2}(|k_{[123]}|^{\kappa/2}+|k_{2}|^{\kappa/2}+|k_{1}|^{\kappa/2})\frac{e^{-(|k_{[123]}|^2+|k_1|^2+|k_2|^2)\bar{c}_f(t-\sigma)}}{|k_1|^2|k_2|^2}
\\&(t-\sigma)^{\kappa/2}dk_{12}d\sigma\bigg)^2
+\frac{\varepsilon^{\kappa}}{|k_3|^{2-2\kappa}}\bigg(\int_0^t\int\frac{e^{-(|k_{[123]}|^2+|k_1|^2+|k_2|^2)(t-\sigma)}}{|k_1|^2|k_2|^2}
(t-\sigma)^{\kappa/4}dk_{12}d\sigma\bigg)^2\bigg{]}dk_3
\\\lesssim &\varepsilon^\kappa 2^{q(1+3\kappa)}.\endaligned$$
Here in the last inequality we used that $\sup_{a\geq0} a^re^{-a}\leq C$ for $r\geq0$ and Lemma 6.1.
Moreover, by Lemmas 6.2 and 6.3 we obtain that
$$\aligned &E[|\Delta_q(\tilde{I}_t^3-\tilde{\bar{I}}_t^3-u_1^{\varepsilon}(t) C^{\varepsilon}(t)+\bar{u}_1^{\varepsilon}(t) \bar{C}^{\varepsilon}(t))|^2]\\\lesssim& \int\frac{1}{|k_3|^2}\theta(2^{-q}k_3)\bigg{(}\int\int_0^t|k_{[12]}|^{-\kappa}|k_3|^\kappa\endaligned$$
$$\aligned& (\varepsilon^{\kappa/2}|k_2|^{\kappa/2}+\varepsilon^{\kappa/2}|k_{1}|^{\kappa/2}+\varepsilon^{\kappa/2}|k_3|^{\kappa/2})
\frac{e^{-|k_1|^2(t-\sigma)\bar{c}_f-|k_2|^2(t-\sigma)\bar{c}_f}}{|k_1|^2|k_2|^2}dk_{12}d\sigma\bigg{)}^2dk_3\\&+\int\frac{\varepsilon^\kappa|k_3|^\kappa}{|k|^2}\theta(2^{-q}k_3)^2
(\int\int_0^t\frac{e^{-|k_2|^2(t-\sigma)-|k_1|^2(t-\sigma)}}{|k_1|^2|k_2|^2}|k_3|^{\kappa}|k_{[12]}|^{-\kappa}dk_{12}d\sigma)^2dk_3
\\\lesssim&\varepsilon^{\kappa}\int\theta(2^{-q}k_3)\frac{1}{|k_3|^{2-3\kappa}}dk_3
\lesssim\varepsilon^{\kappa}2^{q(1+3\kappa)}.\endaligned$$

For $J_t^3$ we  have
$$\aligned E[|\Delta_qJ_t^3|^2]
\lesssim& \int\frac{1}{|k_3|^2}\theta(2^{-q}\tilde{k}_3)\bigg{(}\int
\frac{1_{|k_1|\leq N, |k_2|\leq N}}{|k_1|^2|k_2|^2(|k_1|^2+|k_2|^2+|\tilde{k}_{[12]}|^2)}dk_{12}\bigg{)}^2dk_3
\lesssim\varepsilon^{\kappa}2^{q(1+3\kappa)}.\endaligned$$
Here we used $2^q\backsimeq N$ in the last inequality.

\no\textbf{Terms in the third chaos:}
Now we focus on the bounds for $I_t^2$. We obtain the following inequalities:
$$\aligned &E|\Delta_q (I_t^2-\bar{I}_t^2)|^2\\\lesssim &\int\theta(2^{-q} k_{[234]})\psi_0( k_{[123]},k_4-k_1)\psi_0( k_{[235]},k_4-k_5)\\&\Pi_{i=1}^5\frac{1}{|k_i|^2}\frac{|k_{[123]}|^{\kappa/2}
|k_{[235]}|^{\kappa/2}\varepsilon^\kappa+\sum_{i=2}^4(\varepsilon|k_i|)^{\kappa}}{(|k_1|^{2}+|k_{[123]}|^{2}
+|k_2|^{2})(|k_5|^{2}+|k_{[235]}|^{2})}dk_{12345}\\\lesssim &\int2^{-q(1-\kappa)}\frac{\varepsilon^\kappa\theta(2^{-q} k_{[234]})}{|k_2|^{3-2\kappa}|k_3|^{2}|k_4|^{2}}
dk_{234}
\lesssim \varepsilon^\kappa2^{q(1+4\kappa)},\endaligned$$
where we used Lemma 6.1 in the last inequality. For $J_t^2$ by a similar calculation as above and using the fact that $|k_{[235]}|\backsimeq N\gtrsim |\tilde{k}_{[235]}|$, we know that
 $$\aligned &E|\Delta_qJ_t^2|^2
\lesssim \int2^{-q(1-\kappa)}\theta(2^{-q}\tilde{k}_{[234]})^2\frac{1}{|k_2|^{3}|k_3|^{2}|k_4|^2}
dk_{234}
 \lesssim\varepsilon^\kappa2^{(1+2\kappa) q}.\endaligned$$
Here $\kappa>0$ is small enough and in the last inequality we used Lemmas 6.1 and 6.6.

\no\textbf{Term in the fifth chaos:}
 Now we focus on the bounds for $I_t^1$. We obtain the following inequalities:
$$\aligned &E|\Delta_q (I_t^1-\tilde{I}_t^1)|^2\\\lesssim&\int\theta(2^{-q}{k}_{[12345]})^2\psi_0( k_{[123]},k_{[45]})^2\Pi_{i=1}^5\frac{1}{|k_i|^2}\frac{(\sum_{i=1}^5|\varepsilon k_i|^\kappa+|\varepsilon k_{[123]}|^\kappa)}{|k_{[123]}|^{2}(|k_1|^{2}+|k_2|^{2}+|k_{[123]}|^2)}dk_{12345}\\\lesssim &\varepsilon^\kappa2^{q(1+2\kappa)}.\endaligned$$
For $J_t^1$ by similar calculations for  $I_t^1$ and using the fact that $|k_{[123]}|\backsimeq N\gtrsim |\tilde{k}_{[123]}|$ we also obtain that
 $$\aligned &E|\Delta_qJ_t^1|^2
 \lesssim\varepsilon^\kappa2^{q(1+2\kappa)}.\endaligned$$
 By a similar calculation as above  we also obtain that there exist $\kappa,\epsilon,\gamma>0$ small enough such that
$$\aligned &E[|\Delta_q(\pi_{0,\diamond}(u_2^\varepsilon, (u_1^\varepsilon)^{\diamond2})(t_1)-\pi_{0,\diamond}(u_2^\varepsilon, (u_1^\varepsilon)^{\diamond2})(t_2)-\pi_{0,\diamond}(\bar{u}_2^\varepsilon, (\bar{u}_1^\varepsilon)^{\diamond2})(t_1)+\pi_{0,\diamond}(\bar{u}_2^\varepsilon, (\bar{u}_1^\varepsilon)^{\diamond2})(t_2))|^2]\\\lesssim& \varepsilon^\gamma|t_1-t_2|^{\kappa}2^{q(1+\epsilon)},\endaligned$$
which by Gaussian hypercontractivity and Lemma 2.1 implies that for every $\delta>0, p>1$, $\pi_{0,\diamond}(u_2^\varepsilon, (u_1^\varepsilon)^{\diamond2})-\pi_{0,\diamond}(\bar{u}_2^\varepsilon, (\bar{u}_1^\varepsilon)^{\diamond2})\rightarrow0$ in $L^p(\Omega;C_T\mathcal{C}^{-1/2-\delta/2})$ .

\subsection{Convergence of the error terms}
In this subsection we prove $E_W^\varepsilon\rightarrow^P0$ as $\varepsilon\rightarrow0$.

\no\textbf{Convergence for $\pi_{0}(K^\varepsilon,e_N^{i_1i_2i_3}(u_1^\varepsilon)^{\diamond,2})$}

Now we have the following identity: for $t\in[0,T]$,
$$ \pi_{0}(K^\varepsilon,e_N^{i_1i_2i_3}(u_1^\varepsilon)^{\diamond,2})=I_t^1+4I_t^2+2I_t^3,$$
 $$\aligned I_t^1=&2^{-6}\int e_{\tilde{k}_{[1234]}}\psi_0(k_{[12]},\tilde{k}_{[34]})\int_0^td\sigma P^\varepsilon_{t-\sigma}(k_{[12]})P^\varepsilon_{\sigma-s_1}(k_1)P^\varepsilon_{\sigma-s_2}(k_2)P^\varepsilon_{t-s_3}(k_3)P_{t-s_4}^\varepsilon(k_4)W(d\eta_{1234}),
 \\ I_t^2=&2^{-6}\int\int e_{\tilde{k}_{[23]}}\psi_0({k}_{[12]},\tilde{k}_{3}-k_1)\int_0^td\sigma P^\varepsilon_{t-\sigma}(k_{[12]})P^\varepsilon_{\sigma-s_2}(k_2)P^\varepsilon_{t-s_3}(k_3)V_{t-\sigma}^\varepsilon(k_1)dk_1W(d\eta_{23})
 \\I_t^3=&2^{-6}\int e_N^{i_1i_2i_3}\psi_0({k}_{[12]},-\tilde{k}_{[12]})\int_0^td\sigma P^\varepsilon_{t-\sigma}(k_{[12]})V^\varepsilon_{t-\sigma}(k_1)V_{t-\sigma}^\varepsilon(k_2)dk_{12}.\endaligned$$
\textbf{ Term in the $0$-th chaos:}
We have
$$\aligned E[|\Delta_qI_t^3|^2]\lesssim&  \bigg{(}\int
\frac{1_{|k_{[12]}|\backsimeq N\backsimeq 2^q}\psi_0 ({k}_{[12]},\tilde{k}_{[12]})}{|k_{[12]}|^3}dk_{[12]}\bigg{)}^2
\lesssim\varepsilon^{\kappa}2^{q(3\kappa)}.\endaligned$$

\no\textbf{Term in the second chaos:}
Now we consider $I_t^2$ and we have the following calculations:
$$\aligned E|\Delta_qI_t^2|^2
\lesssim&\int\psi_0( k_{[12]},\tilde{k}_3-k_1)\psi_0( k_{[24]},\tilde{k}_3-k_4)\theta(2^{-q}\tilde{k}_{[23]})^2\\&\frac{1}{|k_2|^2|k_3|^2|k_1|^2(|k_1|^{2}+|k_{12}|^{2})|k_4|^2
(|k_4|^{2}+|k_{[24]}|^{2})}dk_{1234}
\\\lesssim& \int2^{(-2+\kappa)q}\theta(2^{-q}\tilde{k}_{[23]})^2\frac{1}{|k_2|^{2}|k_3|^{2}}dk_{23}\lesssim\varepsilon^\kappa2^{q2\kappa},\endaligned$$
where $\kappa>0$ is small enough. Here we used $|k_{[12]}|\gtrsim 2^q$ in the second inequality and used
Lemmas 6.1, 6.6 in the  third inequality.

\no\textbf{Term in the fourth chaos}: Now for $I_t^1$ we have the following calculations:
$$\aligned E[|\Delta_qI_t^1|^2]\lesssim&\int\theta(2^{-q}\tilde{k}_{[1234]})^2\psi_0( {k}_{[12]},\tilde{k}_{[34]})
\frac{1}{|k_{[34]}||{k}_{[12]}|^{5-\kappa}}dk_{[12][34]}\\\lesssim&2^{-2q}\int\theta(2^{-q}\tilde{k}_{[1234]})^2
\frac{1}{|k_{[34]}||{k}_{[12]}|^{3-\kappa}}dk_{[12][34]}
\lesssim\varepsilon^\kappa2^{q\kappa},\endaligned$$
where we used Lemmas 6.1, 6.6 in the last inequality.
By a similar calculation as above, Gaussian hypercontractivity and Lemma 2.1  we obtain that for $\delta>0$ small enough, $p>1$
$$\pi_{0}(K^\varepsilon,e_N^{i_1i_2i_3}(u_1^\varepsilon)^{\diamond,2})\rightarrow 0\textrm{ in } L^p(\Omega;C_T\mathcal{C}^{-\delta}).$$

\no\textbf{Convergence of $\pi_{0}(K_1^{\varepsilon}, (u_1^{\varepsilon})^{\diamond,2})$}
Now we have the following identity: for $t\in[0,T]$,
$$ \pi_{0}(K_1^{\varepsilon},(u_1^{\varepsilon})^{\diamond,2})=I_t^1+4I_t^2+2I_t^3,$$
 $$\aligned I_t^1=&2^{-6}\int e_{{\tilde{k}}_{[1234]}}\psi_0(\tilde{k}_{[12]},k_{[34]})\int_0^td\sigma P^\varepsilon_{t-\sigma}(\tilde{k}_{[12]})P^\varepsilon_{\sigma-s_1}(k_1)P^\varepsilon_{\sigma-s_2}(k_2)P^\varepsilon_{t-s_3}(k_3)P_{t-s_4}^\varepsilon(k_4)W(d\eta_{1234}),\\ I_t^2=&2^{-6}\int\int e_{\tilde{k}_{[23]}}\psi_0(\tilde{k}_{[12]},k_{3}-k_1)\int_0^td\sigma P^\varepsilon_{t-\sigma}(\tilde{k}_{[12]})P^\varepsilon_{\sigma-s_2}(k_2)P^\varepsilon_{t-s_3}(k_3)V_{t-\sigma}^\varepsilon(k_1)dk_1W(d\eta_{23})\\I_t^3=&2^{-6}\int e_N^{i_1i_2i_3}\psi_0({\tilde{k}}_{[12]},-{k}_{[12]})\int_0^td\sigma P^\varepsilon_{t-\sigma}(\tilde{k}_{[12]})V^\varepsilon_{t-\sigma}(k_1)V_{t-\sigma}^\varepsilon(k_2)dk_{12}.\endaligned$$
$I_t^3, I_t^2$ can be estimated similarly as that for $\pi_0(K^\varepsilon,e_N^{i_1i_2i_3}(u_1^\varepsilon)^{\diamond,2})$ and we only consider
\textbf{Terms in the fourth chaos}: Now for $I_t^1$ we have the following calculations:
$$\aligned E|\Delta_qI_t^1|^2
\lesssim&\int\psi_0( \tilde{k}_{[12]},k_{[34]})\theta(2^{-q}\tilde{k}_{[1234]})^2\frac{1}{|k_2|^2|k_3|^2|k_1|^2(|k_1|^{2}+|k_{2}|^{2})|k_4|^2
|\tilde{k}_{[12]}|^{2}}dk_{1234}
\\\lesssim& \int2^{-2q}\theta(2^{-q}\tilde{k}_{[1234]})^2\frac{1}{|k_{[12]}|^{3-\kappa}|k_{[34]}|}dk_{[12][34]}\lesssim\varepsilon^\kappa2^{q2\kappa},\endaligned$$
where we used Lemmas 6.1, 6.6 in the last inequality.
By a similar calculation as above, Gaussian hypercontractivity and Lemma 2.1  we obtain that for $\delta>0$, $p>1$
$$\pi_{0,\diamond}(K_1^{\varepsilon}, (u_1^{\varepsilon})^{\diamond,2})\rightarrow 0\textrm{ in } L^p(\Omega;C_T\mathcal{C}^{-\delta}).$$
\textbf{Convergence of $\pi_{0,\diamond}( K_1^{\varepsilon},e_N^{i_1i_2i_3} (u_1^{\varepsilon})^{\diamond,2})$}

We have
$$\pi_{0}( K_1^{\varepsilon},e_N^{i_1i_2i_3} (u_1^{\varepsilon})^{\diamond,2})=I_t^1+4I_t^2+2I_t^3.$$
Here $I_t^i, i=1,2$ is defined similarly as that for $\pi_{0}( K^{\varepsilon},e_N^{i_1i_2i_3} (u_1^{\varepsilon})^{\diamond,2})$ with $k_{[12]}$, $e_{\tilde{k}_{[1234]}}$ and $e_{\tilde{k}_{[23]}}$ replaced by  $\tilde{k}_{[12]}$,  $e_{\tilde{\tilde{k}}_{[1234]}}$ and $e_{\tilde{\tilde{k}}_{[23]}}$, respectively and
$$\aligned I_t^3=&2^{-3}\int e_N^{i_1i_2i_3}e_N^{i_1'i_2'i_3'}\psi_0({\tilde{k}}_{[12]}^{i'_1i'_2i'_3},{\widetilde{-k_{[12]}}}^{i_1i_2i_3})\int_0^td\sigma P^\varepsilon_{t-\sigma}(\tilde{k}_{[12]}^{i'_1i'_2i'_3})V^\varepsilon_{t-\sigma}(k_1)V_{t-\sigma}^\varepsilon(k_2)dk_{12},\endaligned$$
 for $i_j, i_j'\in\{-1,0,1\}$ for $j=1,2,3$ with $\sum_{j}i_j^2\neq 0, \sum_j(i_j')^2\neq0$. Choose
 $$C^{\varepsilon,i_1i_2i_3}_{12}=2^{-5}\int_{-\infty}^td\sigma P^\varepsilon_{t-\sigma}({\widetilde{-k_{[12]}}}^{i_1i_2i_3})V^\varepsilon_{t-\sigma}(k_1)V_{t-\sigma}^\varepsilon(k_2)dk_{12},$$
 and  $\varphi^{\varepsilon,i_1i_2i_3}_{2}(t)=-2^{-5}\int_{-\infty}^0d\sigma P^\varepsilon_{t-\sigma}({\widetilde{-k_{[12]}}}^{i_1i_2i_3})V^\varepsilon_{t-\sigma}(k_1)V_{t-\sigma}^\varepsilon(k_2)dk_{12},$ we could easily obtain that
 $$|C^{\varepsilon,i_1i_2i_3}_{12}|\backsimeq1,\quad \sup_{t\in[0,T]}t^\rho|\varphi^{\varepsilon,i_1i_2i_3}_{2}(t)|\lesssim\varepsilon^\kappa,$$
 for every $\rho>\kappa/2>0$.
For the terms in $2I_t^3-C^{\varepsilon,i_1i_2i_3}_{12}-\varphi^{\varepsilon,i_1i_2i_3}_{2}$, we know that $e_N^{i_1i_2i_3}e_N^{i_1'i_2'i_3'}\neq 1$ and we could easily obtain
$$\aligned E[|\Delta_q(2I_t^3-C^\varepsilon_{12}-\varphi^\varepsilon_{2})|^2]
\lesssim\varepsilon^{\kappa}2^{q(3\kappa)}.\endaligned$$

\no\textbf{Term in the second chaos:}
Now we consider $I_t^2$ and we have the following calculations:
$$\aligned E|\Delta_qI_t^2|^2
\lesssim&\int\psi_0( \tilde{k}_{[12]},\tilde{k}_3-k_1)\psi_0( \tilde{k}_{[24]},\tilde{k}_3-k_4)\theta(2^{-q}\tilde{k}_{[23]})^21_{|k_{[12]}|>N, |k_{[24]}|>N}\\&\frac{1}{|k_2|^2|k_3|^2|k_1|^2(|k_1|^{2}+|\tilde{k}_{12}|^{2})|k_4|^2
(|k_4|^{2}+|\tilde{k}_{[24]}|^{2})}dk_{1234}
\\\lesssim&\varepsilon^\kappa \int2^{(-2+2\kappa)q}\theta(2^{-q}\tilde{k}_{[23]})^2\frac{1}{|k_2|^{2-\kappa}|k_3|^{2}}dk_{23}\lesssim\varepsilon^\kappa2^{3q\kappa},\endaligned$$where $\kappa>0$ is small enough. Here we used $|k_i|\backsimeq N$ for some $i\in \{1,2,4\}$ in the second inequality and
Lemmas 6.1, 6.6 in the  third inequality.

\no\textbf{Term in the fourth chaos}: Now for $I_t^1$ we have the following calculations:
$$\aligned E[|\Delta_qI_t^1|^2]\lesssim&\varepsilon^\kappa\int\theta(2^{-q}\tilde{\tilde{k}}_{[1234]})^2\psi_0( {\tilde{k}}_{[12]},\tilde{k}_{[34]})
\frac{1}{|k_{[34]}||{\tilde{k}}_{[12]}|^{5-\kappa}}dk_{[12][34]}\\\lesssim&\varepsilon^\kappa2^{-2q}\int\theta(2^{-q}\tilde{k}_{[1234]})^2
\frac{1}{|k_{[34]}||{\tilde{k}}_{[12]}|^{3-\kappa}}dk_{[12][34]}
\lesssim\varepsilon^\kappa2^{q\kappa},\endaligned$$
where we used Lemmas 6.1, 6.5 in the last inequality.
By a similar calculation as above, Gaussian hypercontractivity and Lemma 2.1  we obtain that for $\delta>0$ small enough, $p>1$
$$\pi_{0,\diamond}(K_1^\varepsilon,e_N^{i_1i_2i_3}(u_1^\varepsilon)^{\diamond,2})\rightarrow 0\textrm{ in } L^p(\Omega;C_T\mathcal{C}^{-\delta}).$$

\no\textbf{Convergence of $(u_1^{\varepsilon})^{\diamond,2}e_N^{i_1i_2i_3}$}
By a similar calculation as that for $(u_1^{\varepsilon})^{\diamond,2}$ we know that
 $$\aligned &E|\Delta_q[(u_1^{\varepsilon})^{\diamond,2}e_N^{i_1i_2i_3}]|^2
\lesssim  \int\theta(2^{-q}\tilde{k}_{[12]})^2\frac{1}{|k_1|^2|k_2|^2}dk_1dk_2
 \lesssim\varepsilon^\kappa2^{(\kappa+2)q}.\endaligned$$
Here $\kappa>0$ is small enough and in the last inequality we used  Lemmas 6.1, 6.6. Then by Gaussian hypercontractivity and  Lemma 2.1 we obtain that for every $\delta>0, p>1,$ $(u_1^{\varepsilon})^{\diamond,2}e_N^{i_1i_2i_3}\rightarrow 0$ in $L^p(\Omega;C_T\mathcal{C}^{-1-\delta})$.

\vskip.10in

\no\textbf{Convergence of $\pi_{0}( u_2^\varepsilon,e_N^{i_1i_2i_3} u^\varepsilon_1)$}
Now we have the following identity: for $t\in[0,T]$,
$$ \pi_0({u}_2^\varepsilon, e_N^{i_1i_2i_3}{u}_1^\varepsilon)(t)=-[I_t^1+3I_t^2+J_t^1+3J_t^2].$$
Here
 $$\aligned I_t^1=&2^{-\frac{9}{2}}\int e_{\tilde{k}_{[1234]}}\psi_0( k_{[123]},\tilde{k}_4)\int_0^td\sigma P_{t-\sigma}^\varepsilon(k_{[123]})
P_{\sigma-s_1}^\varepsilon(k_{1})P_{\sigma-s_2}^\varepsilon(k_{2})P_{\sigma-s_3}^\varepsilon(k_{3})P_{t-s_4}^\varepsilon(k_{4})W(d\eta_{1234}),\\ I_t^2=&2^{-\frac{9}{2}}\int \int e_{\tilde{k}_{[23]}}\psi_0( k_{[123]},\tilde{k}_1)\int_0^td\sigma P_{t-\sigma}^\varepsilon(k_{[123]})
P_{\sigma-s_2}^\varepsilon(k_{2})P_{\sigma-s_3}^\varepsilon(k_{3})V_{t-\sigma}^\varepsilon(k_{1})dk_1W(d\eta_{23}),\endaligned$$
and for $i=1,2,$ $J_t^i$ is  defined similarly as $I_t^3$ with each $k_{[123]}, e_{\tilde{k}_{[1234]}},   e_{\tilde{k}_{[23]}}$ replaced by $\tilde{k}_{[123]}, e_{\tilde{\tilde{k}}_{[1234]}}, e_{\tilde{\tilde{k}}_{[23]}}$.

\no\textbf{Term in the second chaos:}
First we consider $I_t^2$ and by a similar calculations as that for $\pi_0({u}_2^\varepsilon, {u}_1^\varepsilon)$:
$$\aligned &E|\Delta_qI_t^2|^2
\lesssim \varepsilon^\kappa\int2^{-q(2-2\kappa)}\theta(2^{-q}\tilde{k}_{[23]})^2\frac{1}{|k_2|^{2}|k_3|^{2}}
dk_{23}\lesssim \varepsilon^\kappa 2^{q4\kappa},\endaligned$$
where $\kappa>0$ are small enough and we used $|k_{[123]}|\backsimeq |\tilde{k}_1|\backsimeq N$ in the first inequality, we used Lemmas 6.1, 6.5 in the last inequality.
By a similar calculation as above, we know that
 $$\aligned &E|\Delta_qJ_t^2|^2
\lesssim \varepsilon^\kappa\int 2^{-q(2-2\kappa)}\theta(2^{-q}\tilde{\tilde{k}}_{[23]})^2\frac{1}{|k_2|^{2}|k_3|^{2}}
dk_{23}
 \lesssim\varepsilon^\kappa2^{3\kappa q}.\endaligned$$
Here $\kappa>0$ is small enough  and we used $|\tilde{k}_{[123]}|\backsimeq |\tilde{k}_1|\backsimeq N$ in the first inequality, we used Lemmas 6.1, 6.5 in the last inequality.

\no\textbf{Terms in the fourth chaos}: Now for $I_t^1, J_t^1$ we have similar estimates:
 $$\aligned &E[|\Delta_qI_t^1|^2+|\Delta_qJ_t^1|^2]
 \lesssim\varepsilon^\kappa2^{\kappa q}.\endaligned$$
 Here for $I_t^1$ we used $|k_{[123]}|\simeq \tilde{k}_4\backsimeq N$ and for $J_t^1$ we used $|k_{[123]}|\simeq N\gtrsim\tilde{k}_{[123]}$.
Now by a similar calculation as above , Gaussian hypercontractivity and Lemma 2.1 we have that for $\delta>0$, $p>1$
$$\pi_{0}( u_2^\varepsilon,e_N^{i_1i_2i_3} u^\varepsilon_1)\rightarrow 0\textrm{ in } L^p(\Omega;C_T\mathcal{C}^{-\delta}).$$
\textbf{Convergence of $\pi_{0,\diamond}( u_2^\varepsilon,e_N^{i_1i_2i_3}(u_1^\varepsilon)^{\diamond,2})$}
Now we have the following identity: for $t\in[0,T]$,
$$ \pi_{0,\diamond}({u}_2^\varepsilon, e_N^{i_1i_2i_3}({u}_1^\varepsilon)^{\diamond2})=-[I_t^1+6I_t^2+6I_t^3+J_t^1+6J_t^2+6J_t^3]$$
 $$\aligned I_t^1=&2^{-6}\int e_{\tilde{k}_{[12345]}}\psi_0( k_{[123]},\tilde{k}_{[45]})\int_0^td\sigma P^\varepsilon_{t-\sigma}(k_{[123]})\Pi_{i=1}^3P^\varepsilon_{\sigma-s_i}(k_{i})
 \Pi_{i=4}^5P^\varepsilon_{t-s_i}(k_{i})W(d\eta_{12345}),\\ I_t^2=&2^{-6}\int\int e_{\tilde{k}_{[234]}}\psi_0( k_{[123]},\tilde{k}_{4}-k_1)\int_0^td\sigma P^\varepsilon_{t-\sigma}(k_{[123]})\Pi_{i=2}^3P^\varepsilon_{\sigma-s_i}(k_{i})
 P^\varepsilon_{t-s_4}(k_{4})V_{t-\sigma}^\varepsilon(k_1)dk_1W(d\eta_{234}),\\I_t^3=&2^{-6}\int\int e_{\tilde{k}_3}\psi_0 (k_{[123]},\tilde{k}_{[12]})\int_0^td\sigma P_{\sigma-s_3}^{\varepsilon}(k_3)V^\varepsilon_{t-\sigma}(k_1)V^\varepsilon_{t-\sigma}(k_2)P_{t-\sigma}^\varepsilon(k_{[123]})dk_{12}W(d\eta_3),\endaligned$$
  and for $i=1,2,3$, $J_t^i$ is  defined similarly as $I_t^i$ with each $k_{[123]}, e_{\tilde{k}_{[12345]}},   e_{\tilde{k}_{[234]}}, e_{\tilde{k}_{3}}$ replaced by $\tilde{k}_{[123]}, e_{\tilde{\tilde{k}}_{[12345]}}, e_{\tilde{\tilde{k}}_{[234]}}, e_{\tilde{\tilde{k}}_{3}}$.

\no \textbf{Terms in the first chaos} We consider $J_t^3$ and $I_t^3$ can be estimated similarly:
 $J_t^3=J_t^{31}+J_t^{32}$, with $J_t^{31}, J_t^{32}$ associated with the terms that $\tilde{\tilde{k}}_3\neq k_3$ and $\tilde{\tilde{k}}_3= k_3$, respectively.
 For $J_t^{31}$ we have
  $$\aligned &E[|\Delta_qJ_t^{31}|^2]\\\lesssim& \int\frac{1_{|k_3|\lesssim N}}{|k_3|^2}\theta(2^{-q}\tilde{\tilde{k}}_3)\bigg{(}\int
\frac{1_{|k_{1}|\lesssim N, |k_2|\lesssim N}}{|k_1|^2|k_2|^2(|k_1|^2+|k_2|^2+|\tilde{k}_{[123]}|^2)}dk_{12}\bigg{)}^2dk_3
\lesssim\varepsilon^{\kappa}2^{q(1+3\kappa)}.\endaligned$$
Here we used $|\tilde{\tilde{k}}_3|\backsimeq2^q\backsimeq N$ in the last inequality.
 For $J_t^{32}$ we have
 $$J_t^{32}-\tilde{J}_t^{32}+\tilde{J}_t^{32}-C_2^\varepsilon(t) u_1^{\varepsilon}(t),$$
 with $\tilde{J}_t^{32}$ defined as $J_t^{32}$ with $P_{\sigma-s_3}^{\varepsilon}(k_3) $ replaced by $P_{t-s_3}^{\varepsilon}(k_3)$ and with $C_2^\varepsilon(t)=\frac{1}{2}(C_{12}^\varepsilon+\varphi_{2}^\varepsilon(t))$.

Since $\int|P_{t-s_3}^\varepsilon(k_3)-P_{\sigma-s_3}^\varepsilon(k_3)|^2ds_3\leq C\frac{(t-\sigma)^{\kappa/2}}{|k_3|^{2-\kappa}}$, by a straightforward calculation we obtain that for $\kappa>0$ small enough
$$\aligned& E[|\Delta_q(J_t^{32}-\tilde{J}_t^{32})|^2]
\\\lesssim &\int\theta(2^{-q}k_3)^2\frac{1}{|k_3|^{2-4\kappa}}\bigg(\int_0^t\int
\frac{e^{-(|\tilde{k}_{[123]}|^2+|k_1|^2+|k_2|^2)\bar{c}_f(t-\sigma)}}{|k_1|^2|k_2|^2}
\\&(t-\sigma)^{\kappa}dk_{12}d\sigma\bigg)^2dk_3
\lesssim \varepsilon^\kappa 2^{q(1+5\kappa)}.\endaligned$$
Here in the last inequality we used that $|k_{123}|\backsimeq N$ implies that $|k_i|\backsimeq N$ for some $i\in\{1,2,3\}$ and  $\sup_{a\geq0} a^re^{-a}\leq C$ for $r\geq0$ and Lemma 6.1.
Moreover, by Lemmas 6.2 and 6.3 we obtain that
$$\aligned &E[|\Delta_q(\tilde{J}_t^{32}-u_1^{\varepsilon}(t) {C}_2^{\varepsilon}(t))|^2]\\\lesssim& \int\frac{1}{|k_3|^2}\theta(2^{-q}k_3)\bigg{(}\int\int_0^t|\tilde{k}_{[12]}|^{-\kappa}|k_3|^\kappa
\frac{e^{-|k_1|^2(t-\sigma)\bar{c}_f-|k_2|^2(t-\sigma)\bar{c}_f}}{|k_1|^2|k_2|^2}dk_{12}d\sigma\bigg{)}^2dk_3
\\\lesssim&\int\theta(2^{-q}k_3)\frac{1}{|k_3|^{2-2\kappa}}dk_3\bigg[\int_{|k_{[12]}|\leq N}\frac{1}{|\tilde{k}_{[12]}|^\kappa|k_{[12]}|^3}dk_{[12]}+\varepsilon^{\kappa/2}\int_{|k_{[12]}|>N}\frac{1}{|\tilde{k}_{[12]}|^{3+\kappa/2}}dk_{[12]}\bigg]^2
\lesssim\varepsilon^{\kappa}2^{q(1+3\kappa)},\endaligned$$
where in the last inequality we used that if $|k_{[12]}|\leq N$, $|\tilde{k}_{[12]}|\backsimeq N$.

Terms in the third and fifth chaos can be estimated similarly as that for $\pi_{0,\diamond}(u_2^\varepsilon, (u_1^\varepsilon)^{\diamond2})$
 and we also obtain that there exist $\kappa,\epsilon,\gamma>0$ small enough such that
$$\aligned &E[|\Delta_q(\pi_{0,\diamond}(u_2^\varepsilon, e_N^{i_1i_2i_3}(u_1^\varepsilon)^{\diamond2})(t_1)-\pi_{0,\diamond}(u_2^\varepsilon, e_N^{i_1i_2i_3}(u_1^\varepsilon)^{\diamond2})(t_2))|^2]\\\lesssim& \varepsilon^\gamma|t_1-t_2|^{\kappa}2^{q(1+\epsilon)},\endaligned$$
which by Gaussian hypercontractivity and Lemma 2.1 implies that for every $\delta>0, p>1$, $\pi_{0,\diamond}(u_2^\varepsilon, e_N^{i_1i_2i_3}(u_1^\varepsilon)^{\diamond2})\rightarrow0$ in $L^p(\Omega;C_T\mathcal{C}^{-1/2-\delta/2})$ .

\subsection{Convergence of the random operator}
The purpose of this subsection is to prove that $A_N$ defined in Proposition 4.3 converges to zero in probability. Here we follow essentially the same arguments as [GP15, Section 10.2].
\vskip.10in
\th{Theorem 6.7} For every $T\geq0$, $0<\eta<\delta,r\geq1$ we have
$$E[(A_N)^r]^{1/r}\lesssim N^{-\delta+\eta}.$$
\vskip.10in
By a similar argument as [GP15, Section 10.2] we have the following lemmas.
\vskip.10in
\th{Lemma 6.8} We have
$$\aligned &A_N^1(K^\varepsilon+K_1^\varepsilon,(u_1^{\varepsilon})^{\diamond,2}+e_N^{i_1i_2i_3}(u_1^{\varepsilon})^{\diamond,2})
(f)+A_N^2(\tilde{K}^\varepsilon+\tilde{K}_1^\varepsilon,(u_1^{\varepsilon})^{\diamond,2}+e_N^{i_1i_2i_3}(u_1^{\varepsilon})^{\diamond,2})(f)
\\=&\sum_{p,q}\int_{\mathbb{T}^3}g^N_{p,q}(t,x,y)\Delta_pf(y)dy\endaligned$$
with $$\mathcal{F}g^N_{p,q}(t,x,\cdot)(k)=\sum_{k_1,k_2}\Gamma^N_{p,q}(x,k,k_1,k_2)\mathcal{F}(\tilde{K}^\varepsilon+\tilde{K}_1^\varepsilon)(t,k_1)\mathcal{F}((u_1^{\varepsilon})^{\diamond,2}+e_N^{i_1i_2i_3}(u_1^{\varepsilon})^{\diamond,2})(t,k_2).$$
Here $$\aligned\Gamma^N_{p,q}(x,k,k_1,k_2)=&2^{-9/2}e^{\imath (k_1+k_2-k)\pi x}\theta_q(k_1+k_2-k)\tilde{\theta}_p(k)\psi_<(k,k_1)\psi_0(k_1-k,k_2)\\&(1_{|k_1-k|_\infty>N}1_{|k_1|_\infty\leq N}+1_{|k_1-k|_\infty\leq N}1_{N<|k_1|_\infty\leq3 N}),\endaligned$$
with $\tilde{\theta}_p$ be a smooth function supported in an annulus $2^p\mathcal{A}$ such that $\tilde{\theta}_p\theta_p=\theta_p$.
\vskip.10in
\th{Lemma 6.9} For all $r\geq1$ we have
$$\aligned &E[\|(A_N^1(K^\varepsilon+K_1^\varepsilon,(u_1^{\varepsilon})^{\diamond,2}+e_N^{i_1i_2i_3}(u_1^{\varepsilon})^{\diamond,2})+A_N^2(\tilde{K}^\varepsilon+\tilde{K}_1^\varepsilon,(u_1^{\varepsilon})^{\diamond,2}+e_N^{i_1i_2i_3}(u_1^{\varepsilon})^{\diamond,2}))(t)
\\&-(A_N^1(K^\varepsilon+K_1^\varepsilon,(u_1^{\varepsilon})^{\diamond,2}+e_N^{i_1i_2i_3}(u_1^{\varepsilon})^{\diamond,2})+A_N^2(\tilde{K}^\varepsilon+\tilde{K}_1^\varepsilon,(u_1^{\varepsilon})^{\diamond,2}+e_N^{i_1i_2i_3}(u_1^{\varepsilon})^{\diamond,2}))(s)\|_{L(\mathcal{C}^{1-\delta},B^{-1/2-2\delta+\kappa}_{r,r})}^r]
\\\lesssim&\sum_{p,q}2^{qr(-1/2-2\delta+\kappa)}2^{-pr(1-\delta)}(\sup_{x\in\mathbb{T}^3}\sum_kE[|(\mathcal{F}g^N_{p,q}(t,x,\cdot)-\mathcal{F}g^N_{p,q}(s,x,\cdot))(k)|^2])^{r/2}.\endaligned $$
\vskip.10in
\th{Lemma 6.10} For all $p,q\geq-1$, all $0\leq t_1<t_2,$ and all $\lambda,\kappa\in(0,1]$ we have
$$\sum_kE[|(\mathcal{F}g^N_{p,q}(t_2,x,\cdot)-\mathcal{F}g^N_{p,q}(t_1,x,\cdot))(k)|^2]\lesssim1_{2^p,2^q\lesssim N}(2^{3p}2^{2q}+2^{2p}2^{3q})N^{-2+2\lambda+\kappa}|t_1-t_2|^\lambda.$$

\proof We only prove the estimate for $\sum_kE[|\mathcal{F}g^N_{p,q}(t,x,\cdot)(k)|^2])$ and the result can be obtained by essentially the same arguments. First we consider $\mathcal{F}\tilde{K}^\varepsilon(t,l_1)\mathcal{F}(u_1^{\varepsilon})^{\diamond,2}(t,l_2)$. We have the following chaos decomposition:
$$\aligned&\mathcal{F}\tilde{K}^\varepsilon(t,l_1)\mathcal{F}(u_1^{\varepsilon})^{\diamond,2}(t,l_2)=I_t^1+I_t^2+I_t^3.\endaligned$$
Here $$\aligned I_t^1=&2^{-3}\int 1_{k_{[12]}=l_1,k_{[34]}=l_2}\int_0^td\sigma p^\varepsilon_{t-\sigma}(k_{[12]})\varphi(\varepsilon k_{[12]})P^\varepsilon_{\sigma-s_1}(k_1)P^\varepsilon_{\sigma-s_2}(k_2)P^\varepsilon_{t-s_3}(k_3)P_{t-s_4}^\varepsilon(k_4)W(d\eta_{1234}),
\\ I_t^2=&2^{-3}\int1_{k_{[12]}=l_1,k_{3}-k_1=l_2}\int_0^td\sigma p^\varepsilon_{t-\sigma}(k_{[12]})\varphi(\varepsilon k_{[12]})P^\varepsilon_{\sigma-s_2}(k_2)P^\varepsilon_{t-s_3}(k_3)V_{t-\sigma}^\varepsilon(k_1)dk_1W(d\eta_{23}),
 \\I_t^3=&2^{-3}\int
\int_0^td\sigma 1_{k_{[12]}=l_1,-k_{[12]}=l_2}V^\varepsilon_{t-\sigma}(k_1)V^\varepsilon_{t-\sigma}(k_2)p^\varepsilon_{t-\sigma}(k_{[12]})\varphi(\varepsilon k_{[12]})dk_{12},\endaligned$$
\no\textbf{Term in the chaos of order 0}
By a similar calculation as in Section 6.1 we have
$$\aligned&\sum_k\big{|}\sum_{k_1,k_2}\Gamma^N_{p,q}(x,k,k_1,k_2)1_{k_1+k_2=0}I_t^3\big{|}^2
\\\lesssim&\sum_k\big{|}\sum_{k_1}\Gamma^N_{p,q}(x,k,k_1,-k_1)\frac{1}{|k_1|^3}\big{|}^2
\\\lesssim &\sum_k\tilde{\theta}_p(k)^2\theta_q(-k)^2\big{|}\sum_{k_1}(1_{|k_1-k|_\infty>N}1_{|k_1|_\infty\leq N}+1_{|k_1-k|_\infty\leq N}1_{N<|k_1|_\infty\leq3 N})\psi_<(k,k_1)\psi_0(k_1-k,k_1)\frac{1}{|k_1|^3}\big{|}^2.\endaligned$$
In the first case without loss of generality we assume that $|k_1^i-k^i|>N$ for some $i$. Then there are at most $|k^i|$ values of $k_1^i$ with $|k_1^i|<N$ and $|k_1^i-k^i|>N$. In the second case without loss of generality we assume that $|k_1^i|>N$ for some $i$. Then there are at most $|k^i|$ values of $k_1^i$ with $|k_1^i|>N$ and $|k_1^i-k^i|\leq N$. Moreover observe that $|k_1|\simeq N$ on the support of $(1_{|k_1-k|_\infty>N}1_{|k_1|_\infty\leq N}+1_{|k_1-k|_\infty\leq N}1_{|k_1|_\infty>N})\psi_0(k-k_1,k_1)$ and that $|k|\lesssim N$ whenever $1_{|k_1|_\infty\leq3N}\psi_<(k,k_1)\neq0$, which implies that the above term is bounded by
$$\sum_k\tilde{\theta}_p(k)^2\theta_q(-k)^2|k|^21_{|k|\lesssim N}N^{-2}\lesssim 1_{2^p,2^q\lesssim N}2^{3p}2^{2q}N^{-2}.$$

\no\textbf{Term in the second chaos} By a similar calculation as in Section 6.1 we have
$$\aligned&\sum_kE\big{|}\sum_{l_1,l_2}\Gamma^N_{p,q}(x,k,l_1,l_2)I_t^2\big{|}^2
\\\lesssim&\sum_k1_{2^p,2^q\lesssim N}\tilde{\theta}_p(k)^2\int\theta_q(k_{[23]}-k)^2\Pi_{i=2}^3\frac{1}{|k_i|^2}\bigg[\int\psi_<(k,k_{[12]})\frac{1}{(|k_{[12]}|^2+|k_1|^2)|k_1|^2}
\\&(1_{|k_{[12]}-k|_\infty>N,|k_{[12]}|_\infty\leq N}+1_{|k_{[12]}-k|_\infty\leq N,N<|k_{[12]}|_\infty\leq3 N})dk_1\bigg]^2dk_3
\\\lesssim &\sum_k1_{2^p,2^q\lesssim N}\tilde{\theta}_p(k)^2\int\theta_q(k_{[23]}-k)^2\frac{1}{|k_{[23]}|}N^{-2+\kappa}dk_{[23]}
\\\lesssim &\sum_k1_{2^p,2^q\lesssim N}\tilde{\theta}_p(k)^2N^{-2+\kappa}2^{2q}\lesssim1_{2^p,2^q\lesssim N}2^{3p}2^{2q}N^{-2+\kappa}.\endaligned$$
Here in the second inequality we used that $|k_{[12]}|\simeq N$ on the support of $1_{|k_{[12]}-k|_\infty>N,|k_{[12]}|_\infty\leq N}\psi_<(k,k_{[12]})$
and in the third inequality we used Lemma 6.5.

\no\textbf{Term in the forth chaos} We have
$$\aligned&\sum_kE\big{|}\sum_{l_1,l_2}\Gamma^N_{p,q}(x,k,l_1,l_2)I_t^1\big{|}^2
\\\lesssim&\sum_k\tilde{\theta}_p(k)^2\int\theta_q(k_{[1234]}-k)^2\psi_<(k,k_{[12]})^2\psi_0(k_{[12]}-k,k_{[34]})^2\frac{1_{2^p,2^q\lesssim N}}{|k_2|^2|k_3|^2|k_1|^2|k_4|^2|k_{[12]}|^4}
\\&(1_{|k_{[12]}-k|_\infty>N,|k_{[12]}|_\infty\leq N}+1_{|k_{[12]}-k|_\infty\leq N,N<|k_{[12]}|_\infty\leq3 N})dk_{1234}
\\\lesssim &\sum_k\tilde{\theta}_p(k)^2\int\theta_q(k_{[1234]}-k)^2\frac{1_{2^p,2^q\lesssim N}}{|k_{[1234]}|^{1-\kappa}}dk_{[1234]}N^{-2-\kappa}
\\\lesssim &1_{2^p,2^q\lesssim N}\sum_k\tilde{\theta}_p(k)^2N^{-2-\kappa}2^{(2+\kappa)q}\\\lesssim&1_{2^p,2^q\lesssim N}2^{3p}2^{2q}N^{-2},\endaligned$$
where we used Lemma 6.1 and $|k_{[12]}|\simeq N$ in the second inequality.

Moreover we consider $$\mathcal{F}\tilde{K}_1^\varepsilon(t,k_1)\mathcal{F}(u_1^{\varepsilon})^{\diamond,2}(t,k_2)=J_t^1+J_t^2+J_t^3.$$
Here $J_t^i, i=1,2,$ is defined similar as $I_t^i, i=1,2,$ with $k_{[12]},k_{[1234]},k_{[23]},e_{k_{[23]}},e_{k_{1234}}$ replaced by $\tilde{k}_{[12]},\tilde{k}_{[1234]},\tilde{k}_{[23]},e_{\tilde{k}_{1234}}$, respectively and $$J_t^3=2^{-3}\int
\int_0^td\sigma 1_{\tilde{k}_{[12]}=l_1,-k_{[12]}=l_2}V^\varepsilon_{t-\sigma}(k_1)V^\varepsilon_{t-\sigma}(k_2)p^\varepsilon_{t-\sigma}(\tilde{k}_{[12]})\varphi(\varepsilon \tilde{k}_{[12]})dk_{12}.$$

\no\textbf{Terms in the chaos of order 0}
We have
$$\aligned&\sum_k\big{|}\sum_{k_1,k_2}\Gamma^N_{p,q}(x,k,\tilde{k}_1,k_2)1_{k_1+k_2=0}J_t^3\big{|}^2
\\\lesssim&\sum_k\big{|}\sum_{k_1}\Gamma^N_{p,q}(x,k,\tilde{k}_1,-k_1)1_{|k_1|_\infty\lesssim N}\frac{1}{|k_1||\tilde{k}_1|^2}\big{|}^2
\\\lesssim &\sum_k\tilde{\theta}_p(k)^2\theta_q(\tilde{k})^2\big{|}\sum_{k_1}(1_{|\tilde{k}_1-k|_\infty>N}1_{|\tilde{k}_1|_\infty\leq N}+1_{|\tilde{k}_1-k|_\infty\leq N}1_{N<|\tilde{k}_1|_\infty\leq3 N})\psi_<(k,\tilde{k}_1)1_{|k_1|_\infty\lesssim N}\frac{1}{|k_1||\tilde{k}_1|^2}\big{|}^2.\endaligned$$
Similarly as above we obtain there are at most $|k^i|$ values of $\tilde{k}_1^i$ with $|\tilde{k}_1^i|>N$ and $|\tilde{k}_1^i-k^i|\leq N$ or $|k_1^i|>N$ and $|k_1^i-k^i|\leq N$. Moreover observe that $|\tilde{k}_1|\simeq N$ on the support of $1_{|\tilde{k}_1-k|_\infty>N}1_{|\tilde{k}_1|_\infty\leq N}\psi_<(k,\tilde{k}_1)$ and that $|k|\lesssim N$ whenever $1_{|k_1|_\infty<3N}\psi_<(k,k_1)\neq0$, which implies that the above term is bounded by
$$\sum_k\tilde{\theta}_p(k)^2\theta_q(\tilde{k})^2|k|^21_{|k|\lesssim N}N^{-2}\lesssim 1_{2^p,2^q\lesssim N}2^{2p}2^{3q}N^{-2}.$$
For the terms in the second chaos by a similar calculation as above we obtain the estimates.

\no\textbf{Terms in the forth chaos} We have
$$\aligned&\sum_kE\big{|}\sum_{l_1,l_2}\Gamma^N_{p,q}(x,k,l_1,l_2)J_t^3\big{|}^2
\\\lesssim&\sum_k\tilde{\theta}_p(k)^2\int\theta_q(\tilde{k}_{[1234]}-k)^2\psi_<(k,\tilde{k}_{[12]})^2\psi_0(\tilde{k}_{[12]},k_{[34]})^2\frac{1_{2^p,2^q\lesssim N}}{|k_2|^2|k_3|^2|k_1|^2|k_4|^2|\tilde{k}_{[12]}|^4}1_{|k_{[12]}|\lesssim N,|k_{[34]}|\lesssim N}dk_{1234}
\\&(1_{|\tilde{k}_{12}-k|_\infty>N,|\tilde{k}_{12}|_\infty\leq N}+1_{|\tilde{k}_{12}-k|_\infty\leq N,N<|\tilde{k}_{12}|_\infty\leq3 N})
\\\lesssim &1_{2^p,2^q\lesssim N}\sum_k\tilde{\theta}_p(k)^2\int\theta_q(\tilde{k}_{[1234]}-k)^2\frac{1}{|k_1|^2|k_2|^2|k_3|^2|k_4|^2}N^{-4}1_{|k_{[12]}|\lesssim N,|k_{[34]}|\lesssim N}dk_{1234}
\\\lesssim &1_{2^p,2^q\lesssim N}\sum_k\tilde{\theta}_p(k)^2\int\theta_q(\tilde{k}_{[1234]}-k)\frac{1}{|k_{[12]}|^2|k_{[34]}|^2}dk_{[12][34]}N^{-2}\\\lesssim&1_{2^p,2^q\lesssim N}2^{3p}2^{2q}N^{-2}.\endaligned$$
Here in the second inequality we used $|\tilde{k}_{[12]}|\backsimeq N$ and in the third inequality we used $N^{-1}\lesssim |k_{[12]}|^{-1},N^{-1}\lesssim |k_{[34]}|^{-1}$
and in the last inequality we used Lemma 6.5.

Furthermore for the terms associated with $\tilde{K}^\varepsilon(u_1^{\varepsilon})^{\diamond,2}e^{i_1i_2i_3}_N$
and $\tilde{K}_1^\varepsilon(u_1^{\varepsilon})^{\diamond,2}e^{i_1i_2i_3}_N$ we can also obtain similar estimates by similar arguments.
Thus the result follows.$\hfill\Box$

\vskip.10in
\no\emph{Proof of Theorem 6.10} For $t,s\geq0$ we have
$$\aligned &E[\|(A_N(K^\varepsilon+K_1^\varepsilon,(u_1^{\varepsilon})^{\diamond,2}+e_N^{i_1i_2i_3}(u_1^{\varepsilon})^{\diamond,2})+B_N(\tilde{K}^\varepsilon+\tilde{K}_1^\varepsilon,(u_1^{\varepsilon})^{\diamond,2}+e_N^{i_1i_2i_3}(u_1^{\varepsilon})^{\diamond,2}))(t)
\\&-(A_N(K^\varepsilon+K_1^\varepsilon,(u_1^{\varepsilon})^{\diamond,2}+e_N^{i_1i_2i_3}(u_1^{\varepsilon})^{\diamond,2})+B_N(\tilde{K}^\varepsilon+\tilde{K}_1^\varepsilon,(u_1^{\varepsilon})^{\diamond,2}+e_N^{i_1i_2i_3}(u_1^{\varepsilon})^{\diamond,2}))(s)\|_{L(\mathcal{C}^{1-\delta},\mathcal{C}^{-1/2-2\delta})}^r]^{1/r}
\\\lesssim & E[\|(A_N(K^\varepsilon+K_1^\varepsilon,(u_1^{\varepsilon})^{\diamond,2}+e_N^{i_1i_2i_3}(u_1^{\varepsilon})^{\diamond,2})+B_N(\tilde{K}^\varepsilon+\tilde{K}_1^\varepsilon,(u_1^{\varepsilon})^{\diamond,2}+e_N^{i_1i_2i_3}(u_1^{\varepsilon})^{\diamond,2}))(t)
\\&-(A_N(K^\varepsilon+K_1^\varepsilon,(u_1^{\varepsilon})^{\diamond,2}+e_N^{i_1i_2i_3}(u_1^{\varepsilon})^{\diamond,2})+B_N(\tilde{K}^\varepsilon+\tilde{K}_1^\varepsilon,(u_1^{\varepsilon})^{\diamond,2}+e_N^{i_1i_2i_3}(u_1^{\varepsilon})^{\diamond,2}))(s)\|_{L(\mathcal{C}^{1-\delta},B^{-1/2-2\delta+\kappa}_{r,r})}^r]^{1/r}
\\\lesssim&\bigg[\sum_{p,q}2^{qr(-1/2-2\delta+\kappa)}2^{-pr(1-\delta)}1_{2^p,2^q\lesssim N}[(2^{3p}2^{2q}+2^{2p}2^{3q})|t-s|^\lambda N^{-2+2\lambda+\kappa}]^{r/2}\bigg]^{1/r}
\\\lesssim&|t-s|^{\lambda/2}N^{-\delta+2\kappa+\lambda}.\endaligned $$
Here $\delta>2\kappa+\lambda>0$. Thus the result follows by using Kolmogorov's continuity criterion.$\hfill\Box$

\subsection{Convergence of $D^N$}
In this subsection we prove that $D^N\rightarrow^P0$ as $\varepsilon\rightarrow0$.
Now we have the following identity: for $t\in[0,T]$,
$$ \pi_0((I-P_N)\pi_<({u}_2^\varepsilon, K^\varepsilon),(u_1^{\varepsilon})^{\diamond,2})(t)+\pi_0(P_N\pi_<({u}_2^\varepsilon, (P_{3N}-P_N)\tilde{K}^\varepsilon),(u_1^{\varepsilon})^{\diamond,2})(t)=\sum_{i=1}^4(I_t^i+J_t^i).$$
Here
 $$\aligned I_t^1=&2^{-9}\int e_{k_{[1234567]}}\psi_0( k_{[12345]},k_{[67]})\psi_<(k_{[123]},k_{[45]})(1_{|k_{[12345]}|_\infty>N}1_{|k_{[45]}|_\infty\leq N}+1_{|k_{[12345]}|_\infty\leq N}1_{N<|k_{[45]}|_\infty\leq 3N})\\&\int_0^t\int_0^td\sigma d\bar{\sigma} P^\varepsilon_{t-\sigma}(k_{[123]})\Pi_{i=1}^3P^\varepsilon_{\sigma-s_i}(k_{i})
 p^\varepsilon_{t-\bar{\sigma}}(k_{[45]})\varphi(\varepsilon k_{[45]})\Pi_{i=4}^5P^\varepsilon_{\bar{\sigma}-s_i}(k_{i})\Pi_{i=6}^7P^\varepsilon_{t-s_i}(k_{i})W(d\eta_{1234567})\\:=&\int G(t,x,\eta_{1234567})W(d\eta_{1234567}),\\ I_t^{2}=&\sum_{i=1}^3I_t^{2i}, I_t^{21}=6\int \int G(t,x,\eta_{123(-3)567})d\eta_3W(d\eta_{12567}), \\I_t^{22}=&6\int G(t,x,\eta_{12345(-3)7})d\eta_{3}W(d\eta_{12457}), I_t^{23}=4\int G(t,x,\eta_{12345(-5)7})]d\eta_5W(d\eta_{12347})  ,\\ I_t^{3}=&\sum_{i=1}^6I_t^{3i}, I_t^{31}=6\int \int G(t,x,\eta_{123(-3)(-2)67})d\eta_{23}W(d\eta_{167}),\\I_t^{32}=&24\int\int G(t,x,\eta_{123(-3)5(-2)7})d\eta_{23}W(d\eta_{157}),I_t^{33}=12\int \int G(t,x,\eta_{123(-3)5(-5)7})d\eta_{35}W(d\eta_{127})\\I_t^{34}=&6\int\int G(t,x,\eta_{12345(-2)(-3)})d\eta_{23}W(d\eta_{145}),I_t^{35}=12\int \int G(t,x,\eta_{12345(-3)(-4)})d\eta_{34}W(d\eta_{125})\\I_t^{36}=&2\int \int G(t,x,\eta_{12345(-4)(-5)})d\eta_{45}W(d\eta_{123}),\\I_t^{4}=&\sum_{i=1}^3I_t^{4i}, I_t^{41}=12\int\int G(t,x,\eta_{123(-1)(-2)(-3)7})d\eta_{123}W(d\eta_{7}),\\I_t^{42}=&12\int\int G(t,x,\eta_{123(-3)5(-2)(-1)})d\eta_{123}W(d\eta_{5}),I_t^{43}=24\int\int G(t,x,\eta_{123(-3)5(-5){-2}})d\eta_{235}W(d\eta_{1}),\endaligned$$
and  $J_t^1$ is defined similarly as $I_t^1$ with  $k_{[123]}, k_{[12345]}, e_{k_{[1234567]}}$ replaced by $\tilde{k}_{[123]}, \tilde{k}_{[12345]}, e_{\tilde{k}_{[1234567]}}$ respectively and $J_t^i,i=2,3,4$ is defined similarly as $I_t^i$ with the $G$ replaced by that associated with $J^1$.

\no\textbf{Terms in the seventh chaos} Now we have
$$\aligned E|\Delta_q I_t^1|^2\lesssim&\int \theta(2^{-q}k_{[1234567]})\psi_0( k_{[12345]},k_{[67]})\psi_<(k_{[123]},k_{[45]})(1_{|k_{[12345]}|_\infty>N,|k_{[45]}|_\infty\leq N}+1_{|k_{[12345]}|_\infty\leq N,|k_{[45]}|_\infty> N})\\&1_{|k_{[1234567]}|\lesssim N}\Pi_{i=1}^7\frac{1}{|k_i|^2}\frac{1}{|k_{[123]}|^2|k_{[45]}|^2(|k_{[123]}|^2+\sum_{i=1}^3|k_i|^2)(|k_{[45]}|^2+\sum_{i=4}^5|k_i|^2)}dk_{1234567},\endaligned$$
 Observe that $|k_{45}|_\infty\backsimeq N$ on the support of $\psi_<(k_{[123]},k_{[45]})1_{|k_{[12345]}|_\infty> N}$, which combining with  Lemma 6.1 implies that the above term can be bounded by $$\aligned&\int \theta(2^{-q}k_{[1234567]})1_{|k_{[45]}|_\infty\backsimeq N,2^q\lesssim N}\frac{1}{|k_{[123]}|^4|k_{[45]}|^5|k_{[67]}|}dk_{[123][45][67]} \\\lesssim&\int 1_{2^q\lesssim N} \theta(2^{-q}k_{[1234567]})\frac{N^{-2-\kappa}}{|k_{[12345]}|^{3-\kappa}}\frac{1}{|k_{[67]}|}dk_{[12345][67]}\lesssim\varepsilon^\kappa 2^{2q\kappa}.\endaligned$$

\no\textbf{Terms in the fifth chaos}:
Consider $I_t^{21}$ first: by the formula we know that $|k_{5}-k_3|\backsimeq N$
$$\aligned &E|\Delta_q I_t^{21}|^2\\\lesssim&1_{2^q\lesssim N}\int \theta(2^{-q}k_{[12567]})\Pi_{i=5}^7\frac{1}{|k_i|^2}\frac{1}{|k_1|^2|k_2|^2}
\bigg[\bigg(\int\frac{1}{|k_3|^2(|k_{5}-k_3|^2+|k_3|^2)(|k_{[123]}|^2+|k_{5}-k_3|^2)}dk_{3}\bigg)^2\\&+\bigg(
\int\frac{1}{|k_3|^2(|k_{5}-k_3|^2+|k_{[123]}|^2)(|k_{[123]}|^2+|k_3|^2)}dk_{3}\bigg)^2\bigg]1_{\{|k_5-k_3|\backsimeq N,|k_5|\leq N\}}dk_{12567}\\
\lesssim& \int\theta(2^{-q}k_{[12567]}) 1_{2^q\lesssim N}\frac{N^{-4+2\kappa}}{|k_{[12]}|^{3-\kappa}|k_5|^{2+2\kappa}|k_{[67]}|}dk_{[12]5[67]}\lesssim\varepsilon^\kappa 2^{2q\kappa}.\endaligned$$
Here in the first inequality we consider $\sigma\leq\bar{\sigma}$ and $\sigma\geq \bar{\sigma}$ separately and we used $|k_{[123]}|^2+|k_3|^2\gtrsim |k_{[12]}|^2$ in the second inequality.

Now we consider $I_t^{22}$ and in this case we have $|k_{[45]}|\backsimeq N$, which implies that
$$\aligned E|\Delta_q I_t^{22}|^2\lesssim&1_{2^q\lesssim N}\int \theta(2^{-q}k_{[12457]})1_{|k_{45}|\backsimeq N}\frac{1}{|k_1|^2|k_2|^2|k_{[45]}|^5|k_7|^2}(\int\frac{1}{(|k_{[123]}|^2+|k_3|^2)|k_3|^2}dk_3)^2dk_{12[45]7}
\\\lesssim&1_{2^q\lesssim N}\int \theta(2^{-q}k_{[12457]})\frac{N^{-2-\kappa}}{|k_{[12]}|^{3-\kappa}|k_{[45]}|^{3-\kappa}|k_7|^2}dk_{[12][45]7}\lesssim\varepsilon^\kappa 2^{2q\kappa}.\endaligned$$
Here in the second inequality we used $|k_{[123]}|^2+|k_3|^2\gtrsim |k_{[12]}|^2$.

For $I_t^{23}$ we have $|k_{[45]}|\backsimeq N$
$$\aligned E|\Delta_q I_t^{23}|^2\lesssim&1_{2^q\lesssim N}\int \theta(2^{-q}k_{[12347]})1_{\{|k_{[45]}|_\infty\backsimeq N\}}\Pi_{i=1}^4\frac{1}{|k_i|^2}\frac{1}{|k_7|^2}\\&\frac{1}{(|k_{[123]}|^2+\sum_{i=1}^3|k_i|^2)|k_{[123]}|^2}(\int\frac{1}{(|k_{[45]}|^2+|k_4|^2)|k_5|^2}dk_5)^2dk_{12347}
\\\lesssim&\int1_{2^q\lesssim N} \theta(2^{-q}k_{[12347]})\frac{1}{|k_{[1234]}|^{2-\kappa}}\frac{N^{-2+\kappa}}{|k_7|^2}dk_{[1234]7}\lesssim\varepsilon^\kappa 2^{2q\kappa},\endaligned$$
where we used $|k_{[45]}|^2+|k_4|^2\gtrsim |k_{5}|^2$ in the second inequality.

\no\textbf{Terms in the third chaos}:
For $I_t^{31}$, we have
$$\aligned E|\Delta_q I_t^{31}|^2\lesssim&1_{2^q\lesssim N}\int \theta(2^{-q}k_{[167]})\psi_0( k_{1},k_{[67]})\psi_<(k_{[123]},k_{[23]})(1_{|k_{1}|_\infty>N,|k_{[23]}|_\infty\leq N}+1_{|k_{1}|_\infty\leq N,N<|k_{[23]}|_\infty\leq 3N})\\&\frac{1}{|k_1|^2}\frac{1}{|k_6|^2|k_7|^2}(\int\frac{1}{(|k_{[123]}|^2+|k_{[23]}|^2)|k_{[23]}|^2|k_2|^2|k_3|^2}dk_{[23]})^2dk_{167}
\\\lesssim&1_{2^q\lesssim N}\int \theta(2^{-q}k_{[167]})\frac{N^{-3-\kappa}}{|k_1|^{3-\kappa}|k_{[67]}|}dk_{1[67]}\lesssim\varepsilon^\kappa 2^{2q\kappa}.\endaligned$$
Here we used $|k_{[23]}|^2\lesssim |k_2|^2+|k_3|^2$ in the first inequality and in the second inequality we used $|k_{[23]}|\backsimeq N$.

For $I_t^{32}$ we have$$\aligned E|\Delta_q I_t^{32}|^2\lesssim&1_{2^q\lesssim N}\int \theta(2^{-q}k_{[157]})1_{\{|k_5-k_3|_\infty\backsimeq N\}}\frac{1}{|k_1|^2|k_5|^2|k_7|^2}1_{|k_5|\leq N}\\&\bigg[\bigg(\int\frac{1}{(|k_5-k_3|^2+|k_3|^2)(|k_{[123]}|^2+|k_2|^2+|k_5-k_3|^2)|k_2|^2|k_3|^2}dk_{23}\bigg)^2
\\&+\bigg(\int\frac{1}{(|k_{[123]}|^2+|k_2|^2+|k_3|^2)(|k_{[123]}|^2+|k_2|^2+|k_5-k_3|^2)|k_2|^2|k_3|^2}dk_{23}\bigg)^2\bigg]
dk_{157}\\\lesssim&1_{2^q\lesssim N}\int \theta(2^{-q}k_{[157]})N^{-3}\frac{1}{|k_1|^2|k_5|^{3-\kappa}|k_7|^2}dk_{157}\lesssim\varepsilon^\kappa 2^{2q\kappa},\endaligned$$
Here  in the first inequality we consider $\sigma\leq\bar{\sigma}$ and $\sigma\geq \bar{\sigma}$ separately.
For $I_t^{33}$ we have
$$\aligned E|\Delta_q I_t^{33}|^2\lesssim &1_{2^q\lesssim N}\int \theta(2^{-q}k_{[127]})1_{|k_5-k_3|\backsimeq N,|k_{[12]}|\lesssim N}\frac{1}{|k_1|^2|k_2|^2|k_7|^2}\\&\bigg[\bigg(\int\frac{1}{(|k_3|^2+|k_5-k_3|^2+|k_5|^2)(|k_{[123]}|^2+|k_5-k_3|^2+|k_5|^2)|k_3|^2|k_5|^2}dk_{35}\bigg)^2
\\&+\bigg(\int\frac{1}{(|k_3|^2+|k_{[123]}|^2)(|k_{[123]}|^2+|k_5-k_3|^2+|k_5|^2)|k_3|^2|k_5|^2}dk_{35}\bigg)^2\bigg]dk_{127}
\\\lesssim&1_{2^q\lesssim N}\int \theta(2^{-q}k_{[127]})\frac{N^{-2+2\kappa}}{|k_{[12]}|^{3-\kappa}|k_7|^2}dk_{[12]7}\lesssim\varepsilon^\kappa 2^{2q\kappa}.\endaligned$$
For $I_t^{34}$ we have
$$\aligned E|\Delta_q I_t^{34}|^2\lesssim&1_{2^q\lesssim N}\int \theta(2^{-q}k_{[145]})1_{|k_{[45]}|\backsimeq N}\frac{1}{|k_1|^2|k_{[45]}|^5}\\&\bigg(\int\frac{1_{|k_{[23]}|\lesssim N}}{(|k_{[123]}|^2+\sum_{i=2}^3|k_i|^2)|k_2|^2|k_3|^2}dk_{23}\bigg)^2dk_{145}\\\lesssim&1_{2^q\lesssim N}\int \theta(2^{-q}k_{[145]})1_{|k_{[45]}|_\infty\backsimeq N}\frac{N^\kappa}{|k_{[45]}|^5|k_1|^2}dk_{1[45]}\lesssim\varepsilon^\kappa 2^{2q\kappa}.\endaligned$$
For $I_t^{35}$ we have
$$\aligned E|\Delta_q I_t^{35}|^2\lesssim&1_{2^q\lesssim N}\int \theta(2^{-q}k_{[125]})1_{|k_{[45]}|_\infty\backsimeq N}\frac{1}{|k_1|^2|k_2|^2|k_5|^2}\\&(\int\frac{1}{(|k_{[45]}|^2+|k_4|^2)(|k_{[123]}|^2+|k_3|^2)|k_3|^2|k_4|^2}dk_{34})^2dk_{125}\\\lesssim&1_{2^q\lesssim N}\int \theta(2^{-q}k_{[125]})\frac{N^{-2+\kappa}}{|k_{5}|^2|k_{[12]}|^{3-\kappa}}dk_{[12]5}\lesssim\varepsilon^\kappa 2^{2q\kappa}.\endaligned$$
Here  in the second inequality we used $|k_{[123]}|^2+|k_3|^2\gtrsim |k_{[12]}|^2$.

For $I_t^{36}$ we have $$\aligned E|\Delta_q I_t^{36}|^2\lesssim&1_{2^q\lesssim N}\int  \theta(2^{-q}k_{[123]})\frac{1}{|k_{[123]}|^4}\Pi_{i=1}^3\frac{1}{|k_i|^2}\\&\bigg(\int\frac{(1_{|k_{[12345]}|_\infty>N}1_{|k_{[45]}|_\infty\leq N}+1_{|k_{[12345]}|_\infty\leq N}1_{N<|k_{[45]}|_\infty\leq 3N})}{(|k_{[45]}|^2+\sum_{i=4}^5|k_i|^2)|k_4|^2|k_5|^2}dk_{45}\bigg)^2dk_{123},\endaligned$$
Now we use similar argument as Section 6.3. In the first case without loss of generality we assume that $|k_{[123]}^i+k_{[45]}^i|>N$ for some $i$. Then there are at most $|k_{[123]}^i|$ values of $k_{[45]}^i$ with $|k_{[12345]}^i|>N$ and $|k_{[45]}^i|\leq N$. In the second case similarly we obtain that there are at most $|k_{[123]}^i|$ values of $k_{[45]}^i$ with $|k_{[45]}^i|>N$ and $|k_{[12345]}^i|\leq N$. Thus we obtain
$$\aligned E|\Delta_q I_t^{36}|^2\lesssim&1_{2^q\lesssim N}\int \theta(2^{-q}k_{[123]})N^{-2+\kappa}\frac{1}{|k_{[123]}|^2}dk_{[123]}\lesssim\varepsilon^\kappa 2^{q\kappa}.\endaligned$$

\no\textbf{Terms in the first chaos}:
For $I^{41}$ we obtain that
$$\aligned E|\Delta_q I_t^{41}|^2\lesssim&\int \theta(2^{-q}k_7)\frac{1}{|k_7|^2}\bigg[\int\frac{1_{|k_{[12]}|\backsimeq N}}{|k_2|^2|k_3|^2|k_1|^2(|k_{[123]}|^2+|k_3|^2)|k_{[12]}|^2}dk_{123}\bigg]^2dk_{7}\\\lesssim&\int \theta(2^{-q}k_{7})\frac{N^{-2+\kappa}}{|k_7|^2}dk_{7}\lesssim\varepsilon^\kappa 2^{2q\kappa}.\endaligned$$
For $I^{42}$  we obtain that
$$\aligned E|\Delta_q I_t^{42}|^2\lesssim&\int \theta(2^{-q}k_{5})\bigg(\int1_{|k_5-k_3|\backsimeq N }\Pi_{i=1}^3\frac{1}{|k_i|^2}\frac{1_{\{|k_i|_\infty\leq N, i=1,2,3\} }}{(|k_{[123]}|^2+\sum_{i=1}^2|k_i|^2)(|k_{[5-3]}|^2+|k_3|^2)}dk_{123}\bigg)^2dk_5\\\lesssim&\int \theta(2^{-q}k_{5})\frac{N^{-2+\kappa}}{|k_5|^2}dk_{5}\lesssim\varepsilon^\kappa 2^{2q\kappa}.\endaligned$$
For $I^{43}$ we obtain that
$$\aligned E|\Delta_q I_t^{43}|^2\lesssim&\int \theta(2^{-q}k_1)\frac{1}{|k_1|^2}\bigg(\int\frac{1_{|k_{5-3}|\backsimeq N}}{|k_2|^2|k_3|^2|k_5|^2(|k_{[123]}|^2+|k_2|^2+|k_3|^2)(|k_{[5-3]}|^2+|k_5|^2)}dk_{235}\bigg)^2dk_{1}\\\lesssim&\int \theta(2^{-q}k_{1})\frac{N^{-2+\kappa}}{|k_1|^2}dk_{1}\lesssim\varepsilon^\kappa 2^{2q\kappa}.\endaligned$$

Moreover for $J_t^i$ and other terms in $D^N$ we could use similar calculations and Lemma 6.5 to obtain the same estimates.
Then by using Gaussian hypercontractivity, Lemma 2.1 and Kolomogorov continuity criterion we obtain that $D_N\rightarrow^P0$ as $\varepsilon\rightarrow0$.

\th{Acknowledgement.} We are very grateful to Professor Nicolas Perkowski for telling us the technique of random operators, which helped us to simplify the arguments and improve the results of this paper. We would also like to thank Professor Michael R\"{o}ckner for his encouragement and suggestions for this work.

\end{document}